%% file: nullst.tex
\documentstyle{article}

\input macrosus.tex

\begin{document}

\input{seccion0}

\input{seccion1}

\input{seccion2}

\input{seccion3}

\input{seccion4}

\input{ref}

\input{direc}

\end{document}

%% file: macrosus.tex
\normalbaselines
\def\ifm#1#2{\relax\ifmmode#1\else#2\fi} 

\def\char {{\rm \mbox{char\,}}}

\def\Id {{\rm \mbox{Id}}}
\def\Ch {{{\cal C}{\it h}}}

\def\Conv {{\rm \mbox{Conv}}}
\def\Vol {{\rm \mbox{Vol}}}

\def\Tr {{\rm \mbox{Tr}}}
\def\No {{\rm \mbox{N}}}

\def\Supp {{\rm \mbox{Supp}}}
\def\Hom {{\rm \mbox{Hom}}}
\def\ord {{\rm \mbox{ord}}}
\def\GL {{\rm \mbox{GL}}}

\def\init {{\rm \mbox{init}}}
\def\discr {{\rm \mbox{discr}}}

\def\Gal {{\rm \mbox{\scriptsize Gal}}}

\newcommand{\svsp} {\vspace*{2mm}}

\newcommand{\Strichq}  {\vrule height0.65em width0.05em depth0em \,}
\def \N{{\rm I\kern -2.1pt N\hskip 1pt}} 
\def \Z{{\rm Z\!\!Z}}
\newcommand{\Q}        {\ifm {\mbox{\rm Q}\hspace{-0.54em}\Strichq\>\,}
			     {$\mbox{\rm Q}\hspace{-0.54em}\Strichq\>\,$}}
\def \R{{\rm I\kern -2.2pt R\hskip 1pt}} 
\newcommand \C        {{\ifm {\mbox{\rm C}\hspace{-0.45em}\Strichq\>}
			     {$\mbox{\rm C}\hspace{-0.45em}\Strichq\>$}}}
\def\A{{\it I}\!\!{\rm A}}
\def \P{{\rm I\kern -2.2pt P\hskip 1pt}} 
\def \F{{\rm I\kern -2.2pt F\hskip -1.5pt}} 

\def\cA {{\cal A}}
\def\cB {{\cal B}}
\def\cE {{\cal E}}
\def\cM {{\cal M}}
\def\cN {{\cal N}}
\def\cO {{\cal O}}
\def\cP {{\cal P}}

\def\cV {{\cal V}}
\def\cX {{\cal X}}

\def\Qbarra {\overline{\Q}}

\newtheorem{defn}{Definition}[section]
\newtheorem{lem}[defn]{Lemma}
\newtheorem{mainlem}[defn]{Main Lemma} 
\newtheorem{prop}[defn]{Proposition}
\newtheorem{teo}[defn]{Theorem}
\newtheorem{cor}[defn]{Corollary}

\newtheorem{rem}[defn]{Remark}
\newtheorem{exmpl}[defn]{Example}
\newtheorem{assumption}[defn]{Assumption}

\newenvironment{undef}[1]%
		   {\vspace{3.3mm}
		   \noindent{\bf #1}\it}%
		   {\vspace{3.3mm}}

\newenvironment{proof}[1]{
  \trivlist \item[\hskip \labelsep{\it #1}]}{\hfill\mbox{$\Box$}
  \endtrivlist}

%% file: seccion0.tex
\noindent{\LARGE{ \bf  Sharp estimates for  the arithmetic  
Nullstellensatz}}

\vspace{6mm}

\noindent {\large 
Teresa Krick\footnote{T. Krick and M. Sombra were
partially supported
	   by CONICET, UBACyT and ANPCyT
(Argentina), and by the Mathematical Sciences Research
Institute at Berkeley (USA).

\ \  M. Sombra was also partially supported by grant NSF-DMS 97-29992 to
the Institute for 
Advanced Study at Princeton (USA).}, 
Luis Miguel Pardo\footnote{L. M. Pardo was partially  supported
by PB 96-0671-C02-02 (Spain), and by CNRS 1026 MEDICIS (France).}, 
and Mart\'\i n Sombra\footnotemark[1] }

\vspace{3mm}


\noindent {\small {\bf Abstract.}  
We present  sharp estimates for the degree and the height of 
the polynomials in the Nullstellensatz over $\Z$.
The result improves previous work of Philippon, 
Berenstein-Yger and Krick-Pardo. 

We also present degree and height estimates of intrinsic type, 
which depend  mainly on the degree and the height of the input 
polynomial system. 
As an application, we derive an effective arithmetic Nullstellensatz 
for sparse polynomial systems.

The proof of these results relies heavily on the notion of local height 
of an affine variety defined over a number field. 
We introduce this notion and   study its basic properties. 
}

\vspace{1mm}


\vspace{1mm}



\vspace{-2mm}

\typeout{contenido}

\tableofcontents




\typeout{Introduccion}

\section*{Introduction}

\addcontentsline{toc}{section}{Introduction}

\vspace{2mm}

Hilbert Nullstellensatz is a cornerstone of algebraic geometry. 
Under a simplified form, its statement is 
the following: 
\begin{quote}{\em
Let $f_1,\ldots,f_s \in \Z[x_1,\ldots,x_n]$ be polynomials 
such that the equation system 
\begin{equation} \label{eq0.1}
f_1(x)=0,\ldots, \, f_s(x)=0           
\end{equation}
has no solution in $\C^n$.
Then  there exist $a\in \Z \setminus \{0\}$ and 
$g_1,\ldots,g_s \in \Z[x_1,\ldots,x_n]$ satisfying the B\'ezout 
identity 
\begin{equation} \label{eq0.2}
a = g_1 \,f_1+ \cdots + g_s \,f_s .
\end{equation}
}
\end{quote}
As for many central results in commutative algebra and
algebraic geometry, it is an
existential non-effective statement.
The estimation of both the degree and the height of
polynomials satisfying Identity (\ref{eq0.2}) became an important and 
widely considered  question.
Effective versions of Hilbert Nullstellensatz apply to a wide range of
situations in number theory and theoretical computer
science. In particular, they decide 
 the consistency of a given 
polynomial system. 
In their arithmetic presentation, they apply to 
Lojasiewicz inequalities \cite{Solerno89}, \cite{JiKoSh92}
and to the consistency problem 
over finite fields \cite{Koiran96}, \cite{HaMoPaSo98}.

\smallskip 

We recall that the {\it height} $h(f)$ of a polynomial 
$f\in \Z[x_1,\dots,x_n]$ is defined as the logarithm of the 
maximum modulus of its
coefficients.
The main result of this paper is the following effective arithmetic 
Nullstellensatz: 

\begin{undef}{Theorem 1} 
Let  $f_1,\ldots,f_s \in \Z[x_1,\ldots,x_n]$
be polynomials without 
common zeros in  $\C^n$. 
Set  $d:= \max_i \deg f_i $ and $h:=\max_i h(f_i)$.

Then there exist 
$ a\in \Z \setminus \{ 0\} $ and $g_1,\dots,g_s \in \Z[x_1,\dots,x_n]$
 such that 
\begin{itemize}
\item $a = {g}_1 \, {f}_1 + \cdots + {g}_s \, {f}_s $, 
\item $ \deg g_i \le 4\, n \, d^n  $, 
\item $ h(a) , h(g_i) \leq 
4\, n\,(n+1) \, d^n \, (h + \log s +
(n+7) \, \log (n+1) \, d   )$.
\end{itemize}
\end{undef}

As we will see below, this result substantially improves  all  previously
known estimates for the arithmetic Nullstellensatz. 

\smallskip

The following variant 
of a well-known example due to Masser and Philippon \cite{Brownawell87} 
yields a lower bound for any 
general degree and height estimate.
Set 
$$
f_1:= x_1^d\, ,\, \, f_2:= x_1 \, x_n^{d-1} - x_2^d\, ,\, 
\ldots \, ,\, \, f_{n-1}:= x_{n-2} \, x_n^{d-1}  -
x_{n-1}^{d}\, ,\, 
\, 
f_n:= x_{n-1} \, x_n^{d-1} - H 
$$
for any $n\, ,\, d\, ,\, H \in \N$.
These  are polynomials of degree $d$ and height bounded by $h:= \log H$
without common zeros in $\C^n$.
Let $a \in \Z \setminus  \{0\}$ and 
$g_1, \ldots, g_n \in \Z[x_1, \ldots, x_n]$  such that
$$
a= g_1 \,f_1 + \cdots + g_n\,f_n .
$$
Specializing this identity at
$ x_1:= H^{d^{n-2}} \, t^{d^{n-1}-1}, \dots, x_{n-1}:= H \, t^{d-1}, 
x_n := 1/t $
 we obtain 
$$ 
a = g_1( H^{d^{n-2}} \, t^{d^{n-1}-1}, \dots, 
H \, t^{d-1}, 1/t) \, H^{d^{n-1}} \, t^{d^{n}-d}.
$$
We conclude that $ \ \deg g_1 \ge d^{n}-d \ $ and
$ \  h(a) \ge d^{n-1} \, h$. 
In fact, a modified version of this example gives the improved lower bound
$ \  h(a) \ge d^{n} \, h $ (Example \ref{d^n h}).
This shows that  our estimate is essentially optimal. 

\svsp 

The earlier work on 
the effective Nullstellensatz dealt with the  degree bounds. 
Let $k$ be a field  and let 
$f_1, \ldots, f_s \in k[x_1, \ldots, x_n]$ 
be polynomials of degree bounded by $d$ without common zeros in
$\overline k^n$.

In 1926, Hermann \cite{Hermann26}  (see also \cite{Heintz83}, 
\cite{MaWu83}) proved   that there exist 
$g_1, \ldots, g_s \in k[x_1, \ldots, x_n]$
such that 
$$
1=  g_1 \, f_1 + \cdots + g_s \, f_s
$$
with $ \deg g_i \, f_i \leq  2 \, (2d)^{2^{n-1}}$.

After a conjecture of Keller and Gr\"obner, this estimate was 
dramatically improved by Brownawell \cite{Brownawell87}
to $\deg g_i f_i \leq n^2 d^n+ n\, d $ in case $\char(k)=0$, 
while 
Caniglia, Galligo and Heintz \cite{CaGaHe88} showed that 
$\deg g_if_i \leq d^{n^2}$ holds in the general case. 

These results were  then independently 
refined by Koll\'ar \cite{Kollar88} and by Fitchas 
and Galligo \cite{FiGa90} to 
$$
\deg g_i \, f_i \leq \max \{3,d\}^n, 
$$ 
which  is optimal in case $ d \geq 3$. 
For  $d=2$, Sombra \cite{Sombra98a}  recently showed that the bound
$\ \deg g_i f_i \le 2^{n+1} \ $ holds.

\bigskip

Now, let  us consider the height aspect: 
assume $f_1,\dots, f_s\in\Z[x_1,\dots,x_n]$ are polynomials of
degree and height bounded by $d$ and $h$, respectively. 
The previous degree bound   reduces B\'ezout identity (\ref{eq0.2}) 
to a system of $\Q-$linear equations, which  can be solved by 
Cramer rule. 
The obtained estimate  for the height of the integer $a$
and the polynomials 
$g_i$ is of type $\ s\, d^{n^2}\, (h+ \log s + d)$. 

\smallskip

However, it was soon 
conjectured that the true height bound
should be much smaller. 

Philippon \cite{Philippon90} obtained the following sharper estimate for the 
denominator $a$ in the B\'ezout equation: 
$$
\deg g_i \le (n+2) \, d^n \quad \quad , \quad \quad h(a) 
\le \kappa (n) \, d^n ( h+ d),
$$
where $\kappa (n)$ depends exponentially on $n$.  

The first essential progress on height estimates for all 
the polynomials
$g_i$
was achieved by Berenstein and Yger 
\cite{BeYg91}, who obtained 
$$
\deg g_i \leq n\, (2\, n+1) \, d^{ n} \quad ,  \quad h(a), h(g_i) \leq
 \lambda(n) \, d^{8\, n+3}\, ( h + \log s + d\, \log d) , 
$$
where $\lambda (n)$ is a (non-explicit) constant which depends 
exponentially on $n$. 
Their proof relies on the previous work of Philippon and on 
techniques from complex analysis.

Later on, Krick and Pardo \cite{KrPa94}, \cite{KrPa96} 
obtained 
$$
\deg g_i \leq (n\, d)^{c \, n }  \quad \quad,  \quad \quad  h(a), h(g_i)
\leq (n\, d)^{c \, n} ( h+ \log s + d) , 
$$
where $c$ is a universal constant ($c \leq 35$).
Their proof,
based on duality theory for Gorenstein algebras, is completely algebraic. 

Finally, Berenstein and Yger \cite{BeYg96} improved 
their height bound to 
$\ \lambda (n)\, d^{4 \, n+2}\, ( h + \log s + d\,) $, and extended 
it to the case when $\Z$ is replaced by an arbitrary 
diophantine ring. 
It should be said, however, that the possibility 
of such an extension was already clear  from the arguments of
\cite{KrPa96}.

\smallskip

We refer the reader to the surveys \cite{Teissier91}, \cite{BeSt91}, 
\cite{Pardo95}
 for a broad introduction 
to the history of the effective Nullstellensatz, 
main results and open questions. 
Aside from  degree and height estimates, 
there is a strong current area of research on  
computational issues 
\cite{GiHeSa93}, \cite{FiGiSm95}, \cite{KrPa96},
\cite{GiHeMoMoPa98}, \cite{GiHaHeMoPaMo97}, 
\cite{HaMoPaSo98}. 
There are other results in the recent research papers
\cite{SaSo95}, \cite{Kollar98}, \cite{EiLa99}.

\bigskip

With respect to  previous work, in this paper  we improve in an 
almost optimal way the 
dependence of the height estimate  on $d^n$ and we eliminate the extraneous
exponential constants depending on $n$. 
We remark that the polynomials arising in Theorem 1 
are a slight variant of the polynomials which appear 
in \cite{KrPa96} and can thus  be
 effectively computed by their algorithm. 

\bigskip

Although the  exponential behavior of the degree and height estimates is
--- in the worst-case --- unavoidable, 
it has been observed that there are many particular 
instances in which these estimates can be essentially  improved. 
This has motivated the introduction of 
parameters associated to the input system which identify special families 
whose behavior with respect to our problem
is polynomial instead of exponential.

\smallskip
In this spirit, Giusti et al. \cite{GiHeMoMoPa98} introduced the 
notion of {\em degree of a polynomial system} $f_1, \dots, f_s$.
Roughly speaking, this parameter measures
the degree of the varieties cut  out by 
$f_1, \ldots, f_i$  for $i=1, \dots, s-1$. 
It was soon realized that the degrees in the Nullstellensatz can be 
controlled in terms of this parameter, giving rise to the 
so-called ``intrinsic Nullstellens\"atze''
\cite{GiHeMoMoPa98}, \cite{KrSaSo97},  \cite{GiHaHeMoPaMo97}, 
 \cite{Sombra97}. 

Recently H\"agele, Morais, Pardo and Sombra 
\cite{HaMoPaSo98} (see also \cite{Hagele98}) 
obtained an arithmetic 
analogue of these intrinsic Nullstellens\"atze. 
To this aim, they introduced the notion of {\em height
of a polynomial system}, the
arithmetic analogue of the degree of the 
system. 
They obtained degree and height estimates which depend 
{\em polynomially } on the
number of variables and on the degree, height and complexity
of the input system.
This result followed from their study of the computational complexity 
of the Nullstellensatz. 

\smallskip
In this paper we obtain a dramatical improvement over this result, 
bringing it to an (apparently) almost optimal form.
In particular, we show that the dependence on the degree and the height
of the system is {\em linear}, and 
we eliminate the influence of the complexity of the input. 

\begin{undef}{Theorem 2}
Let $f_1,\ldots,f_s \in \Z[x_1,\ldots,x_n]$
be polynomials without 
common zeros in  $\C^n$. 
Set  $d:= \max_i \deg f_i $ and $h:= \max_i h(f_i)$. 
Let  $\delta$ and $\eta$ denote the degree and the height
of the polynomial system $f_1, \dots, f_s$. 

Then there exist $a\in \Z\setminus \{0\}$ and
$g_1,\ldots,g_s \in \Z[x_1,\ldots,x_n]$  such that 
\begin{itemize}
\item
$a = {g}_1 \, {f}_1 
+ \cdots + {g}_s \,
{f}_s 
$,
\item
$
 \deg g_i  \le  2\,  n^2 \, d\, \delta$,
\item
$h(a), h(g_i) \leq  (n+1)^2 \,  d\, (
2\, \eta +  (h +\log s) \, \delta + 
21\,(n+1)^2\, d\, \log (d+1) \, \delta ).
$
\end{itemize}
\end{undef}

Since  $\ \delta \leq  d^{n-1} \ $ 
and $ \ \eta \le 
n\,d^{n-1} \, 
(h + \log s + 3\,n\, (n+1) \, d) $ 
(Lemma \ref{d^n}) one recovers  from this statement 
essentially the same 
estimates of Theorem 1. 
However, we remark that Theorem 2 is a  more flexible result,  
as there are many situations in which the 
degree and the height of the input system 
are smaller than the B\'ezout bounds. 
When this is the case, it yields a much more accurate estimate
(Subsection \ref{Proof of Theorem 2}).

\bigskip

As an application of Theorem 2 we derive an arithmetic 
effective Nullstellensatz for sparse polynomial systems. 
To state this  result, we first need to introduce some standard notation. 

 The support $\Supp (f_1, \ldots, f_s)$ of a polynomial system 
$f_1, \ldots, f_s \subset \C[x_1,\dots,x_n]$
is defined as the set of exponents of all the non-zero monomials 
of all  $f_i$'s, and 
the {\em Newton 
polytope} $\cN (f_1, \ldots , f_s)$ is  the convex hull of 
this support. 
The {\em (normalized) volume}  of $f_1, \dots, f_s$
equals $n!$ times the volume of the corresponding  Newton polytope.

\smallskip

The notions of Newton polytope and volume of a polynomial system 
give a sharper characterization of its 
monomial structure than  the degree alone.
These concepts were introduced in the context of root counting by Bernstein
\cite{Bernstein75} and Kushnirenko \cite{Kushnirenko76}, and 
are now in the 
basis of sparse elimination theory (see e.g. \cite{Sturmfels93}). 

\smallskip

We obtain the following result:

\begin{undef}{Corollary 3}
Let  $f_1,\ldots,f_s \in \Z[x_1,\ldots,x_n]$
be polynomials without 
common zeros in  $\C^n$. 
Set $d:= \max_i \deg f_i $ and $h:= \max_i h(f_i)$.
Let $\cV$ denote the volume of the polynomial 
system $1, x_1, \dots, x_n, f_1, \dots , f_s$.

Then there exist 
$ a\in \Z \setminus \{0\}$ and 
$g_1,\ldots,g_s \in \Z [x_1,\ldots,x_n]$ 
such that 
\begin{itemize}
\item $a = {g}_1 \, {f}_1 + \cdots + {g}_s \, {f}_s $, 
\item $ \deg g_i \le  2\, n^2 \, d\, \cV$, 
\item $ h(a) , h(g_i) \leq  
2\, (n+1)^3\,  d\, \, \cV\,(\,h
+\log s    +  2^{2n+4} d\,  \log(d+1) )
.$
\end{itemize}
\end{undef}

The crucial observation here is that both the degree and the height 
of a polynomial system are essentially controlled by the
normalized volume. 
This follows from an adequate arithmetic version of 
the Bernstein-Kushnirenko
theorem (Proposition \ref{bernstein}).
Our result follows then from Theorem 2 in 
a straightforward way.

As before, we can apply the worst-case bound 
$ \ \cV \le d^n$ to recover from this result 
an estimate similar to the one presented in  Theorem 1. 
However, this result gives sharper estimates for both the degree and the 
height  when  the input 
system is sparse (Example \ref{ejemplosparse}).

\smallskip

The sparse aspect in the Nullstellensatz was
previously  considered by Canny 
and Emiris \cite{CaEm96}
for the case 
of $n+1$  $n$-variate Laurent polynomials without common roots 
at toric infinity. 
Their result  is the sparse analogue of Macaulay effective 
Nullstellensatz (see e.g. \cite{Lazard77}).
The first general sparse Nullstellensatz was obtained
by Sombra 
\cite{Sombra98a}.
In both cases the authors give bounds for the Newton polytopes of the 
output polynomials in terms of the Newton polytopes of the 
input ones. 
We refer to the original papers for the exact statements. 

It is quite difficult to make a definite comparison between these results 
and ours. The latter 
does not give sharp bounds for Newton polytopes. 
But on the other hand, our 
degree estimate for the general case is  better, 
 while  the height estimate 
is completely new. 

\bigskip

The key  ingredient in  our
treatment of the arithmetic Nullstellensatz is
the notion of  local height  of a variety defined
over a number field $K$. 

\smallskip

Let $V \subset \A^n(\Qbarra)$ be an equidimensional 
affine variety defined over $K$. 
For each absolute value $v$ over
 $K$, we introduce the {\em local height}
$h_v(V)$ of $V$ at $v$  as a Mahler measure 
of a suitable normalized Chow form of $V$. 
This definition is consistent  with the Faltings height $h(V)$ of $V$, 
namely:
$$
h(V) = {1\over[K:\Q]} \, \sum_{v\in M_K} N_v \, h_v(V), 
$$
where $M_K$ denotes the set of canonical absolute values of $K$, and 
$N_v$ the multiplicity of $v$. 

We study the basic properties of this notion. 
In particular we are able to estimate the local 
height of the trace and the norm
of a polynomial $f \in K[x_1, \ldots, x_n]$ with respect to an 
integral extension $K[\A^r] \hookrightarrow K[V]$. 
We also obtain local analogues of many of the global results 
of Bost, Gillet and Soul\'e 
\cite{BoGiSo94} and Philippon \cite{Philippon919495}.

\smallskip

Our proof of the arithmetic Nullstellensatz 
is based on duality theory for Gorenstein 
algebras (Tate trace formula).
This technique was introduced in the context of the 
effective Nullstellensatz  in \cite{GiHeSa93}, 
\cite{FiGiSm95}.
Here, we follow mostly  the  lines of Sabia-Solern\'o \cite{SaSo95}
and Krick--Pardo \cite{KrPa96}.

The trace formula allows   to 
perform division modulo complete intersection ideals, with good
control of the degree and  height of the involved polynomials. 
The local arithmetic intersection theory 
plays, with respect to the
height estimates, 
the role of the classical intersection theory
with respect to the degree bounds.

\bigskip

Finally, we remark that all of our results are valid
not just for $\Q$
but for arbitrary number fields. 
Moreover, 
it  is quite evident from our arguments that they can be extended to 
any product formula field.

In fact,  the general analysis 
over number fields is 
necessary to obtain 
the sharpest estimates for the case $K:= \Q$. 
We also remark that the estimates in the general version of Theorem 1 do not
depend on the involved number field. 

\bigskip

The outline of the paper is the following:

\smallskip 

In Chapter 1, we recall the basic definitions and
properties of the  height of polynomials, and we introduce the notion of
local height of a variety
defined over a number field. 

\smallskip

In Chapter 2, we  derive useful estimates for the local heights 
of the trace and the norm of a polynomial in $K[V]$, and we 
study the
behavior of the local heights of the 
intersection of a variety with a hypersurface. 

\smallskip

In Chapter 3, we recall the basic facts of  duality
theory which will be useful in our context, 
and we prove Theorem 1.

\smallskip

In Chapter 4, we focus on the intrinsic and sparse versions of the 
arithmetic Nullstellensatz.

%% file: seccion1.tex
\typeout{Seccion 1}

\section{Height of polynomials and varieties}
\label{Height of polynomials and varieties}

\setcounter{equation}{0}

\renewcommand{\theequation}{\thesection.\arabic{equation}}

Throughout this paper $\Q$ denotes the field of rational numbers,
$\Z$ the ring of rational integers,
$K$ a number field, and
$\cO_K$ its ring of integers.
We also denote by $\R$ the field of
real numbers, $\C$ the field of
complex numbers,
$k$ an arbitrary field, and
$\overline k$ an algebraic closure of $k$.
As usual, $\A^n$ and $\P^n$ will denote the affine and the projective
space of $n$ dimensions over $\overline{k}$, respectively.

For every rational prime $p$  we denote by $ | \cdot |_p $ the
corresponding $p$-adic absolute value over $\Q$. We also denote
the ordinary absolute value by $|\cdot|_\infty$ or simply by $
| \cdot | $. These form a complete set of independent absolute
values over $\Q$: we identify the set $M_Q$ of these absolute
values with the set $\{ \infty, p \,;\, \, p \ \mbox{prime} \}  $.

For $v\in M_Q $ we denote by $\Q_v$ the completion of
$\Q$ with respect to the absolute value $v$.
In case $v=\infty$ we have
$\Q_\infty=\R$, while in case $p$ is prime, we have that
$\Q_p$ is the $p$-adic field.
There exists a unique extension of $v$ to an  absolute
value over the algebraic closure $\C_v$  of
$\Q_v$, which we denote by $v$.
In case $v=\infty$ we have $\C_\infty = \C$.

\typeout{Sub-Seccion 1.1}

\subsection{Height of polynomials}
\label{Height of polynomials}

In this  section we  introduce the different measures for the
size of a multivariate polynomial, both over
$\C_v$ and over a number field.
We establish the link between the different notions and study their
basic properties.

\typeout{Sub-Sub-Seccion 1.1.1}

\subsubsection{Height of polynomials over
${\mbox{\rm C}\hspace{-0.45em}\Strichq\>}_v$ }

\label{ Height of polynomials over }

We fix an absolute value $v\in \{ \infty, p \,;\, \, p \ \mbox{prime} \}  $
for the rest of this chapter.
Let ${\cal A}\subset \C_v$ be a finite set.
The {\em absolute value} of
$\cal A$ is defined as
$$
|{\cal A}|_v:= \max \{\,|a|_v\, , {a\in {\cal A}} \} ,
$$
and the {\em (logarithmic) height} of  $\cal A$   as
$$
h_v({\cal A}):= \max \{\, 0,\log|{\cal A}|_v \,\}.
$$
For a polynomial
$f = \sum_\alpha a_\alpha \, x^\alpha \in \C_v[x_1,\dots,x_n]$, we define
the {\em absolute value} of
$f$ (denoted by $|f|_v$) as the absolute value of its
set of coefficients, and the
{\em height} of  $f$  (denoted by $h_v(f)$) as the
 height of its set of
coefficients.
That is
\begin{eqnarray*}
|f|_v&:=& \max_\alpha \{\,|a_\alpha|_v\,\}, \\
h_v(f)&:=& \max \{\,0,\log |f|_v\,\} .
\end{eqnarray*}

When $v=\infty$, i.e. when $f$ has complex coefficients,
we shall make use of the
{\em (logarithmic) Mahler measure} of $f$ defined as
$$
m(f):=\int_{0}^1\cdots \int_0^1  \log
|f(e^{2\pi \,i t_1},\dots,e^{2\pi \,i t_n})| \, dt_1\dots dt_n.
$$
This integral is well-defined, as $\log |f|$ is a plurisubharmonic function on
$\C^n$ \cite[Appendix I]{LeGr86}.

The Mahler measure was introduced by Lehmer \cite{Lehmer33}
for the case of a univariate
polynomial $f:=a_d\prod_{i=1}^d (x-\alpha_i)\in \C[x]$ as
$$m(f)=\log |a_d| + \sum_{i=1}^d \max\{ 0  ,\log |\alpha_i|\,\}.$$
The link between both expressions of $m(f)$ is given by Jensen formula.
The general case was introduced and studied by Mahler \cite{Mahler62}.

The key property of the Mahler measure is its additivity:
$$m(f \, g )=m(f)+m(g).$$

We have the following
relation between $\log | f |$ and $m(f)$:
\begin{equation} \label{eq1}
-  \log (n+1)\,\deg f \le m(f) -  \log |f| \le  \log (n+1)\,\deg f .
\end{equation}
The right inequality follows  from the definition of $m$ and
the fact that the number of monomials of $f$ is bounded
by ${n+\deg f\choose{n}} \le (n+1)^{\deg f}$.
For the left inequality, we refer to  \cite[Lemme 1.13]{Philippon86} and 
its proof.

When $f$ has total degree bounded
by $1$, the  inequality is refined to  $\log
|f| \le m(f)$.
Also, for any degree,
 $m(f(x_1, \dots, x_{n-1},0)) \le m(f)$.

\smallskip

We shall make frequent use of the following more
precise relation:

\begin{lem} \label{eq2}
Let $f\in \C[X_1, \ldots, X_r]$ be a polynomial in
$r$ groups of $n_i$ variables each.
Let $d_i$ denote the degree of $f$ in the group of variables $X_i$.
Then
$$
- \sum_{i=1}^r \log(n_i+1) \, d_i \le m(f)  - \log|f| \le
\sum_{i=1}^r \log(n_i+1) \, d_i .
$$

\end{lem}

\begin{proof}{Proof.--}

The right inequality follows directly from the
definition of $m(f)$
and the fact that we can bound by
$\prod_i (n_i+1)^{d_i}$
the number of monomials of $f$.
Thus we only consider the left inequality.

\smallskip

Let $f_{\alpha_1 \cdots \alpha_i}\in \C[X_{i+1}, \ldots, X_r]$ denote
 the coefficient of $f$ with respect to  the monomial $X_1^{\alpha_1}
 \cdots X_i^{\alpha_i}$.
Applying Inequality (\ref{eq1}) we obtain
for all $(\xi_{i+1},\dots,\xi_r) \in \C^{n_{i+1} + \cdots + n_r}$:
$$
\log|f_{\alpha_1 \cdots \alpha_{i-1}}(X_i,\xi_{i+1},\dots,\xi_r)|
\le
m(f_{\alpha_1 \cdots \alpha_{i-1}}(X_i,\xi_{i+1},\dots,\xi_r))
+ \log(n_i+1) \, d_i
.
$$

We have $|f_{\alpha_1 \cdots \alpha_{i-1}}(X_i,\xi_{i+1},\dots,\xi_r)| =
\max_{\alpha_i}  |f_{\alpha_1 \cdots \alpha_{i}}(\xi_{i+1},\dots,\xi_r)|  $.
We integrate both sides of the last inequality on
$S_1^{n_{i+1} + \cdots + n_r}$ and we deduce
$$\max_{\alpha_i} \ m(f_{\alpha_1 \cdots \alpha_{i}})\le
m(f_{\alpha_1 \cdots \alpha_{i-1}}) + \log(n_i+1) \, d_i $$

The statement follows then by induction and the fact that 
$\log |f|_v = \max \{ m(f_{\alpha_1 \cdots \alpha_r});
\ \alpha_i \in \Z^{n_i}\} $.

\end{proof}

Let $f \in \C[X_1,\dots,X_r]$
be a multihomogeneous polynomial in $r$ groups of $n_i+ 1$ each,
and set $f^a$ for a deshomogenization of $f$ with respect to
these groups of variables.
Then $ m(f^a) =m(f)$, $ \log |f^a| = \log |f|$.
Thus
the estimates of the preceding lemma also hold for $f$.

\smallskip

Next we introduce the
{\em (logarithmic) $S_n$-Mahler measure} of a
polynomial $f\in \C[x_1,\dots,x_n]$ as
$$
m(f;S_n):= \int_{S_n} \log |f| \ \mu_n ,
$$
where $S_n:=\{ (z_1,\dots,z_n)\in \C^{n}: |z_1|^2+\cdots+|z_n|^2=1\}$
is the unit sphere in $\C^{n}$, and $\mu_n$ is the measure of total mass
$1$,
invariant with
respect to the unitary group $U(n)$.

More generally, let  $f\in \C[X_1, \ldots, X_{r}]$ be a polynomial
  in $r$ groups
 of $n$ variables each.
Its
{\em $S_{n}^{r}$-Mahler measure} is then defined as
$$
m(f;S_{n}^{r}):= \int_{S_{n}^{r}} \log |f| \ \mu_{n}^{r} ,
$$
with  $S_n^{r}:=S_n\times \cdots \times S_n$.
This alternative Mahler measure was introduced by Philippon
\cite[I]{Philippon919495}.

With this notation, the ordinary Mahler measure $m(f)$
of $f\in \C[x_1,\dots,x_n]$
coincides with $m(f;S_1^n)$.

The $S_n^r$-Mahler measure is
 related to the ordinary Mahler measure by the following inequalities
\cite{Lelong94}:
\begin{equation} \label{emes}
0 \le  m(f) -  m(f;S_n^{r}) \le r \,d\,\sum_{i=1}^{n-1} {1\over 2\,i} ,
\end{equation}
where $d$ is a bound for the degree of $f$ in each group of
variables.

\svsp
Finally,
we summarize in the following lemma
the basic properties of the notion of height of
polynomials in $\C_v[x_1,\dots,x_n]$.

\begin{lem} \label{hprod}
Let $v\in M_Q$ and $f_1,\dots, f_s\in \C_v[x_1,\dots,x_n]$.
\begin{enumerate}
\item If $v=\infty $ then
\begin{enumerate}
\item $h_\infty(\sum_i f_i)\le \max_i\{ h_\infty(f_i)\} + \log s$ .
\item $ h_\infty(\prod_{i=1}^s f_i) \le
\sum_{i=1}^s h_\infty(f_i) + \log(n+1)\sum_{i=1}^{s-1} \deg f_i$ .

$h_\infty(f_1\,f_2)\le h_\infty(f_1) + h_\infty(f_2) +
\log  (n+1)\,\min\{\deg f_1,\deg f_2\}$
\item Let $g\in \C[y_1,\dots,y_s]$. Set  $d:=\max_i\{\deg f_i\}$ and
 $h_\infty:=\max_i\{ h_\infty(f_i)\}$.
Then

\smallskip

$
h_\infty(g(f_1,\dots,f_s))\le h_\infty(g)
+ \deg g \, (h_\infty + \log(s+1) + \log(n+1)\,d).
$
\item $\log |\prod_i f_i |_\infty
\ge \sum_i \log |f_i|_\infty - 2\,\log(n+1)\sum_i \deg f_i   $.
\end{enumerate}
\item If $v= p$ for some prime $p$ then
\begin{enumerate}
\item $h_p(\sum_i f_i)\le \max_i\{ h_p(f_i)\} $ .
\item $ h_p(\prod_if_i) \le
\sum_i h_p(f_i) $ .
\item Let $g\in \C_p[y_1,\dots,y_s]$. Set $d:=\max_i\{\deg f_i\}$ and
 $h_p:=\max_i\{ h_p(f_i)\}$.
Then

$
h_p(g(f_1,\dots,f_s))\le h_p(g)
+ \deg g \, h_p.
$
\item $\log |\prod_i f_i |_p
=  \sum_i \log |f_i|_p   $.
\end{enumerate}
\end{enumerate}
\end{lem}

\begin{proof}{Proof of Lemma \ref{hprod}.--}
(1.a), (1.b), (2.a) and (2.b) are immediate from the definition
of $h_v$.

\smallskip

(1.c) and (2.c):

Let us consider the case $v=\infty$.
Set $c(n):=\log(n+1)$.

First we compute $h_v(f_1^{\alpha_1}\cdots f_s^{\alpha_s})$ for
the exponent
$(\alpha_1,\dots,\alpha_s)$ of a monomial of $g$.
Applying (1.b) we obtain
$$h_\infty(f_1^{\alpha_1}\cdots f_s^{\alpha_s})\le (c(n)\,d  +
 h_\infty ) \, \sum_i {\alpha_i} \le (c(n)\,d  +
 h_\infty ) \, \deg g.$$

The polynomial $g$ has at most $(s+1)^{\deg g}$ monomials
and so
$$
h_\infty(g(f_1,\dots,f_s))\le h_\infty(g) + (c(n)\,d + h_\infty )\deg g   
+ c(s)\deg g
.
$$
The case $v\neq \infty$ follows in a similar way.

\smallskip

(1.d) and (2.d):

In case $v= \infty$, we apply directly  Inequality (\ref{eq1}):
\begin{eqnarray*}
\sum_i \log |f_i|_\infty&\le& \sum_i(m(f_i)+c(n)\deg f_i) \\
&=  & m(\prod_if_i) + c(n)\sum_i \deg f_i \\
&\le& \log |\prod_i f_i|_\infty + 2\,c(n) \sum_i \deg f_i .
\end{eqnarray*}

In case $v=p$, Gauss Lemma implies that
$\sum_i \log |f_i|_p= \log |\prod_i f_i|_p$.
\end{proof}

We shall make frequent use of the following particular case of the previous
lemma:

Let $(f_{i j})_{i j}$ be a $s\times s$-matrix of polynomials in
$\C_v[x_1,\dots,x_n]$ of degrees and heights bounded by $d$ and
$h_v$ respectively.
{}From Lemma \ref{hprod}(a,b)  we obtain:
 \begin{itemize}
 \item
 $h_\infty(\det(f_{i j})_{i j} )\le
s\,(h_\infty + \log s + d\,\log(n+1))$ ,
 \item
 $h_p(\det(f_{i j})_{i j})\le s\,h_p $ .
 \end{itemize}


\typeout{Sub-Sub-Seccion 1.1.2}

\subsubsection{Height of polynomials over a number  field }

\label{ Height of polynomials over a number  field }

The set $M_K$ of absolute values over $K$ which extend
the absolute values in $M_Q$ is called the
{\em canonical set}.
We denote by $M_K^\infty$ the set of archimedean absolute values in $M_K$.

If $v\in M_K$ extends an  absolute value
${v_0}\in M_Q$ (which is denoted
by $v\,|\, {v_0}$)
there exists a (non necessarily unique)  immersion
$
\sigma_v:K \hookrightarrow \C_{v_0} $
 corresponding to $v$, i.e. such that
$|\xi |_v=|\sigma_v(\xi)|_{v_0}$ for every $\xi \in K$.

\smallskip

In the $p$-adic case,
there is a one-to-one correspondence $ \cP \mapsto v(\cP) $ between
the prime ideals of ${\cal O}_K$ which divide $p$,
and absolute values extending $p$, defined by
$$
|\xi|_{v(\cP)} := {p}^{- \mbox{\scriptsize ord}_{\cal P}(\xi)/e_{\cal P}}
= {\No (\cP)}^{- \mbox{\scriptsize ord}_{\cal P}(\xi)/e_{\cal P} \, f_{\cal P}}
$$
for $\xi\in K^*$.
Here $e_{\cal P}$ denotes  the
ramification index of $\cal P$,
 $\ord_{\cal P}(\xi )$   the order of
$\cal P$ in the factorization of $\xi $,
and $\No(\cP)$   the norm of the ideal $\cP$.
Thus $\xi $ lies in $\cO_K$ if and only if
$\log|\xi |_v\le 0 $ for every
non-archimedean absolute value $ v$.

\smallskip

We denote by $K_v$ the completion of $K$ in $ \C_{v_0}$.
The {\em local degree} of $K$ at $v$ is defined as:
$$
N_v:= [K_v:\Q_{v_0}],
$$
and it coincides with the number of different immersions $\sigma:K
\hookrightarrow \C_{v_0}$ which correspond to $v$.

When $v$ is archimedean,  $K_v$ is either
$\R$ or $\C$, and
$N_v$ equals $1$ or $2$ accordingly.
In the non-archimedean case we have
$$
N_v= e_{\cal P}f_{\cal P}
$$
where $f_{\cal P}:=[{\cal O}_K/{\cal P}:\Z/(p)]$
denotes the residual degree
of the prime ideal $\cal P$ which corresponds to $v$.

We have
$$
[K:\Q]=\sum_{v\,|\,{v_0}} N_v
$$
for $v_0\in M_Q$.
The canonical set $M_K$ satisfies the {\em product formula} with multiplicities
$N_v$:
\begin{equation}
\label{product formula}
\prod_{v\in M_K} |\,x\,|_v^{N_v}=1,    \qquad \forall
x\in K^*.
\end{equation}

\smallskip

Let ${\cal A}\subset K$ be a finite set. Let $v\in M_K$ be an
absolute value which extends $v_0\in M_Q$, and let
$\sigma_v$ be an immersion corresponding
to $v$.  The
{\em local absolute value} and the
{\em local height } of ${\cal A}$ {\em at $v$} are
defined as
\begin{eqnarray*}
 |{\cal A}|_v  &:=&  |\sigma_v( {\cal A})|_{v_0} , \\
h_v({\cal A})  &:=&  h_{v_0}(\sigma_v( {\cal A})),
\end{eqnarray*}
respectively.
For a polynomial
$f = \sum_\alpha a_\alpha \, x^\alpha \in K[x_1,\dots,x_n]$, we define
the {\em local absolute value} of
$f$ {\em at} $v$ (denoted by $|f|_v$) as the absolute value at $v$ of its
set of coefficients, and the
{\em local height} of  $f$ {\em at} $v$   (denoted by $h_v(f)$) as the
logarithmic height at $v$ of its set of
coefficients.

\medskip

Finally the {\em (global) height} of  a finite set $\cA\subset K
$ is defined as
$$
h({\cal A}):= {1\over [K:\Q]} \sum_{v\in M_K} N_v \,h_v({\cal A}) .
$$

In the same way,  the
{\em (global) height} of
$f_1,\dots,f_s \in K[x_1,\dots, x_n]$ is defined 
as the global height of its set
of coefficients, that is
\begin{equation} \label{globalheight}
h(f_1,\dots,f_s):= {1 \over [K:\Q]} \sum_{v\in M_K} {N_v}\,\max_i h_v(f_i).
\end{equation}

These quantities do not depend on the field $K$ in which
we consider the set
$\cal A$ or the polynomials $f_1,\dots,f_s$.
This allows us to extend the definition of $h$ to
the polynomial ring $\overline{\Q} [x_1,\dots,x_n]$.

\smallskip

We have $h(\cA)\ge 0$ and $\max_{a\in \cA} h(a) \le h(\cA)$.
In case $\cA \subset {\cO}_K$,
 $h_v (a) =0$ for every  $a\in \cA$ and $v\notin M_K^\infty$, and so
$h(\cA)\le [K:\Q] \max_{a\in \cA} h(a)$.

\medskip

Let $q=m/n\in \Q^*$ be a rational number, where $m\in \Z$ and $n\in \N$ are
coprime.
Then
$h(q) = \max\{ |m|, n\} $, that is, the height of $q$ controls both the
size of the minimal numerator and denominator of $q$.
More generally,  let $\cA \subset \Q$ be a finite set, and
let $a \in \N$ be a minimal common denominator for all
the elements of $\cA$.
Then
 $h(\cA)=\log \max \{\,  |a\,\cA|, a \,\}$.
The following is the analogous statement for the general case:

\begin{lem} \label{hacer entero}
Let $\cA \subset K$ be a finite set. Then there exist
$b\in \Z\setminus\{0\}$ and $\cB \subset \cO_K$ such that
$$
b\, \cA= \cB \quad \quad ,  \quad \quad h(\cA)\le h(\{b\}\cup \cB)
 \le [K:\Q] \, h(\cA).$$
\end{lem}

\begin{proof}{Proof.--}

Let $v$ be a non-archimedean absolute value, and
set $\cP$ for  the corresponding prime ideal of $\cO_K$.
Then $h_v (\cA) =
 c(\cP)\,\log \No(\cP) / e(\cP)\,  f(\cP)$ for some $c(\cP) \ge 0$.
We set
$$
b:= \prod_{\cP} \No(\cP)^{c(\cP)} \quad \quad
,  \quad
\quad
\cB := \{ b\,a : a \in \cA\}.
$$

Clearly $b\in \Z$,
 and $\log|b\,a|_v \le h_v(\cA) -c(\cP)\,\log \No(\cP)    \le 0$
for every
$v\notin M_K^\infty$, that is $\cB\subset \cO_K$.

We have
$$
 h_v (b) = \sum_{\cP} c(\cP) \,\log \No(\cP) =
\sum_{v\notin M_K^\infty} N_v h_v (\cA)
$$
for $v \in M_K^\infty$, and  also
$ \ h_v (\cA) + \log |b|_v \le h_v(\{b\}\cup \cB) \le h_v (\cA) + h_v(b)$
for all $v \in M_K$.
Thus
\begin{eqnarray*}
h(\{b\} \cup \cB )
&\le  &
{1\over [K:\Q]} \sum_{v\in M_K^\infty} N_v \, (h_v(\cA)+ h_v(b)) \\
[-1mm]
&+ & {1\over [K:\Q]} \, \sum_{v\in M_K^\infty} N_v \, h_v(\cA)
+\sum_{v\notin M_K^\infty} N_v h_v (\cA) \\
[1mm]
&\le & [K:\Q] \,h(\cA) .
\end{eqnarray*}

On the other hand
$$
h(\cA)= {1\over [K:\Q]} \sum_{v}N_v
(h_v(\cA) + \log |b|_v) \le {1\over [K:\Q]}\sum_{v}N_v
h_v(\{b\} \cup \cB) = h(\{b\}\cup \cB) .
$$
\end{proof}

Finally, let $\alpha \in \Qbarra^*$ be a non-zero algebraic number, and set
$p_\alpha \in
\Z[t]$ for its primitive minimal polynomial.
We have
$h(\alpha) =
 m(p_\alpha)/\deg \alpha$.
More generally, the height of a finite set
can be seen as the height of the minimal polynomial
of a generic linear combination of its elements.
This gives a partial motivation for the notion of global
height of a finite set.

\begin{lem}
Let $\cA:= \{ a_1, \ldots , a_N \} \subset K$ be a finite set and set 
$$
p_\cA:= \prod_\sigma (u_0 + \sigma(a_1) \, u_1+
\cdots + \sigma(a_N) \, u_N) \in \Q[u_0, \ldots, u_N],
$$
where the product is taken over all $\Q$-immersions
$\sigma:K\hookrightarrow \Qbarra$.
Then
$$
- \log (N+1) \le h(\cA)- {h(p_\cA) /[K : \Q]} \le \log (N+1) .
$$
\end{lem}

\begin{proof}{Proof.--}

Set $ L(u):= u_0 + a_1 u_1 + \cdots + a_N u_N\in K[u]$.
We have $\log|L|_v = h_v(\cA)$ and so
$$
h(p_\cA) \le [K:\Q] \, (h(\cA) + \log (N+1))
$$
by Lemma \ref{hprod}(b).
On the other hand we have
$ \ \log |L|_v \le m(\sigma_v(L)) \ $ for $v \in M_K^\infty$
and thus
$$
\begin{array}{rcl}
[K:\Q] \, h(\cA) & = & \sum_{v\in M_K^\infty } N_v \, h_v(\cA)
+\sum_{v\notin M_K^\infty } N_v \, h_v(\cA) \\[3mm]
&\le & m(p_\cA) +\sum_{v\notin M_K^\infty } N_v \, h_v(\cA) \\[3mm]
& \le &
h(p_\cA) + [K:\Q] \, \log (N+1) .
\end{array}
$$
by application of Lemma \ref{hprod}(d) and Inequality (\ref{eq1}).
\end{proof}


\typeout{Sub-Seccion 1.2}

\subsection{Height of varieties}
\label{Height of varieties}

In this section we introduce the notions of local and global height
of an affine variety defined over a number field.
For  this aim, we  recall the basic facts of
the degree and  Chow form of varieties.
As an important particular case, we study the height 
of an affine toric
variety.

\typeout{Sub-Sub-Seccion 1.2.1}

\subsubsection{Degree of varieties }

\label{ Degree of varieties }

Let $k$ be an arbitrary field and
  $V \subset \A^n$  be an affine equidimensional
variety of dimension $r$.
We recall that the  degree of $V$ is  defined as the number of points
in the intersection
of $V$ with a generic linear space of  dimension $n-r$.
This  coincides with the sum of the degrees of
its irreducible components.

\smallskip

For an arbitrary variety $ V \subset \A^n$ we set
$V = \cup_i \, V_i$ for its decomposition into
equidimensional varieties.
Following Heintz \cite{Heintz83}, we define
the {\em degree} of $V$ as
$$
\deg V:= \sum_i \deg V_i.
$$
For $V=\emptyset$ we agree $\deg V :=1$.

\smallskip

This is a positive integer, and we have
$\deg V=1$ if and only $V$ is a linear variety.

The degree of a hypersurface equals the degree of any generator of its
defining ideal.
The degree of a finite variety equals its cardinal.

For a linear morphism $\varphi : \A^n \to \A^m$ and a variety
$V \subset \A^n$ we have $\deg \overline {\varphi(V)} \leq \deg V$.

\smallskip

The basic aspect of this notion of
degree is its behavior with respect to intersections.
It verifies the {\em B\'ezout inequality}:
$$
\deg (V \cap W) \leq \deg V \, \deg W
$$
for $V,W \subset \A^n$, without any restriction on the intersection type
of $V$ and $W$ \cite{Heintz83}, \cite{Fulton84}.

\typeout{Sub-Sub-Seccion 1.2.2}

\subsubsection{Normalization of Chow forms}

\label{  Normalization of Chow forms}

Let $V\subset \A^n$
be an affine equidimensional  variety of dimension $r$
defined over a field $k$.
Let $f_V$ be a Chow form  of   $ V$, that is a
Chow form of its projective
closure $\overline V \subset \P^n$.
This is  a squarefree
polynomial   over $k$ in $r+1$ groups
$U_0, \ldots, U_{r}$ of $n+1$ variables each.
It is  multihomogeneous
of degree $D:= \deg V$ in each group of variables, and
is  uniquely determined  up to
a scalar factor.
In case $V$ is irreducible, $f_V$ is an irreducible polynomial, and
in the general case of an equidimensional variety,
the product of  Chow forms of its irreducible
components is a Chow form of $V$.

\smallskip
In order to avoid this indeterminacy of $f_V$,
we are going to fix one of its coefficients under some
assumption on the variety $V$.
For purpose of reference,
we resume it in the following:

\begin{assumption} \label{assumption}
We assume that the projection $\pi_V: V\to \A^r$ defined by
 $x \mapsto (x_1,\dots,x_r)$ verifies $\# \pi_V^{-1}(0)= \deg V $.
\end{assumption}

\smallskip

This assumption implies that $\pi_V: V\to \A^r$ is a dominant map of degree
$\deg V$, by the theorem of dimension of fibers.
Later on, we will prove that in fact the projection
$\pi_V$
is finite, that is, the variables
$x_1,\dots, x_r$ are in Noether normal position
 with respect to $V$ (Lemma \ref{evaluacion}).
We remark that the previous condition is satisfied by any variety
under a generic linear change of variables.

\smallskip

Each group of variables $U_i$ is associated to the coefficients of a
generic linear
form $L_i(U_i):= U_{i\,0}+ U_{i \, 1} \, x_1 +\cdots+U_{i\,n}x_n$.
The  main feature of a Chow form is that
$$
f_V(\nu_0,\dots,\nu_{r})=0 \Leftrightarrow  \overline V \cap \, \{
L^h(\nu_0) =0\} \, \cap \dots \cap \, \{L^h(\nu_r) =0\} \ \not= \emptyset
$$
holds for $\nu_i \in \overline{k}^{n+1}$.
Here $L_i^h:= U_{i\,0} \, x_0 +\cdots+U_{i\,n}x_n$ stands
for the homogenization of $L_i$.

\smallskip

Assumption \ref{assumption} implies
that $ \overline{V} \cap \{x_1=0\} \cap \dots \cap \{x_r=0\} $
is a zero-dimensional
variety of $ \P^n$ lying in the affine space $\{x_0\neq 0\}$.
Set $e_i$ for the the
$(i+1)$-vector of the canonical basis of $k^{n+1}$.
Then $f_V(e_0,\dots,e_{r})$ --- that is, the coefficient
of the monomial $U_{0\,0}^D \cdots U_{r\,r}^D  $ --- is non-zero.

\smallskip

We then define {\em the  (normalized) Chow form }
$\Ch_V$  of $V$ by
fixing the election of $f_V$ through the condition
$$
\Ch_V(e_0,\dots,e_r)=1 .
$$

Under this normalization, $\Ch_V$ equals the product
of the normalized Chow forms of the irreducible components of $V$.

\typeout{Sub-Sub-Seccion 1.2.3}

\subsubsection{Height of varieties over
${\mbox{\rm C}\hspace{-0.45em}\Strichq\>}_v$ }

\label{ Height of varieties over }

Let $v$ be an absolute value over
$\Q$, and $V\subset \A^n(\C_v)$ an equidimensional  variety
of dimension $r$
which satisfies Assumption \ref{assumption}.
We introduce the height of $V$ as a Mahler measure
of its normalized Chow form.
This notion is closely related
to Philippon local height of a projective variety
\cite[II]{Philippon919495}.

\medskip

\begin{defn} \label{localh}
The {\em height} of the affine variety $V\subset \A^n(\C_v)$
is defined as
$$h_\infty(V) := m(\Ch_V;S_{n+1}^{r+1}) \,+\,
(r+1)\, (\sum_{i=1}^{n}
 {1 / 2\,i})  \,\deg V  $$
in case $v=\infty$ is archimedean,
and as
$$h_p(V):= h_p(\Ch_V)$$
in case $v$ is non-archimedean, that is $v=p$ for some prime $p$.
\end{defn}

\medskip

Let us consider some examples:
\begin{itemize}

\item
We have that  $ \ h_\infty(\A^n) \ $ equals the Stoll
number $\sum_{i=1}^n\sum_{j=1}^i
{1/2\, j} $,  while $h_p(\A^n)=0$.
This follows from \cite[Lem.3.3.1]{BoGiSo94} and the fact that
$\Ch_{\A^n} = \det (U_0,\dots, U_n)$.

\item
Let $V\subset \A^n(\C_v)$ be a hypersurface verifying Assumption
\ref{assumption}, defined by a squarefree polynomial
 $f\in \C_v[x_1,\dots,x_n]$.
Then the
coefficient of the monomial $x_n^{\deg V}$ is
non-zero, and we can suppose without loss of generality that it equals 1.
Then
\begin{eqnarray*}
h_\infty(V)&=&m(f^h;S_{n+1}) +  (\sum_{i=1}^{n-1}\sum_{j=1}^i
 {1 / 2\,j})\,\deg V , \\
h_p(V)&=&h_p(f),
\end{eqnarray*}
where $f^h$ denotes the homogenization of $f$
\cite[I, Cor. 4]{Philippon919495}.

\item
In case $V=\{\xi\}$ for some $\xi \in \A^n$, we have
(see e.g. \cite[I, Prop. 4]{Philippon919495})
\begin{eqnarray*}
h_\infty (V)&=& {1\over 2} \log (1+|\xi_1|^2+\cdots+|\xi_n|^2), \\
h_p(V)&=&h_p(\xi).
\end{eqnarray*}

\end{itemize}


\typeout{Sub-Sub-Seccion 1.2.4}

\subsubsection{Height of varieties over a number field }

\label{ Height of varieties over a number field}

Let $V\subset \A^n(\Qbarra)$ be
an equidimensional variety of dimension $r$
defined over a number field $K$.

We define the {\em (global) height} $h(V)$ of $V$ as the Faltings
height \cite{Faltings91}
of its projective closure $\overline{V} \subset \P^n$.
Following Philippon \cite[III]{Philippon919495},
we introduce $h$ --- without appealing to
Arakelov theory --- through the
identity
$$
h(V)=
{1\over [K:\Q]} (\,\sum_{v\in M_K^\infty} N_v\,m (\sigma_v(f_V);S_{n
+1}^{r+1}) \ + \
 \sum_{v\notin M_K^\infty} N_v\,\log
 |f_V|_v \, ) \,+\,
(r+1)  \, (\sum_{i=1}^{n}
 {1/2\,i})\,\deg V ,
$$
where $f_V$ denotes any Chow form of $V$
\cite{Soule91}, \cite[I]{Philippon919495}.
This coincides
with the sum of the heights of the irreducible components of $V$.

For an arbitrary affine variety, we define its
{\em (global) height } as the sum of the heights
of its equidimensional components.
We agree that $h(\emptyset):= 0$.

\medskip

We introduce the local counterpart of this notion.
Let $v \in M_K$ be absolute value over $K$,
and suppose that $V$ satisfies Assumption \ref{assumption}.
Let $v_0 \in M_Q $ such that $v|v_0$, and let $\sigma_v:K_v\to \C_{v_0}$
be an immersion corresponding to $v$.
We define the {\em local height} of $V$ at $v$ as
$$
h_v(V):= h_{v_0} (\sigma_v(V)).
$$
This definition
is consistent with the global height,
namely
$$ h(V)=
{1\over [K:\Q]} \,\sum_{v\in M_K} N_v\,h_v(V). $$

This notion  is related to the  height
$\widehat{h}$ of
Bost, Gillet and Soul\'e,  by the formula \cite[Prop. 4.1.2 (i)]{BoGiSo94}:
$$
h(V)= \widehat{h} (V) + (\sum_{i=1}^r\sum_{j=1}^i
{1/ 2\, j}) \, \deg V.
$$

It is  also related to the  height $\widetilde{h}$ introduced
by Giusti et al. \cite{GiHaHeMoPaMo97} in terms of the so-called
geometric solution of a variety.
They are {\em polynomially } equivalent \cite[1.3.4]{Sombra98b}, namely
$$
h(V) \le (n\, \deg V \widetilde{h}(V) )^c \quad \quad ,
\quad
\quad
\widetilde{h}(V) \le (n\, \deg V h(V) )^c,
$$
for some constant $c>0$.

\smallskip

We have $h(V) \ge 0$. Moreover
$h(V) \ge (\sum_{i=1}^r\sum_{j=1}^i
{1/ 2\, j} )\, \deg V$, with
equality only in case $V$ is defined by the
vanishing of $n-r$ standard coordinates.
\cite[Th. 5.2.3]{BoGiSo94}.
For instance $h(\A^n)= \sum_{i=1}^{n}\sum_{j=1}^i {1/2\,j}$.

\medskip

This notion of height satisfies the {\em arithmetic Bezout inequality}
\cite[Th. 5.5.1 (iii)]{BoGiSo94},
\cite[III, Th. 3]{Philippon919495}:
$$
h(V\cap W)\  \le \ h(V)\,\deg W + \deg V\,h(W) + c\,\deg V\,\deg W,
$$
for $V, W \subset \A^n(\Qbarra)$, with
$ \ c:= (\, \sum_{i=0}^{\dim V}\sum_{j=0}^{\dim W}
 {1/ 2(i+j+1)} \, )
+ (\, n-{(\dim V + \dim W)/2} \, ) \log 2$.


\typeout{Subsubsection 1.2.5}

\subsubsection{Height of affine toric varieties}
\label{Height of affine toric varieties}

Now we consider the case of affine toric varieties.
The obtained height estimate is crucial in our treatment of the
sparse arithmetic Nullstellensatz (Corollary \ref{Sparse arithmetic 
Nullstellensatz}).

\smallskip

Let
$\cA=\{ \alpha_1, \dots, \alpha_N \} \subset
\Z^n$ be a finite
set of integer vectors such that $0\in \cA$.
Let $r:= \dim \cA$ denote the {\em dimension } of
$\cA$, that is, the dimension of the free $\Z-$module $\Z  \cA$.
We normalize the volume form of $\R \cA$ in order that any elementary
simplex of the lattice $\Z \cA$ has volume 1.
The {\em (normalized) volume}  $\Vol(\cA) $ of $\cA$ is defined as
the volume of the convex hull
$\Conv(\cA)$ with respect to this volume form.
In case $\Z \cA = \Z^n $, then $\Vol (\cA) $ equals $n!$ times the
volume of $\Conv(\cA)$ with respect to the Euclidean
volume form of
$\R^n$.

\smallskip

We associate to the set $\cA$ a map
$ \, (\Qbarra^*)^n \to \Qbarra^N \ $ defined by
$ \, \xi \mapsto (\xi^{\alpha_1}, \dots, \xi^{\alpha_N})$.
The Zariski closure of the image of this map is
the {\it affine toric variety}  $X_\cA \subset \A^N $.
This is an irreducible variety
of dimension $r$  and degree
$\Vol(\cA)$.

\smallskip

For $i=0, \dots, r$, let $U_i$
denote a group of variables indexed by the elements of $\cA$ and set
$$F_i := \sum_{\alpha \in \cA} U_{i \alpha }  \, x^{\alpha}$$
for the generic Laurent polynomial
with support contained in $\cA$.
Let $  W \subset (\P^{N-1})^{r+1} \times (\Qbarra^*)^n  $
be the incidence variety of $ F_0, \dots, F_r$ in
$ (\Qbarra^*)^n $, that is
$$
W= \{ (\nu_0, \dots, \nu_r; \xi) ; \ \ \ F_i (\nu_i)(\xi)=0 \ \ \forall i\},
$$
and let $ \ \pi : (\P^{N-1})^{r+1}
\times (\Qbarra^*)^n \to (\P^{N-1})^{r+1} \ $
be the canonical projection.
Then
$ \overline{\pi (W)} $
is an irreducible variety of codimension 1.
Its defining polynomial $R_\cA\subset \Q[U_0, \dots, U_r]$
is called the {\it $\cA$-resultant} or {\it sparse resultant}, and
it coincides with the Chow form of the affine
toric variety $X_\cA$ \cite{KaStZe92}.
It is a multihomogeneous polynomial of degree
$\Vol (\cA) $ in each group of variables,
and it is uniquely defined up
to its sign, if we assume it to be a primitive polynomial with integer
coefficients.
Basic references for affine toric varieties and sparse resultants are
\cite{GeKaZe94}, \cite{Sturmfels96}.

\smallskip

We obtain the following bound for the height of
$X_\cA$.
Our argument relies on the
Canny-Emiris determinantal formula for the sparse resultant
\cite{CaEm96}.

\begin{prop} \label{height of toric variety}
Let $\cA \subset \Z^n$ be a finite set of dimension $r$
and cardinality $\# \cA \ge 2$.
Then  $ \  h(X_\cA) \le  2^{2\, r +2} \, \log (\# \cA) \, \Vol (\cA)$.
\end{prop}

\begin{proof}{Proof.--}
Let $\cM$ be the Canny-Emiris matrix associated to
the generic polynomial system $F_0, \dots,
F_r$.
This is a non-singular square matrix of order $M$, where
$M$ denotes the cardinality of the set
$$
\cE:= ((r+1) \, Q+ \varepsilon) \, \cap \, \Z^n .
$$
Here $Q:= \Conv(\cA)$, and  $\varepsilon\in \R^n$ is a
vector such that each point
in $\cE$ is contained in the interior of a cell in a given triangulation
of the polytope $(r+1) \, Q$.

Every non-zero entry of $\cM$ is a variable $U_{i \alpha}$.
In fact, each row has exactly $N$ non-zero entries, which consist
of the variables in some group $U_i$.
We refer to \cite{CaEm96} for the precise construction.

\smallskip

Thus $\det \cM  \in \Z[U_0, \dots, U_r]$ is a multihomogeneous
polynomial of total degree $M$ and height bounded by $M \, \log N$.
This polynomial is a non-zero multiple of the sparse resultant
$R_\cA$  \cite{CaEm96}. 
The assumption that $R_\cA$ is primitive  implies that 
$\det \cM /R_\cA$ lies in $\Z[U_0, \dots, U_r]$, and so 
$ \ m(R_\cA) \le m(\det \cM)$.

\smallskip

Let $\{ T_j\}_{j\in I} $ be a unimodular triangulation of $Q$,
so that
$\{ (r+1) \, T_j \}_{j\in I}$ is a
triangulation of $(r+1) \, Q$.
For every $\varepsilon\in \R^n $, the set of integer points 
in $(r+1) \, T_j + \varepsilon$ is in correspondence with a subset of those
of $ (r+1) \, T_j $.
Moreover, for a generic choice of $\varepsilon$ we loose --- at least ---
the set of integer points in a facet of codimension 1.
Thus 
$$
\# ( (r+1) \, T_j + \varepsilon) \cap \Z^n=  {{2\, r}\choose{r}}
\le 2^{2\, r}
$$
and so $ \ M \le 2^{2\, r}\,  \Vol  (\cA) $.
Applying Lemma \ref{eq2} we obtain 
\begin{eqnarray*}
m(R_\cA) &\le &  \log |\det \cM|
+ \deg (\det \cM)  \, \log N  \\
&\le & 2\, M \, \log N
\\ &\le &  2^{2\, r+1} \, \log N  \, \Vol (\cA).
\end{eqnarray*}

We conclude
\begin{eqnarray*}
h(X_\cA)  & = & m(R_\cA; S_{N+1}^{r+1})
+ (r+1) \,  ( \sum_{i=1}^N 1/2\, i  ) \, \Vol (\cA)\\
&\le & m(R_\cA) + 2\, (r+1) \, \log N \, \Vol (\cA)
\\ [2mm]
&\le &  2^{2\, r+2} \, \log N  \, \Vol (\cA),
\end{eqnarray*}
as $N= \# \cA \ge 2$.
\end{proof}

In case $\cA \subset (\Z_{\ge 0})^n$ --- that is, when $F_0, \dots, F_r$
are polynomials --- we set $ \, d:= \max\{ |\alpha|: \alpha\in \cA\}
= \deg F_0$.
We have then
$ \ N\le (n+1)^d $ and so
$$
h(X_\cA) \le 2^{2\, r + 2} \, \log (n+1) \, d \, \Vol (\cA).
$$

%% file: seccion2.tex
\section{Estimates for local and global heights }

\setcounter{equation}{0}

In this chapter we study the basic properties of
local and global heights that  we will need for our purposes.
The key result is
a precise estimate for the local height
of the trace and the norm of a polynomial $f\in K[x_1, \dots, x_n]$
with respect to
an integral extension $K[\A^r] \hookrightarrow K[V]$.

We also study some of the basic properties of the height of a variety,
in particular its behavior
under intersection with hypersurfaces and under affine maps.


\subsection{Estimates for Chow forms }

In this section we
recall the notion of generalized Chow form of a variety in the sense
of Philippon \cite{Philippon86}, and we prove a technical 
estimate for its local height. 


\subsubsection{Generalized Chow forms}
\label{Generalized Chow forms}

Let $V\subset \A^n$
 be an affine equidimensional  variety of dimension $r$ and
degree $D$ defined over a field $k$.

\smallskip

For $d\in \N$ we denote by  $U(d)_{ 0}$
a group of $d+n \choose {n}$ variables.
Also, for $1\le i \le r$ 
we denote by $U_i$ a group of $n+1$ variables, and we set 
$ \ U(d) := \{ U(d)_0, U_1, \dots, U_r\} $.
Set
$$
F:=\sum_{|\alpha|\le d} U(d)_{0 \alpha} \, x^\alpha \quad \quad \quad ,
\quad 
\quad
\quad
L_i:=U_{i\, 0}+U_{i\, 1}\, x_1+\cdots+U_{i\, n}\, x_n 
$$
for the generic polynomial in $n$ variables of degree $d$ and
$1$ associated to $U(d)_0$ and $U_i$
respectively.

\smallskip

Set $N:= {{d+n}\choose{n}} + r\, (n+1) $ and 
let $W\subset \A^N \times V $ be the incidence variety
of $F, \, L_1, \dots, \, L_r$  with respect to $V$, that is
$$
W:= \{ (\nu(d)_0, \nu_1, \dots, \nu_r; \xi) \ ; \ \
 \xi\in V, \ \ F(\nu(d)_0)(\xi) =0,
\  \ \ L_i(\nu_i)(\xi) =0,
\   1\le i\le r \}  .
$$
Let $\pi : \A^N \times \A^n \to \A^N$ denote the canonical projection.
Then $\overline{\pi (W) } \subset \A^N $ is a hypersurface
\cite[Prop. 1.5]{Philippon86} and any
of its defining equations $f_{d, V} \in k[U(d)]$
is called a {\em generalized Chow form } or
a {\em $d$-Chow form }
of $V$.

\smallskip

A $d$-Chow form  is uniquely defined up to a scalar factor.
It shares many properties with the usual Chow form,
which corresponds to the case $d=1$.
We have
$$
f_{d,V}(\nu(d)_0 , \nu_1, \dots,\nu_r )=0
\Leftrightarrow  \overline V \, \cap \, \{
F^h (\nu(d)_0) =0 \}
\, 
\, \cap \, \{ L_1^h(\nu_1) =0\} \, \cap \cdots \cap \, \{L_r^h(\nu_r) =0\}
\, \not= \emptyset
$$
for $\nu(d)_0 \in \overline{k}^{{d+n}\choose{n}}$
and $\nu_i \in \overline{k}^{n+1}$.
Here $\overline{V} \subset \P^n$ denotes the projective closure
of $V$, while $F^h$ and  $L_i^h$ stand
for the homogenization of $F$ and $L_i$ respectively.

\smallskip

A $d$-Chow form  $f_{d, V}\in k[U(d)]$ is  a multihomogeneous
polynomial
of degree
$D$ in the group of variables
$U(d)_0$ and of degree $d\, D$ in each group $U_i$
\cite[Lem. 1.8]{Philippon86}.
When  $V$ is an irreducible variety, $f_{d, V}$
is an irreducible polynomial of $k[U(d)]$.
When  $V$ is equidimensional, it coincides
with
the product of $d$-Chow forms of its irreducible components.

\medskip

We will appeal to the following relation between a
$d$-Chow form $f_{d,V}$  and the usual one: 

Let $U_0$ be another group of $n+1$ variables, and consider
the morphism 
$$
\varrho_{d}: k[U(d)] \to k[U_0,U_1,\dots,U_r]
$$
defined by
$ \ \varrho_d(F)=L_0^{d} \ $ and $ \ \varrho_d(L_i)=L_i \ $ for
$i =1, \dots , r$, where
 $L_0$ stands for the generic linear form associated to $U_0$.
In other terms
$$\varrho_d(U(d)_{0  \alpha })= {d \choose \alpha}
U_{00}^{d-|\alpha|} \, U_{01}^{\alpha_1} \, \cdots
\, 
U_{0n}^{\alpha_n} \quad \quad \quad \mbox{where}\ \ {d \choose \alpha}
:={ d! \over (d-|\alpha|)! \, \alpha_1! \, \cdots \, \alpha_n!}
$$
for
$|\alpha|\le d$, and $\varrho_d(U_{i\, j})
= U_{i\, j}$ for $i=1, \dots, r$ and $j=0, \dots, n$.
Then $\varrho_d(f_{d,V}) = \lambda\,
f_V^d$
for some $\lambda \in k^*$ \cite[Prop. 1.4]{Philippon86}.

\medskip

Now assume that $V$ satisfies
Assumption \ref{assumption}.
Then
$
 \overline{V} \cap \{x_0^d=0\}
\cap \{x_1=0\} \cap \dots \cap \{x_r=0\} = \emptyset.
$

Setting $e(d)_\alpha $  and
$e_i$ for the $\alpha$-vector and the
$(i+1)$-vector of the canonical bases of $k^{{d+n}\choose{n}}$ and
$k^{n+1}$ respectively, 
 we infer that
  $f_{d,V}(e(d)_0, e_1,\dots,e_{r})$ --- that is, the coefficient
of the monomial $U(d)^D_{00} \, U_{1 1}^{d\, D}
\cdots U_{rr}^{d\, D}  $ --- is non-zero.

\smallskip
We define {\em the (normalized) $d$-Chow form }
$\Ch_{d,V}$ of $V$ by
fixing the election of $f_{d, V}$ with the condition
$\Ch_V(e(d)_0,e_1, \dots,e_r)=1 $.

\smallskip

In the previous construction,
$U(d)_{00}^D \,U_{11}^{d\,D}\dots \,U_{rr}^{d\,D}$
is the only monomial
of $k[U(d)]$
which maps through $\varrho_d$ to
$U_{00}^{d\, D} \cdots U_{rr}^{d\, D}$.
The imposed normalizations
imply  then 
$$
\varrho_d(\Ch_{d,V}) =\Ch_V^d.
$$


\subsubsection{An estimate for generalized Chow forms}
\label{an estimate for gen chow forms}

The following technical result is crucial to our local height estimates for the
trace and the norm of a polynomial (Subsection
\ref{Estimates for norms and traces}), as
well as for the intersection of a variety with an hypersurface
(Subsection \ref{intersection}).
The proof follows   the lines of \cite[Prop. 2.8]{Philippon86}.

\smallskip

We adopt the following convention:

Let $f\in k[x_1, \dots, x_n]$ be a polynomial of degree $d$.
We denote by $f_{d, V}(f)$ and $\Ch_{d,V}(f)$  the specialization of
$U(d)_0$ into the coefficients of $f$ in $f_{d, V}$ and $\Ch_{d,V}$  
respectively.

\begin{lem} \label{import}
Let $V\subset \A^n(\C_v)$ be an equidimensional  variety of dimension $r$
which satisfies Assumption \ref{assumption}.
Let $f \in \C_v[x_1, \ldots, x_n] $.
Then

\begin{itemize}

\item 
$ m(\Ch_{\deg f,V}(f) ;  S_{n+1}^r )
+ r \, (\sum_{i=1}^n 1/2\, i  ) \, \deg f \, \deg V
$\\[2mm]
\hspace*{34mm} $ 
\le \ \deg f \, h_\infty(V)
+  h_\infty(f) \, \deg V +  \log(n+1) \, \deg f \,\deg V
$ \\[0mm]
if $v=\infty$, 

\item $ 
h_p(\Ch_{\deg f ,V}(f) ) \le
\deg f \,h_p(V) + h_p(f) \, \deg V $ \ \ if $v=p$ for some prime $p$.

\end{itemize} 

\end{lem}

We will need the following lemma in order to treat the non-archimedean case:

\begin{lem} \label{abierto}
Let
$f\in \C_p[x_1,\dots,x_n]$, and  let $\Omega\subset \A^n(\C_p)$ be a
Zariski open set. Then
$$
|f|_p=\max\, \{ |f(z )|_p; \ \ z \in \Omega,  \ |z|_p=1 \}.
$$
\end{lem}

\begin{proof}{Proof.--}
For $q \in\N$ we denote by $G_q$ the set of $q$-roots of 1 in
$\Qbarra \subset \C_p$.
Let $\alpha=(\alpha_1, \dots, \alpha_n) \in \Z^n$
such that $ |\alpha_i| < q$. Then
$$\sum_{\omega \in G_q^n} \omega^\alpha
=\left\{ \begin{array}{lcl}
0 & \mbox{\ \ \ if} & \alpha\not= 0, \\
[1mm]
q^n & \mbox{\ \ \ if} & \alpha= 0 . \end{array} \right.
$$

Set $f=\sum_\alpha a_\alpha \, x^\alpha$.
Let  $q > \deg f$ such that $|q|_p=1$, that is 
$ p\hspace{-.7mm}\not \hspace{-.8mm}| \,q$.
Then for any $\nu =(\nu_1, \dots, \nu_n) \in (\C_p^*)^{n}$ we have
$$a_\alpha={1 \over\nu^\alpha \, q^n}
 \sum_{\omega \in G_q^n} f(\nu \, \omega)\, \omega^{-\alpha} .
$$

\smallskip

Let $\alpha \in (\Z_{\ge 0})^n$ such that
$|f|_p=|a_\alpha|_p$.
{}From the previous expression we derive that for each 

$\nu \in S_n:= \{ \nu ;  |\nu_i|_p =1\} $ there exists
$\omega_0\in G_q^n$ such that
$|f|_p= |f(\nu \, \omega_0)|_p$.
The set $S_n $ is Zariski dense
in $\A^n(\C_p)$, and so $S_n \, \cap \, \Omega $ is also dense.
Thus we can take $\nu $ such that $\nu \, G_q \subset \Omega$, and therefore
$$
|f|_p \le \max\{ |f(z )|_p; \ z \in \Omega,  |z|_p=1 \}.
$$

The other inequality is straightforward.
\end{proof}

\begin{proof}{Proof of Lemma \ref{import}.--}
First we consider the case when $V$ is a 0--dimensional
variety.
We may assume
without loss of generality
that $V$ is irreducible, that is $V=\{\xi\} $ for some
$\xi= (\xi_1, \dots, \xi_n) \in \C_v^n$.

Set $d:= \deg f$.
Then 
$$
\Ch_V=L(\xi):=U_0+U_1\xi_1+\cdots +U_n\xi_n \quad \quad
,  \quad
\quad
\Ch_{d,V}=F(\xi):=\sum_\alpha U_\alpha \, \xi^\alpha ,
$$
where $L$ and $F$ denote generic
polynomials in $n$ variables
of degree $1$ and $d$  respectively.
Then
\begin{eqnarray*}
h_\infty(F(\xi)) &=&\log \max_{|\alpha|\le d}
\{ |\xi^\alpha| \}\\
&=&\log \max_i \{ 1, \, |\xi_i|^d \}\\
 &=& d\,h_\infty(L(\xi))\\
[-1mm]
& \le & d\,m(L(\xi); S_{n+1}) + (\sum_{i=1}^n 1/2\, i) \, d
\end{eqnarray*}
and so $ \
h_\infty(\Ch_{d,V}(f) )
\le d\, h_\infty(F(\xi)) + h_\infty(f)+ \log (n+1) \, d 
  \le d\, h_\infty(V) + h_\infty(f)+ \log (n+1) \, d $. 
  
\smallskip

Analogously,
$ h_p(F(\xi)) \le d\,h_p(L(\xi)) $
and so
$h_p(\Ch_{d,V}(f) ) \le d\, h_p(V) + h_p(f)$.

\medskip

Now we consider the general case.
Set $\nu=(\nu_1, \dots, \nu_r) \in \C_v^{r\, (n+1)}$,
$L(\nu_i):= \nu_{i \, 0} + \nu_{i \, 1}\, x_1 + \cdots +
\nu_{i\, n}\, x_n$  and
$$
V(\nu) := V\cap V(L(\nu_1), \dots, L (\nu_r)) \subset
\A^n(\C_v).
$$
Then  $V(\nu)$ is a 0-dimensional variety of
degree $\deg V$ 
for $\nu $ in a 
Zariski
 open set $\Omega$ of $\A^{r\, (n+1)} (\C_v)$. 

\smallskip

Let $\nu\in \Omega$. By  \cite[Prop. 2.4]{Philippon86}
there exist $\lambda (\nu), \theta(\nu) \in k^*$
such that
$$
\Ch_{V(\nu)} = \lambda(\nu)\, \Ch_V(\nu)\quad \quad
, \quad
\quad
\Ch_{d, V(\nu)} = \theta (\nu) \, \Ch_{d, V} (\nu)
,
$$
where $\Ch_V(\nu), \Ch_{d,V}(\nu)$ stand for the
specialization of $U_1,\dots,U_r$ into  $\nu_1,\dots,\nu_r$.
Applying the morphism  $\varrho_d$ linking the
$d$-Chow form with the usual one we obtain 
$$
\Ch_{V(\nu)}^d = \varrho_d(\Ch_{d, V(\nu)})
= \theta (\nu) \, \varrho_d(\Ch_{d, V}  (\nu))
= \theta (\nu) \, \Ch_{V}^d(\nu)
$$
and so $\theta(\nu)= \lambda(\nu)^d$.

\smallskip

We consider the case $v= \infty$. Any Zariski closed set
of $\A^{r\, (n+1)}(\C)$ intersects $S_{n+1}^r$ in a set of
$\mu_{n+1}^r$-measure 0, and so the previous relation
 holds for almost every
$\nu \in S_{n+1}^r$.
Therefore
\begin{eqnarray*}
m(\Ch_{d,V}(f);S_{n+1}^r)& = &
 \int_{S_{n+1}^r} (\log |\Ch_{d,V(\nu)}(f)|
- d\,\log |\lambda(\nu)|)\, \mu_{n+1}^r\\
&\le & \int_{S_{n+1}^r} (d\,h_\infty(V(\nu))
+h_\infty(f) \, \deg V(\nu)+ \log (n+1) \, d\, \deg V(\nu)
- d\,\log |\lambda(\nu)|)\, \mu_{n+1}^r\\
&= &  d\, \int_{S_{n+1}^r}  m( \Ch_{V}(\nu);S_{n+1})
\, \mu_{n+1}^r+
(\sum_{i=1}^n 1/ 2\, i)\, d\, \deg V + h_\infty(f) \, \deg V
\\
[-1mm]
& & + \  \log (n+1)\, d\,\deg V \\
[2mm]
& = & d\,h_\infty(V)
+
h_\infty(f)\,  \deg V +  \log (n+1) \, d\,\deg V  - 
r\,(\sum_{i=1}^n 1/ 2\, i)\, d\, \deg V.
\end{eqnarray*}

The case $v =p$ follows analogously from the 0-dimensional
case and the previous lemma.
\end{proof}

\begin{rem} \label{rem-import}
The only role played by Assumption \ref{assumption} in
the proof of the previous result is in the normalization of the involved 
Chow forms. This is essential in  order to properly define 
local heights of varieties. 
If we disregard normalization, we get altogether the 
following
global result: 

\smallskip
 
Let $V\subset \A^n$ be an equidimensional variety of dimension $r$ defined
over a number field $K$, 
and let 
$f_{d, V}$ be a $d$-Chow form of $V$.
Let $f \in K[x_1 , \dots, x_n] $ be a polynomial of degree $d$.
Then

$
{\displaystyle{1\over [K:\Q]}}({\displaystyle{\sum_{v\in M_K^\infty}}} 
N_v\,m(\sigma_v(f_{d,V}(f)) ;
 S_{n+1}^r )   +  {\displaystyle{\sum_{v\notin M_K^\infty}}} N_v
\log|f_{d,V}(f)|_v) +\    r \, (\sum_{i=1}^n 1/2\, i)  \, d\, \deg V  \  \le  
$\\
\vspace{3mm}
\rightline{$
    d\, h(V)
+   h(f) \, \deg V   \ + \  \log(n+1) \, d\,\deg V. 
$\quad}

\end{rem}


\subsection{Basic properties of the height}

We derive some of the basic properties of the notion of
height of a variety.
In particular, we study the behavior of the height of a variety under
intersection with a hypersurface and under an affine map.

We also obtain 
an arithmetic version of the Bernstein-Kushnirenko theorem.


\subsubsection{Height of varieties under affine maps }
\label{maps}

Let $\varphi: \A^n \to \A^m$ be a regular map defined by
polynomials $\varphi_1, \dots, \varphi_m \in K[x_1,\dots, x_n]$.
We recall that the height
of $\varphi$ is defined as $h(\varphi) := h(\varphi_1, \dots, \varphi_m ) $.

\smallskip

We obtain the following estimate for the height
of the image of a variety under  an affine  map:

\begin{prop}\label{afin}
Let
$V\subset \A^n$ be a variety of dimension
$r$, and  let $\varphi :\A^n\to \A^N$ be an 
affine map.
Then
$$
h({\varphi(V)} ) \le  h(V) +  (r+1)\,
( h(\varphi ) + 8 \, \log (n+N+1)) \, \deg V .
$$
\end{prop}

The proof of this result follows from the study of the particular cases
of a linear projection and an injective affine map.

The
following estimate for the height of a linear projection of a variety
 generalizes \cite[Prop. 2.10]{Faltings91}
and \cite[3.3.2]{BoGiSo94}.
Its proof is essentially based on the description
of the Chow form of such projection variety, due to 
Pedersen and Sturmfels
\cite[Prop. 4.1]{PeSt93}.

\begin{lem} \label{proyeccion}
Let $V\subset \A^n \times \A^m $ be a variety of dimension $r$, and
let
$\pi: \A^n \times \A^m \to \A^n$ denote the projection $(x, y) \mapsto x$.
Then
$$
h(\overline{\pi (V)}) \le h(V) + 3 \, (r+1) \, \log (n+m+1) \, \deg V .
$$
\end{lem}

\begin{proof}{Proof.}
We assume without loss of generality that $V$ is irreducible.
Set $W: = \overline{\pi (V)} \subset \A^n$ and $s:= \dim  W$.

\smallskip

The case $s=r$ follows directly from \cite[Prop. 4.1]{PeSt93}: 
in this case, there exists a partial monomial order $\prec$
such that
$$
f_W \, |\, \init f_V,
$$
where $\init f_V $ denotes the initial polynomial of $f_V$ with respect to
$\prec$.
In particular $\init f_V $ is the sum of some of the terms in the monomial
expansion of $f_V$.

\smallskip

The general case $s\le r$ reduces to the previous one:
we choose standard coordinates
$z_{s+1}, \dots, z_{r} $ of $\A^m$ such that
the projection
$$
\varpi: \A^n\times \A^m \to \A^n\times \A^{r-s}\quad \quad \quad ,
\quad 
\quad \quad (x,y) \mapsto (x , z)
$$
verifies $\dim Z = r$ for 
$Z:= \overline{\varpi(V)} $. 

\smallskip

Let
$\varrho : \A^n\times \A^{r-s} \to \A^n $ denote 
the canonical projection.
Then $f_Z\,|\,\init f_V$, $ \pi = \varrho  \circ \varpi$ 
and  $W = \overline{\varrho (Z)}$.
We have that $\varrho^{-1}(\xi)
=\{\xi\}\times  \A^{r-s}$ for  $\xi \in \varrho (Z)$ by the theorem of 
dimension of fibers. 
Thus $Z =W \times \A^{r-s} $, and in particular
$$
i(W) = Z \cap V(z_{s+1}, \dots, z_{r}) \subset \A^n\times \A^{r-s} ,
$$
where 
$i$ denotes the canonical inclusion $\A^n \hookrightarrow
\A^n\times \A^{r-s}$. 
We have $\deg W = \deg Z$ and so
$ \, f_W:=  f_Z (z_{s+1}, \dots , z_{r})$ is a Chow
form of $W$ \cite[Prop. 2.4]{Philippon86}. 

\medskip

Now we estimate  the height of $f_W$.
Let $K$ be a number field of definition of $V$, and 
set 
$$ 
\init f_V = Q \, f_Z
$$ 
for some polynomial $Q$. 
{}From the proof of \cite[Lem. 1.12(v)]{Philippon86}, there is a
non-zero coefficient $\lambda $ of $Q$
such that
$\log|\lambda|_v \le m(\sigma_v(Q))$
for all $v\in M_K^\infty$. Clearly
 $\log |\lambda|_v \le \log |Q|_v$ also holds  for
all $v\notin M_K^\infty$.
Thus 
$$
m(\sigma_v(f_Z)) \le m(\sigma_v(\init f_V)) - \log |\lambda|_v
$$ 
for $v \in M_K^\infty$, while 
$\log|f_Z|_v \le \log |\init f_V|_v - \log |\lambda|_v$
for $v \notin M_K^\infty$.

\smallskip

Let $v\in M_K^\infty$.
{}From \cite[Lem. 1.13]{Philippon86}
we obtain
$m(\sigma_v(f_W)) \le m(\sigma_v(f_Z))$. Hence
\begin{eqnarray*}
m(\sigma_v(f_W); S_{n+1}^{s+1} ) & \le & m(\sigma_v(f_W)) \\
 &\le & m(\sigma_v(\init f_V))
- \log |\lambda|_v \\  
& \le & \log |\init f_V|_v +  (r+1) \, \log (n+m+1) \,
\deg V - \log |\lambda|_v\\
[2mm]
& \le & \log |f_V|_v +  (r+1) \, \log (n+m+1) \, \deg V - \log |\lambda|_v\\
& \le & m(\sigma_v(f_V ); S_{n+m+1}^{r+1})
+  (r+1) \, (\sum_{i=1}^{n+m} 1/2\, i ) \, \deg V\\
&& + \
 2\, (r+1) \, \log (n+m+1) \, \deg V -  \log |\lambda|_v
\end{eqnarray*}
by application  of Lemma \ref{eq2} and Inequality (\ref{emes}).
In case $v\notin M_K^\infty$ we have analogously
$\log |f_W|_v \le \log |f_V|_v- \log |\lambda|_v$,
and so 
\begin{eqnarray*}
h(W) & \le & h(V) + (s+1) \, (\sum_{i=1}^n 1/2\, i ) \, \deg V+
 2\, (r+1) \, \log (n+m+1) \, \deg V \\
& \le & h(V) + 3\, (r+1) \, \log (n+m+1) \, \deg V .
\end{eqnarray*}
\end{proof}

The following is a variant of \cite[I, Prop. 7]{Philippon919495}:

\begin{lem}\label{inversible} 
Let
$V\subset \A^m$ be a variety of dimension
$r$, and   let $\psi :\A^m\to \A^n$ be an injective 
affine map.
Then
$$
h(\psi(V) ) \le  h(V) + (r+1)\,
( h(\psi ) + 5 \, \log (n+1)) \, \deg V .
$$
\end{lem}

\begin{proof}{Proof.--}
We assume again without loss of generality that $V$ is irreducible.
Let $K$ be a number field of definition of both $V$ and $\psi$, and 
set $\psi (x) = a + A\, x $ for some 
$m\times n-$matrix $A$ of maximal rank and $a\in K^n$.
Then let 
$\ \psi^*: \A^{n+1} \to \A^{m+1} \ $ be the linear map 
$y \mapsto (a,A)^t \, y$ defined by the transpose of the matrix associated to 
$\psi$. 

\smallskip

Set $W:=\psi(V)$, 
and   let $\overline V\subset \P^m$, 
$\overline W\subset \P^n$ denote the projective closures of
$V$ and $W$ respectively. 

\smallskip

For $i=0, \dots, r$ we let $\nu_i \in \Qbarra^{n+1}$, and we set 
$ \ L^h(\nu_i) := \nu_{i\, 0} \, x_0 + \dots + \nu_{i\, n} \, x_n \ $ 
for the homogenization of the associated linear form. 
Then
$f_W(\nu_0,\dots,\nu_r)=0$ if and only if there exists $\xi\in \overline V$
such that $\psi(\xi)$ lies in the linear space determined by $\nu_0, 
\dots, \nu_r$. 
Equivalently $\xi$ lies in the linear space determined by $\psi^*(\nu_0), 
\dots, \psi^*(\nu_r)$. 
We conclude that 
$$
f_W= f_V \circ (\psi^*)^{r+1}.
$$

Let $v\in M_K^\infty$. Then 
\begin{eqnarray*}
m(\sigma_v (f_W), S_{n+1}^{r+1}) & \le & \log |f_W|_v 
+ (r+1)\,\log (n+1)\,\deg V\\
[2mm]
&\le  & \log |f_V|_v + (r+1)\,( h_v(\psi)  + 2\, \log (n+1))\,\deg V
+  (r+1)\,\log (n+1) \, \deg V \\
[2mm]
&\le & m(\sigma_v (f_V)) + \ (r+1)\,\log (m+1) \, \deg V 
+  (r+1)\,( h_v(\psi)  + 3\, \log (n+1))\,\deg V\\
[0mm]
&\le & m(\sigma (f_V), S_{m+1}^{r+1}) + (\sum_{i=1}^m 1 / 2\,i)\,(r+1)\,\deg V\\
&&
+ \ (r+1)\,( h_v(\psi)   +  4\, \log (n+1))\,\deg V. 
\end{eqnarray*}
Here we have applied Lemma \ref{eq2}, Inequality (\ref{emes}) and the proof of
Lemma \ref{hprod}(c), using the fact that the number of monomials of $f_V$
 is bounded by $(n+1)^{(r+1)\,\deg V}$.

\smallskip

In case $v\not\in M_K^\infty$ we obtain analogously 
$ \ \log |f_W|_v \le \log |f_V|_v + (r+1)\, h_v(\psi)\,\deg V $, and hence
$$ 
h(\psi(V) )  \le h(V) + (r+1)\,
( h(\psi ) + 5 \, \log (n+1)) \, \deg V .
$$

\end{proof}

\begin{proof}{Proof of Proposition \ref{afin}.--}

Let $\psi: \A^n \to \A^N\times 
\A^n$ be the injective map 
$x \mapsto (\varphi(x), x)$. 
Then $\varphi$ decomposes as 
$$
\varphi= \pi \circ \psi,
$$
where $\pi: \A^N\times \A^n \to \A^N$ denotes the canonical projection. 
Thus
\begin{eqnarray*}
h({\varphi(V)}) & \le & h(\psi(V)) + 3\, (r+1) \, \log (n+N+1) \, 
\deg \psi(V) \\[2mm]
& \le & 
h(V) + (r+1) \, (h(\psi) + 5\, \log (n+N+1) )\, \deg V + 
3\, (r+1) \, \log (n+N+1) \, 
\deg V \\[2mm]
& = & 
h(V) + (r+1) \, (h(\varphi) + 8\, \log (n+N+1) ) \, \deg V .
\end{eqnarray*}

\end{proof}


\subsubsection{Local height of the intersection of varieties}
\label{intersection}

We obtain the following
estimate for
the local height of the intersection of a variety with an hypersurface. 
This is a consequence of our previous estimate for generalized Chow forms. 
This result can be seen as the local analogue of
\cite[Prop. 2.8]{Philippon86}, and its proof
closely follows it.

\begin{prop} \label{inters}
Let $V\subset \A^n$ be an
equidimensional  variety of dimension $r$
defined over a number field $K$.
Let $f \in K[x_1, \ldots, x_n] $ be a polynomial
which is not a zero-divisor in $K[V]$.
We assume that both $V $ and $V\cap V(f)$ satisfy Assumption
\ref{assumption}.

Then there exists $\lambda\in K^*$ such that

\begin{itemize}\smallskip

\item $
h_v (V\cap V(f)) \le
\deg f \, h_v(V)
+ h_v(f) \, \deg V +
  \log(n+1)  \, \deg f \,\deg V- \log |\lambda|_v
$
\ \ \ for $v\in M_K^\infty$,

\smallskip

\item $h_v (V\cap V(f)) \le
\deg f \, h_v(V)
+ h_v(f) \, \deg V - \log |\lambda|_v $
\ \ \ for $v\notin M_K^\infty$.

\end{itemize}
\end{prop}

\begin{proof}{Proof.--}

Set $d:= \deg f$ and $W:= V\cap V(f) \subset \A^n$.
By \cite[Prop. 2.4]{Philippon86} there
exists  $Q\in
K[U_1,\dots,U_r] \setminus \{ 0\}$ 
such that $\Ch_{d,V}(f)= Q\, \Ch_W$.
Then --- as in the proof of Lemma \ref{proyeccion} --- there exists a
non-zero coefficient $\lambda $ of $Q$
such that
$\log|\lambda|_v \le m(\sigma_v(Q))$
for all $v\in M_K^\infty$ and  
$\log |\lambda|_v \le \log |Q|_v$ for
all $v\notin M_K^\infty$.

\smallskip

Now let $v\in M_K^\infty$. From Inequality (\ref{emes}) we obtain
$$
\log |\lambda|_v \le
m(\sigma_v(Q))\le m(\sigma_v(Q);S_{n+1}^r) + r\,
(\sum_{i=1}^n {1 / 2\,i} )  \, (d\,\deg V-\deg W) 
$$
since $Q$ has  degree
$d\,\deg V-\deg W$ in each group of variables.
Then 
\begin{eqnarray*}
h_v (W)&=& m(\sigma_v(\Ch_W) ;  S_{n+1}^r )
+ r \, (\sum_{i=1}^n 1 / 2\,i ) \,\deg W  \\ 
&= & m(\sigma_v(\Ch_{d,V}(f)) ;  S_{n+1}^r )
+ r\, ( \sum_{i=1}^n {1 / 2\,i} ) \, d \, \deg V \\
[-1mm]
&& - \ m(\sigma_v(Q); S_{n+1}^r) -
r\, (\sum_{i=1}^n {1 /2\,i}) \, (d\,\deg V - \deg W )\\
& \le&
d\, h_v (V)
+ h_v (f) \, \deg V  +  \log(n+1) \, d\,\deg V- \log |\lambda|_v .
\end{eqnarray*}
by straightforward application of
Lemma \ref{import}.
The case $v \notin M_K^\infty$ follows in an analogous way.
\end{proof}

This  result can be immediately generalized to
families of polynomials:

\begin{cor} \label{bezloc}
Let $V \subset \A^n$ be an
equidimensional variety of dimension $r$
defined over $K$.
Let $f_1, \ldots, f_s  \in K[x_1, \ldots, x_n] $ be
polynomials which form a complete intersection in $V$.
We assume that $V\cap V(f_1, \ldots , f_i)$
satisfies
Assumption \ref{assumption} for $i = 0, \dots, s$.
Set $d_i:= \deg f_i$.

Then there exists $\lambda \in K^*$ such that

\begin{itemize} 

\item $
h_v (V\cap V(f_1, \ldots, f_s) ) \le
\prod_i d_i \, ( h_v (V)
+ ( \sum_i h_v(f_i ) / d_i  ) \, \deg V
+  s\,\log(n+1) \, \deg V)
 - \log |\lambda |_v
$

\smallskip 

for $v\in M_K^\infty$,

\item 
$
h_v (V\cap V(f_1, \ldots, f_s) ) \le
\prod_i d_i \, ( h_v (V)
+ ( \sum_i h_v(f_i ) / d_i  ) \, \deg V )
 -  \log |\lambda |_v $
 \  \ for $v\notin M_K^\infty$.

\end{itemize} 
\end{cor}

\begin{proof}{Proof.--}
We just consider the case when $v$ is archimedean, as the other one follows
similarly.
{}From the preceding result we obtain
\begin{eqnarray*}
h_v(V\cap V(f_1,\dots,f_{i}))&\le&
d_{i}\,h_v(V\cap V(f_1,\dots,f_{i-1})) +
h_v(f_i) \, \deg (V\cap V(f_1,\dots,f_{i-1}))\\
[1mm]
& &  +\ \log(n+1) \,  d_i\,\deg
(V\cap V(f_1,\dots,f_{i-1}))- \log |\lambda_i|_v
\end{eqnarray*}
for some $\lambda_i \in K^*$.
For the final estimate we apply iteratively this inequality and we set
$\ \lambda:= \prod_{i=1}^s  \lambda_i^{d_{i+1}\cdots d_s}$.
\end{proof}

\begin{cor}\label{bezloc1}

Let $f_1, \ldots, f_s  \in K[x_1, \ldots, x_n] $ be
polynomials which form a complete intersection in $\A^n$.
We assume that $V(f_1, \ldots , f_i)$
satisfies
Assumption \ref{assumption} for $i=1, \dots , s$.
Set $d_i:= \deg f_i$.

Then there exists $\lambda \in K^*$ such that:

\begin{itemize}

\item $
h_v (V(f_1, \ldots, f_s) ) \le
 \prod_i d_i \, ( \sum_i h_v(f_i ) / d_i  +
 (n+s)\,\log(n+1) ) - \log |\lambda |_v$
\ \  \ for $v\in M_K^\infty$,

\smallskip

\item 
$
h_v (V(f_1, \ldots, f_s) ) \le
 \prod_i d_i \, ( \sum_i h_v(f_i ) / d_i ) - \log |\lambda |_v
$ 
\ \ \ for $v\notin M_K^\infty$.
\end{itemize}
\end{cor}

\begin{proof}{Proof.--}

 We apply  the previous result to $V:= \A^n$, using
 the fact that
$$
h_\infty(\A^n)
= \sum_{i=1}^n \sum_{j=1}^i {1 / 2\, j} \le n \, \log (n+1) \quad
\quad
\quad, \quad
\quad
\quad
h_p(\A^n)=0.
$$
\end{proof}

The following global result is the arithmetic analogue of
\cite[Prop. 2.3]{HeSc82}:

\begin{cor} \label{inters-global}
Let $V \subset \A^n$ be a variety of
dimension $r$, and
let $f_1, \ldots, f_s  \in \Qbarra[x_1, \ldots, x_n] $.
Set $d_i:= \deg f_i$, $h:= h(f_1, \dots, f_s)$ and $n_0 := \min\{ r, s\}$.
We assume that $d_1 \ge  \dots \ge d_s$ holds.
Then:
$$
h (V\cap V(f_1, \ldots, f_s) ) \le
\prod_{i=1}^{n_0} d_i \, \, ( h (V)
+ ( \sum_{i=1}^{n_0} 1/ d_i  ) \, h  \, \deg V + n_0 \,\log(n+1)\, \deg V) .
$$
\end{cor}

\begin{proof}{Proof.--}

We proceed by induction on $(r,s)$ with respect to the product order
of $\N\times \N$. 

\smallskip

The cases $r=0$  or $s=0$ are both trivial.
Thus we assume $r,s \ge 1$. 
Let $V= \cup_C \, C$ be the decomposition
of $V$ into irreducible components.

\smallskip

In case $C \subset V(f_s)  $ we have that 
$ C\cap V(f_1, \dots, f_s) = C\cap V(f_1, \dots, f_{s-1}) $ and so 
$$
h (C\cap  V(f_1, \ldots, f_{s}) ) \le
\prod_{i=1}^{m_0}  d_i \, ( h (C)
+ ( \sum_{i=1}^{m_0}
1/ d_i  ) \, h  \, \deg C + m_0\,\log(n+1)\, \deg C)
$$
with $m_0  := {\min \{r, s-1\}}$.

\smallskip

In case $C \not\subset V(f_s) $ we have either 
$C \cap V(f_s) = \emptyset$ or 
$\dim C \cap V(f_s) \le r-1$. 
The first case is trivial.

\smallskip

In the second case we have 
$$
h(C\cap V(f_s)) \le d_s \, h(C) + h \, \deg C + \log (n+1) \, d \, \deg C .
$$
To obtain this, we proceed as in the proof of Proposition 
\ref{inters}, applying Remark \ref{rem-import} instead of Lemma \ref{import}.

Then 
we apply the inductive hypothesis and we obtain
\begin{eqnarray*}
h (C\cap V(f_1, \ldots, f_{s}) ) &\le&
\prod_{i=1}^{n_0-1}  d_i \, ( h (C\cap V(f_s))
+ ( \sum_{i=1}^{n_0-1}
1/ d_i  ) \, h  \, \deg (C\cap V(f_s))\\
&& + ({n_0-1}) \,\log(n+1)\,
\deg (C\cap V(f_s))) \\
[2mm]
&\le&
\prod_{i=1}^{n_0}  d_i \, ( h (C)
+ ( \sum_{i=1}^{n_0}
1/ d_i  ) \, h  \, \deg C+ {n_0} \,\log(n+1)\,
\deg C) .
\end{eqnarray*}

\end{proof}

With the same notations than Corollary \ref{inters-global},   for  $V:= \A^n$ 
we obtain
$$
h (V(f_1, \ldots, f_s) ) \le
 \prod_{i=1}^{n_0} d_i \, ( ( \sum_{i=1}^{n_0} 1 / d_i  ) \, h +
(n+ {n_0} )\,\log(n+1) ) .
$$


\subsubsection{An arithmetic Bernstein-Kushnirenko theorem}
\label{An arithmetic Bernstein-Kushnirenko theorem}

{}From our estimate for the height of an affine toric variety
(Proposition \ref{height of toric variety}) and the previous results of 
this section
we derive the following arithmetic version of the Bernstein-Kushnirenko
theorem. We refer to Subsection
\ref{Height of affine toric varieties} for the notation.

\begin{prop} \label{bernstein}
Let $f_1, \dots, f_s \in K[x_1, \dots, x_n]$, and
let $\cA \subset (\Z_{\ge 0})^n$ be a finite set such that
$\Supp(1, x_1, \dots, x_n, f_1, \dots, f_s) \subset \cA$.
Set $d:= \max_i \deg f_i$ and $h:=  h(f_1,\dots,f_s)$.
Then

\begin{itemize}

\item $ \deg V(f_1, \dots, f_s) \le \Vol(\cA) $, 

\item $ h(V(f_1, \dots, f_s) ) 
\le (n \, h + 2^{2\, n+3} \, \log (n+1) \, d) \,
\Vol(\cA)$.

\end{itemize}

\end{prop}

\begin{proof}{Proof.--}
Set $\cA:= \{ \alpha_1, \dots , \alpha_N\}$.
The case $N=1$ is trivial, and so we assume $N\ge 2$.
We also assume  that
$\alpha_1, \dots, \alpha_n$ are the vectors of the canonical
 basis of $\R^n$.

\smallskip

The map $\varphi_\cA: \A^n \to \A^N$ induces an isomorphism
between $\A^n$ and the affine toric variety
$X_\cA\subset \A^N$.
The projection map $\pi_\cA: \A^N \to  \A^n$ defined by
$y \mapsto (y_1, \dots, y_n)$  restricted to ${X_\cA}$
 is the inverse map of $\varphi_\cA$.

For $i=1, \dots, s$ we set 
$  f_i= \sum_{j=1}^N a_{i\, j} \, x^{\alpha_j} $ and
we 
let
$$
\ell_i := \sum_{j=1}^N a_{i\,j} \, y_{j} \in K[y_1, \dots, y_N]
$$
be the associated linear form.
Set $V:= V(f_1, \dots, f_s) \subset \A^n$ and
$W:= X_\cA \cap V(\ell_1, \dots, \ell_s) \subset \A^N$.
We have $\varphi_\cA(V) = W$ and so  $V = \pi_\cA(W)$.
Then 
$$
 \deg V \le \deg W \le \deg X_\cA = \Vol(\cA) 
$$ 
and
\begin{eqnarray*}
h(V) & \le & h(W) + 3\, (n+1) \, \log(N+1) \, \deg W \\
[1mm]
& \le & h(X_\cA) + n\, h\, \deg (X_{\cA}) + 4\, (n+1) \, \log(N+1) \,
\deg (X_{\cA})\\
[1mm]
& \le & (n\, h\, + (2^{2\, n +2}\, \log N + 4\, (n+1) \, \log(N+1)) \,
\Vol(\cA).
\end{eqnarray*}
by successive application of Lemma \ref{proyeccion}, 
Corollary \ref{inters-global} and Proposition
\ref{height of toric variety}. 
Finally $N \le {{d+n}\choose{n}}$ and so
$ \ h(V) \le  (n\, h\, + 2^{2\, n +3}\, \log(n+1) \,d \,  ) \, \Vol(\cA)$.
\end{proof}

It seems that the factor $2^{2\, n}$
in the estimate of $h(X_A)$ is superfluous.
If this is the case, the above estimate
can be considerably improved.
We remark that Maillot 
has   recently obtained a much more precise
estimate for the height of the {\em isolated } points of
$V(f_1, \dots, f_s)$ \cite[Th. 8.48]{Maillot99}
.


\subsection{Local height of norms and traces}
\label{norms and traces}

Let $V\subset \A^n $ be
an equidimensional variety of dimension $r$ and degree $D$
defined over a field
 $k$ which satisfies Assumption \ref{assumption}.
As we will see below, this implies that
the projection $ \pi_V: V \to \A^r$ defined by $x\mapsto (x_1, \dots, x_r)$ 
is finite 
(Lemma \ref{evaluacion}).
Set $L:= k(\A^r) $ and $M := L \otimes_{k[\A^r]} k[V]$, so that 
$M$ is a finite $L$-algebra of dimension $D$.

\smallskip

Let $f\in k[x_1,\dots,x_n]$.
We identify $\overline{f}\in k[V]$
with the multiplication map $M\to M$ defined by
$q\mapsto f\, q$.
The {\em characteristic polynomial}
${\cal X}_f \in L[t]$ of $f$ is then defined as
the characteristic polynomial
of this $L$-linear map.
The fact that 
the inclusion $\pi_V^*: k[\A^r] \hookrightarrow k[V]$ is integral implies that 
this polynomial lies in $k[\A^r][t]$, and we have
$\cX_f (\overline f)=0$ in $k[V]$.

\smallskip

Set $\cX_f = t^D + b_{D-1} \, t^{D-1} + \cdots + b_0 \in k[\A^r] [t]$.
Then the {\em norm} $\No_V(f) $ and the {\em trace} $\Tr_V (f) $ of $f$
are defined as
$$
\No_V(f):= (-1)^D \, b_0 \in k[\A^r] \quad
\quad, \quad
\quad
  \Tr_V (f) := -b_{D-1} \in k[\A^r] .
$$
They equal the determinant and the trace of the $L$-linear map
$\overline{f}: M\to M$ respectively.
We also define the {\em adjoint polynomial} $f^*$ of $f$ as
$$
f^*:= (-1)^{D-1} \, (f^{D-1}+ b_{D-1}\, f^{D-2}+ \cdots + b_1)
\in k[x_1, \dots, x_n].
$$
We have $\overline{f^*} \overline{f} = \No_V(f)$ in $k[V]$.

\medskip

The key result of this subsection is a precise  bound for the
 height of the norm and the trace of a
polynomial in case $k$ is a number field.


\subsubsection{Characteristic polynomials}
\label{Characteristic polynomials}

Let $V\subset \A^n$ be an equidimensional
 variety of
dimension $r$ and degree $D$ defined over $k$.
We keep notations as in Subsection
\ref{Generalized Chow forms}:
 for $d\in \N$ we denote by
$ \ F:=\sum_{|\alpha|\le d} U(d)_{0 \, \alpha} \, x^\alpha \ $
and $ \ L_i:=U_{i\, 0}+U_{i\, 1}\, x_1+\cdots+U_{i\, n}\, x_n \ $
the generic polynomial of degree $d$ and 1 
associated to the group of variables
$U(d)_0$ and $U_i$
respectively.

\smallskip

As before, we set $U(d):=\{ U(d)_0,U_1,\dots ,U_r\}$ and
$N:= {{d+n}\choose{n}} + r\, (n+1) $.
Also we introduce an additional  group $T:=\{T_0,\dots,T_r\}$
 of $r+1$ variables which correspond to the
coordinate functions of $\A^{r+1}$.
We consider the map
$$
\psi: \A^N\times \A^{n} \to \A^N\times \A^{r+1} \quad, \quad
(U(d),x) \mapsto (\nu(d), \ F(\nu(d)_0)(\xi), \ L_1(\nu_1)( \xi), \dots, \
 L_r(\nu_r)( \xi)).
$$
where $\nu(d):=(\nu(d)_0,\nu_1,\dots,\nu_r) \in \A^N$ and $\xi \in \A^n$.

Then the Zariski closure 
$\overline{\psi(\A^N\times V)} \subset \A^N\times \A^{r+1}$
is a hypersurface, and any of its defining equations
$P_{d, V} \in k[U(d)][T]$ is called a {\em $d$-characteristic
polynomial } of $V$.
Also we define the {\em characteristic polynomial} of $V$
by $P_{V}:= P_{1, V}$.

\smallskip

A $d$-characteristic polynomial is uniquely defined up
to a scalar factor.
In case $V$ is an irreducible variety, $P_{d, V}$
is an irreducible polynomial. 
When $V$ is  equidimensional, it coincides
 with
the product of $d$-characteristic polynomials of its irreducible components.

\medskip

The following construction links
the  characteristic polynomial of a variety with 
its generalized Chow form.
Set
$$
\zeta(d)_{0 \alpha}:= 
\left\{ \begin{array}{ll}
U(d)_{00} - T_0 \quad & \quad
\mbox{ for } \alpha = 0 \\[1mm]
U(d)_{0 \alpha} 
 \quad & \quad
\mbox{ for } \alpha \ne 0 . 
\end{array}
\right.
$$
Analogously, for $i=1, \dots, r$  we set 
$ \zeta_{i \, 0}:= U_{i \, 0} - T_i $ and 
$ \zeta_{i \, j}:= U_{i \, j} $ for $ j \ne 0$. 
Finally we set
$\zeta(d):= (\zeta(d)_0, \zeta_1, \dots, \zeta_r)$.

\begin{lem} \label{caracteristico}
Let $V\subset \A^n$ be an equidimensional
 variety of
dimension $r$ and degree $D$. Let $f_{d,V}$ be a
$d$-Chow form of $V$.
Then $f_{d,V}\circ \zeta(d)$ is a $d$-characteristic polynomial
of $V$.
\end{lem}

\begin{proof}{Proof.--}

It is enough to consider the case when $V$ is
irreducible.

Let $P_{d,V}$ be a $d$-characteristic polynomial of $V$.
For $(\nu(d) , \xi) \in \A^N\times V$ we set 
$$ 
\vartheta:= 
( F(\nu(d)_0)(\xi), \, L_1(\nu_1)(\xi), \dots, \, L_r(\nu_r)(\xi)) 
\ \in \A^{r+1}, 
$$
so that 
$P_{d, V} (\nu(d))( \vartheta)=0$. 
We observe that 
$$
\xi \in  V\, \cap \, \{F(\nu(d)_0 )( x) =\vartheta_0\} \, \cap
\{ L_1(\nu_1 )( x) =\vartheta_1\} \, \cap \cdots \cap
\, 
\{ L_r(\nu_r )( x) = \vartheta_r)\} \ \subset \A^n.
$$
In particular, this variety is non-empty, and so we infer that
$f_{d,V}\circ \zeta(d)(\nu(d), \vartheta)=0 $.
This implies that $P_{d, V} | f_{d,V}\circ \zeta(d)$ as
$P_{d, V} $ is an irreducible polynomial.

\smallskip

On the other hand 
$f_{d,V}\circ \zeta(d)$ is also irreducible, as it is
multihomogeneous and
$f_{d,V}\circ \zeta(d)(U(d), 0)= f_{d,V}(U(d))$.
We conclude that $P_{d,V}$ and $
f_{d,V}\circ \zeta(d)$ coincide up to a
factor in $ k^*$.
\end{proof}

The previous construction shows that a
$d$-characteristic polynomial of $V$ is multihomogeneous
of degree $D$ in
the group of variables $U(d)_0 \cup \{T_0\}$
and of degree $d\, D$ in each group
$U_i \cup \{T_i\}$.

\medskip

Set $k_d:= \overline{k(U(d))}$, and set
$$
\phi: \A^n(k_d) \to \A^{r+1}(k_d)\quad \quad
\quad
,  \quad
\quad
\quad
x\mapsto (F(x), L_1(x), \dots, L_r(x)).
$$
Then $P_{d, V}\in k_d[T]$
is also a minimal equation for the hypersurface
$\overline{\phi(V)}$, and
by B\'ezout inequality we have also 
$\deg_T P_{d, V} \le d\, D$ (see e.g. \cite{SaSo95}).

\medskip

We assume  from now on that $V$ satisfies Assumption \ref{assumption},
that is  $\# \pi_V^{-1}(0)= \deg V$.
In order to avoid the indeterminacy of the
$d$-characteristic polynomial, we fix it as
$$
P_{d, V}:= (-1)^{D} \, \Ch_{d,V}\circ \zeta(d) .
$$
In particular, we set 
$P_V:=(-1)^D\,\Ch_V\circ \zeta(1)$
for 
{\em the} characteristic polynomial of $V$. 

\smallskip

Set $P_{ V}:= a_D\, T_0^D+ \cdots +a_0 $ for the expansion of $ P_V $
with respect to $T_0$. 
We have that $P_V$ is multihomogeneous of degree 
$D$ in each group $U_i\cup\{ T_i\}$. This implies that 
$a_D$ lies in fact in
$k[U_1, \dots, U_r]$ 
and  is multihomogeneous
of degree $ D$ in each $U_i$ for $i=1, \dots, r$. 

Moreover, $a_D$ coincides  with the coefficient of
$U_{00}^D$ in $\Ch_V$, and
the imposed normalization on $\Ch_{V}$ implies
that 
$$ 
a_D(e_1, \dots, e_r) = \Ch_V(e_0, e_1, \dots, e_r)=1.
$$

We  extend the morphism $ \varrho_{d}$ of 
Subsection \ref{Generalized Chow forms} to
a morphism 
$k[U(d)][T] \to k[U_0,\dots,U_r][T]$ defining
$\varrho_d(U(d)_{00} - T_0): =(U_{00}-T_0)^d$
and $\varrho_d(T_i) := T_i$. 
In other terms
$$
 \varrho_{d}(T_0)= \sum_{j=1}^d (-1)^{j-1} \, {d \choose j} \, 
 U_{0\, 0}^{d-j}\, T_0^j  . 
$$
{}From the previous lemma we obtain 
$$
\varrho_{d}(P_{d , V}) =\varrho_{d}((-1)^D\Ch_{d , V}\circ \zeta(d))
= (-1)^D(\Ch_V\circ\zeta(1))^d = (-1)^{(d+1)D}P_V^d. 
$$
Now set
$$
P_{d, V} = a_{d ,  D} \, T_0^D + \cdots + a_{d, 0}
$$
for the expansion  of $P_{d, V} $ with respect to $T_0$. 
The previous remark implies 
that 
$a_{d,D}=\varrho_d(a_{d, D})= a_D^d$.
In particular $a_{d,D}\in k[U_1,\dots,U_r]$ and
$a_{d,D}(e_1,\dots,e_r)= 1$.

\medskip

The following lemma allows us to obtain a
characteristic polynomial of $f\in k[x_1, \dots, x_n]$
from the $d$-characteristic polynomial of the variety $V$.

\smallskip

We introduce the following convention: 

Given a polynomial $f\in k[x_1, \dots, x_n]$ of degree $d$ and linear 
forms 
$\ell_1, \dots, \ell_r \in k[x_1, \dots,x_n]$, we denote by 
$P_{d, V}(f, \ell_1, \dots, \ell_r) $ the specialization of the variables 
in $U(d)$ into the coefficients of $f, \ell_1, \dots, \ell_r$.

\begin{lem} \label{evaluacion}
Let $V\subset \A^n$ be an equidimensional
 variety of
dimension $r$ and degree $D$ which satisfies Assumption
\ref{assumption}.
Then the projection $\pi_V:V\to \A^r$ is finite.

Moreover, for a polynomial
$f\in k[x_1, \dots, x_n]$  of degree
$d$, the characteristic polynomial of $f$ is
given by
$$
\cX_f = P_{d, V}
(f, e_1, \dots, e_r)(t, x_1, \dots, x_r) \ \in k[\A^r][t]. 
$$

\end{lem}

\begin{proof}{Proof.--}

We have that $P_V(U_0,\dots,U_r)(L_0,\dots,L_r)=0$ in
$k[U]\otimes k[V]$ and so 
$$
P_V(e_j,e_1,\dots,e_r)(t,x_1,\dots,x_r)\in k[\A^r][t]
$$
is a monic equation for $x_j$ in $k[V]$, for 
$j=r+1, \dots, n$. 
Thus the projection $\pi_V$ is finite.

\smallskip

For the second assertion,  set
$$
P_F(t) := P_{d, V}(U(d)_0, e_1, \dots, e_r)(t, x_1,
\dots, x_r)
\ \in k[U(d)_0][\A^r][t] .
$$
This is a polynomial of degree $D$. It is monic 
with respect to $t$, as $a_{d ,
D}\in k[U_1,\dots,U_r]$ and $a_{d,D}(e_1, \dots, e_r)=1$.
We have $P_F(F)=0$ in $k[U(d)_0]\otimes k[V]$.

\smallskip

Now let $m_F$ 
be the monic minimal polynomial of $F$. 
Let $U^\prime(d)_0$ be a  group of ${d+n-r}\choose{n-r}$ variables 
and set $F_0$ for the 
generic polynomial of degree
$d$ in the variables $x_{r+1},\dots,x_n$. 

Then 
$$
m_F(U^\prime(d)_0, 0) \in k[U^\prime(d)_0][t]
$$ 
is an equation for $F_0$
over $\pi_V^{-1}(0)$. 
Since
$\pi_V^{-1}(0)$ is a 0-dimensional variety of degree $D$ and 
$F_0$ separates its points, we infer  that $\deg_{T_0} m_F=D$, 
and so $P_F=m_F$. 

Finally we obtain
$$
\cX_f = \cX_F(f)= P_F(f)=
P_{d, V}(f, e_1, \dots, e_r)(t, x_1, \dots, x_r) .
$$
\end{proof}


\subsubsection{Estimates for norms and traces}
\label{Estimates for norms and traces}

Finally we prove the announced
estimates for the height of the norm and the trace of a polynomial.

\begin{lem} \label{norma}
Let  $V\subset  \A^n$ be an
equidimensional variety of dimension $r$
defined over $K$ which satisfies
Assumption \ref{assumption}.
Let $f\in K[x_1,\dots,x_n]$.
Then

\begin{itemize}
\item $\deg \No_V (f)\le \deg f \,\deg V $,

\smallskip

\item $
h_v(\No_V (f)) \le \deg f \,h_v(V) +
h_v(f) \, \deg V +
 (r+1)\,\log(n+1) \,\deg f \,\deg V
$ \ \ \ for $v \in M_K^\infty$,

\smallskip

\item $
h_v(\No_V (f)) \le \deg f \,h_v(V) +
h_v(f) \, \deg V $
\ \ \ for $v \notin M_K^\infty$.
\end{itemize}
\end{lem}

\begin{proof}{Proof.--}
We keep notations as in Subsection \ref{Characteristic polynomials}.
Set $d:= \deg f$ and $D:=\deg V$.
We have then
$$
\No(f)=
(-1)^D\, P_{d, V} (f, e_1, \dots, e_r)(0, x_1, \dots, x_r)
=  \Ch_{d,V} (f,e_1- e_0 \, x_1,\dots,e_r- e_0 \, x_r)
$$
by Lemmas \ref{evaluacion} and \ref{caracteristico}.
Then 
$$
\deg \No (f) \le \deg_T P_{d, V} \le d\,D.
$$

{}From the previous expression we also obtain that the coefficients
of $\No(f)$ are some of the coefficients of $\Ch_{d, V}(f)$, and so
$|\No(f)|_v\le |\Ch_{d,V}(f)|_v$ for every absolute
value $v$ of $K$.

\smallskip

Let $v\in M_K^\infty$.
Then 
\begin{eqnarray*}
 \log|\No(f)|_v &\le & \log|\Ch_{d,V}(f))|_v \\
& \le & m(\sigma_v(\Ch_{d,V}(f)); S_{n+1}^r) + r\,
( \sum_{i=1}^n {1/ 2\, i} )\, d\, D +
r\, \log (n+1) \, d\, D \\
 & \le & d\,h_v(V) +h_v(f) \, D + (r+1)\,\log(n+1) d\, D
 \end{eqnarray*}
by Inequalities  (\ref{eq2}) and 
(\ref{emes}),
and Lemma \ref{import}.
In a similar way
we obtain
$h_v(\No(f)) \le d\,h_v(V) +h_v(f) \, D$ for $v \notin M_K^\infty$.
\end{proof}

The proof of the following lemma follows closely that of
\cite[Lem. 9]{SaSo95}.
We slightly improve the degree estimate obtained therein,
and we get the corresponding height estimate.

\begin{lem} \label{traza}
Let
 $V\in \A^n$ be an
equidimensional variety of dimension $r$ defined over
$K$ which satisfies Assumption \ref{assumption}.
Let $f,g \in K[x_1,\dots,x_n]$ such that
$\overline{f}$ is not a zero-divisor in $K[V]$.
Set $d:=\max \{ \deg f,\deg g \}$ and
$h_v:=\max \{ h_v(f),h_v(g)\}$
for $v\in M_K$.
Then

\begin{itemize}

\item $ \deg \Tr_V(f^*g) \le d\,\deg V $,

\smallskip

\item $
h_v(\Tr_V (f^*g)) \le d\,h_v(V) +
( h_v + \log 2) \, \deg V +
 (r+1)\, \log(n+1) \,d\,\deg V $
\ \ \ for $v \in M_K^\infty$,

\smallskip

\item $
h_v(\Tr_V (f^*g)) \le d\,h_v(V) +
h_v \, \deg V $
\ \ \ for $v \notin M_K^\infty$.

\end{itemize}

\end{lem}

\begin{proof}{Proof.--}
Let $D:=\deg V$, and let $t$ be a new variable.
Then $K[x_1, \dots, x_r , t] \hookrightarrow K[V\times \A^1]$
is again an integral inclusion.
Set $Q(t):= \No_{V\times \A^1} (t \, f - g)\in K[x_1, \dots, x_r , t] $.
We have then $\No(t-f^*g)={\cal X}_{f^*g}(t)$ and so
$$
\No_V (f)^{D-1} \, Q = \No_V (f^*) \, Q =
{\cal X}_{f^*g}(\No_V (f) t) .
$$
Set $Q= c_D \, t^D+ \cdots + c_0 $ with $c_i \in K[\A^r]$.
The last identity implies then
 $\ \Tr (f^*g)= - c_{D-1}$.

\smallskip

Set $q >  D$, and let $G_q$ denote the group of
$q$-roots of 1.
Then $ Q(\omega) = N_V (\omega \, f -g)$ for $\omega \in G_q$,
and so
$$
\Tr (f^*g)=- {1 \over q} \sum_{\omega \in G_q}
\No_V(\omega \, f- g) \, \omega^{1-D} .
 $$
{}From Lemma \ref{norma} we get 
$\deg \Tr (f^*g)\le d\,D $.

For $v\in M_K^\infty$, we  then obtain
$$
h_v (\Tr (f^* g)) \le \max_{\omega\in G_q}
h_v(\No_V(\omega\, f -g))\le
 d\,h_v(V) + (h_v + \log 2) \,D +
 (r+1)\, \log(n+1) \, d\, D.
$$
Analogously, for $v\notin M_K^\infty$
we take $q >  D$ such that $|q|_v=1$, and we obtain 
$h_v(\Tr (f^* g))\le  d\,h_v(V) + h_v \, D$.
\end{proof}

%% file: seccion3.tex
\typeout{Seccion 3}
  
\section{An effective arithmetic Nullstellensatz}

\setcounter{equation}{0}

In this chapter we obtain the announced estimates for the arithmetic 
Nullstellensatz over 
the ring of integers of a number field $K$. 
Theorem 1 of Introduction corresponds to  the case
$K:= \Q$. 

These estimates depend on the number of variables and on 
the degree 
and height of
the input polynomials.

\subsection{Division modulo complete intersection ideals}
\label{Division modulo complete intersection ideals}

The tool we will use here is Tate trace formula, which
 has already been  used in several papers on elimination
theory. 
One of its outstanding features is that 
it performs effective division modulo 
complete intersection ideals 
\cite{GiHeSa93}, 
\cite{FiGiSm95}, \cite{KrPa96}, \cite{SaSo95}, \cite{GiHaHeMoPaMo97}, 
\cite{HaMoPaSo98}. 
In this section we apply  trace formula to obtain
sharp height estimates in the division procedure.


\subsubsection{Tate trace formula}
\label{Tate trace formula}

We describe in what follows the basic aspects of duality theory for 
complete intersection algebras that we will need in the sequel. 
We refer to 
Kunz \cite[Appendix F]{Kunz86} for a complete presentation of this theory. 

\smallskip

Let $k$ be a perfect field, and set 
$A:= k[t_1,\ldots,t_r]$ and $A[x]:= A[x_1,\ldots,x_n]$.
Let $F:=\{F_1,\ldots,F_n\} \subset A[x]$  be a reduced complete 
intersection which defines a radical ideal $(F)$ of dimension $r$. 

We consider the $A$-algebra 

$$
B:= A[x]/(F)= A[x_1,\dots,x_n]/(F_1,\dots,F_n).
$$

We assume that the inclusion $A\hookrightarrow B$ is finite,  
that is the variables $t_1,\dots, t_r$ are in
Noether normal position with respect to the variety 
$V:=V(F) \subset \A^{r+n}$.  
This is the case, for instance, if $V$ satisfies Assumption 
\ref{assumption}.
Thus $B$ is a
projective $A$-module, which turns to be free  of rank bounded by $\deg V$
by  Quillen-Suslin theorem.

The dual $A$-module $B^*:=\Hom_A(B,A)$ can be seen as a $B$-module with 
scalar
multiplication defined by $f \cdot \tau (g):=\tau(f\,g)$ for $f,g\in B$ and
$\tau \in B^*$. 
It is a free $B$-module of rank $1$ and any 
of its generators 
is called a {\em trace} of $B$. 

\smallskip 

The following construction yields a trace $\sigma$
canonically associated to the
complete intersection $F$. 

We take new variables $y:=\{y_1,\dots, y_n\}$, and we set 
$F_i^{(x)}:=F_i(x) \in A[x]$ and
$F_i^{(y)}:=F_i(y) \in A[y]$. 
Then $F_i^{(y)}-F_i^{(x)}$
belongs to the ideal $(y_1-x_1,\dots,y_n-x_n)$ and so 
there exist (not unique) 
$l_{ij}\in A[x,y]$ such that

$$
F_i^{(y)}-F_i^{(x)}=\sum_{j=1}^n l_{ij}(y_j-x_j), 
$$
for $i= 1, \dots, n$. 
We consider the determinant $\Delta \in A[x,y]$ of the square matrix 
$(l_{ij})_{ij}$, and we write it as
$$
\Delta=\sum_m a_m\,b_m
$$ 
with $a_m \in A[x]$ and $b_m\in A[y]$. 
Again, the polynomials $a_m,b_m$ are not uniquely defined.
The polynomial $\Delta \in A[x,y]$ 
is called a {\em pseudo-Jacobian determinant} of the
complete intersection $F$. 

\smallskip

Set $c_m:=b_m(x) \in A[x]$. Then there exists a
unique trace $\sigma \in B^*$ such that for $g\in A[x]$
$$ 
\overline{g} = \sum_m \sigma(\overline{g}\,\overline{a}_m)\,
\overline{c}_m  
$$
where the bar denotes class modulo $(F)$.

This is what is known as {\em Tate trace formula}.

\smallskip

Let $J:= \det (\partial F_i/ \partial x_j)_{i\, j}$ 
be the Jacobian determinant of the complete intersection
$F$ with respect to the variables $x_1,\dots,x_n$. 
Then the
following identity ---which justifies the name of pseudo-Jacobian for $\Delta$---
holds
$$ 
\overline{J} = \sum_m \overline{a}_m\,\overline{c}_m.
$$
 The standard trace $\Tr_V$ is related to $\sigma$ by the equality 
$$
\Tr_V (\overline g)=\sigma (\overline J\,\overline g) $$
for all $g\in A[x]$.


\subsubsection{A division lemma}

Throughout this subsection we keep notations and assumptions as in 
the previous one. 
In addition we replace $k$ by a number field $K$.

\smallskip

We will choose concrete polynomials $a_m, c_m$ which satisfy  
trace formula and we will estimate their degree and local height. 
Set $d:=\max_{i}\deg F_i$ and 
$h_v:=\max_i h_v(F_i)$ for $v\in M_K$.

\smallskip

First we  choose the polynomials $l_{ij}$. Remarking that $$
F_i^{(y)}-F_i^{(x)}= \sum_{j=1}^n
F_i(x_1,\dots,x_{j-1},y_j,\dots,y_n)-F_i(x_1,\dots,x_j,y_{j+1},\dots,y_n),
$$
 we set 
$$
l_{ij}:= (F_i(x_1,\dots,x_{j-1},y_j,\dots,y_n)-F_i(x_1,\dots,x_j,y_{j+1},
\dots,y_n))/(y_j-x_j)$$
Here we perform the division through the formula 
$$
(y_j^k-x_j^k)/ (y_j-x_j)= y_j^{k-1}+y_j^{k-2}\, x_j + \cdots +
y_j \, x_j^{k-2}+x_j^{k-1}.$$

We set $\Delta:=\det (l_{ij})_{i j}$.
Finally we choose $b_m \in A[y]$ as the monomials in the 
expansion of $\Delta$ 
with respect to $y$, 
$a_m \in A[x]$ as 
the corresponding coefficient, and we set $c_m:=b_m(x)$.

\smallskip

Set $F_i= \sum_{\alpha} a_{i\, \alpha} \, x^\alpha$ with 
$a_{i\, \alpha} \in A$. Then 
$$
l_{i j} = \sum_{\alpha} a_{i\, \alpha} \, 
x_1^{\alpha_1} \cdots x_{j-1}^{\alpha_{j-1}}
y_{j+1}^{\alpha_{j+1}} \cdots y_{n}^{\alpha_{n}}
(y_j^{\alpha_j-1}+ \cdots +x_j^{\alpha_j-1}) \in A[x, y].
$$
{}We deduce that 
$\, \deg l_{ij}\le d-1 \, $ and  $\, h_v(l_{ij})\le h_v$
for  every $v \in M_K$. Then $\deg \Delta \le n\, (d-1) $ and so 
$$
\deg a_m + \deg c_m \le n\, (d-1).
$$
We have also $h_v(c_m) =0$ and $h_v(a_m ) \le h_v(\Delta)$.

{}Finally  we can write $$
l_{i j} = C_0 + \cdots + C_{d-1} \, y_j^{d-1}, 
$$
where each $C_k \in A[x_1 , \dots, x_{j}, y_{j+1}, \dots, y_n] $
is a polynomial in $n+r$ variables of degree bounded by
$d-1$. 
This implies that the number of monomials of $l_{ij}$ is bounded
by $d\,{n+r+d-1\choose n+r} \le d\,(n+r+1)^{d-1}$.

  Therefore, for $v\in M_K^\infty$ we have 
\begin{eqnarray} \label{hvam}
h_v(a_m)&\le& h_v(\Delta) \nonumber \\[1mm]
&\le&  n\,h_v + (n-1)\,(\log d + (d-1)\,\log(n+r+1)) + n\,\log n \nonumber\\[1mm]
&\le& n\,(h_v + d\,\log (n+r+1)  + \log d).
\end{eqnarray}
Analogously we have $h_v(a_m)\le n\,h_v$ for $v \notin M_K^\infty$. 

\medskip
The following is a sharp estimate for the degree and the local height of the 
polynomials in the division procedure. 
It is a substantial improvement over \cite[Thm. 29]{KrPa96}.

\smallskip
We introduce the notation  $\deg_t f$  and $\deg_x f$  for the degree of 
a polynomial
$f\in A[x]$ with respect to the group of 
variables $t$ and $x$, respectively.

\begin{mainlem}{(Division Lemma)} 
\label{Division Lemma}

Set $A:= K[t_1,\ldots,t_r]$ and 
$A[x]:=   A[x_1,\ldots,x_n]$. 
Let $F:=\{F_1,\ldots,F_n\} \subset A[x]$  be a reduced complete 
intersection defining a variety 
$V:= V(F) \subset \A^{r+n}$ which satisfies Assumption \ref{assumption}.
Set $B:= K[V] = A[x]/(F)$.

Let $f, g \in A[x]$  be polynomials 
such that $\overline{f} \in B$ is a non-zero divisor and
$\overline{f}\, |\, \overline{g}$ in $B$. 
Set $d:=\max\{\deg f, \deg F_1, \dots, \deg F_n\}$ 
and  
$h_v:=\max\{ h_v(f), h_v(F_1), \dots, h_v(F_n)\}$ 
for $v\in M_K$. 

Then there exist  $q \in A[x]$ and $\xi \in K^*$ such that
\begin{itemize}
\item
$ \overline{q} \, \overline{f} = \overline g ,
$
\item
$\deg_{x} q  \leq  n \, d ,$
\item $\deg q \le \deg_t g + (n\, d + \max\{ (n+1) \, d, \deg_x g\} )\deg V $, 
\end{itemize}
\begin{itemize}
\item 
$ h_v(q)   \leq   h_v(g) + (n\, d+ \max\{ d, \deg_x g\} ) \,
h_v(V) \\[1mm]
     \hspace*{12mm} + \  ((n+1)\,h_v  
+   (r+6) \, \log(n+r+1) \, 
  (n\, d + \max\{(n+1)\, d, \deg_x g\}  ))\, \deg V \\[1mm]
     \hspace*{12mm} + \  2 \, \log (r+1) \, \deg_t g - \log |\xi|_v   
     $ \\[2mm]
    for $v \in M_K^\infty$,
\item 
$  h_v(q) \le  h_v(g) +  (n\, d + \max\{d, \deg_x g\}) \, 
h_v(V) + (n+1)\,h_v  \,\deg V  - \log |\xi|_v  $\\[2mm]
for $v \notin M_K^\infty$. 

\end{itemize}
\end{mainlem}

\begin{proof}{Proof.--} 
Set $L:= K(t_1, \dots, t_r)$ for the quotient field of  $A$ and 
$M:= L \otimes_A \, B$. 
Then $M$ is a finite $L$-algebra of dimension $\deg V$ and 
$\sigma $ can be uniquely extended to a $L$-linear map 
$\sigma :M\to M$. 

The fact that $B$ is a torsion-free $A$-algebra implies that the canonical map 
$B\to M$ is an inclusion. 

\smallskip

We will only consider the case $n\ge 1$. 
For the case $n=0$ we refer to Remark \ref{n=0}.
Whenever it is clear from the context, 
we will avoid explicit reference to the ring in which  we are considering 
a given 
element of $A[x]$.

\smallskip

Let $q_0 \in A[x]$ be any polynomial such that
$q_0 \, f = g$ in $B$. 
We have that $f$ is a non-zero divisor in $B$, and so it is invertible in 
$M$. 
Then
$q_0= f^{-1} \, g$ in $M$ and therefore
$ \ \sigma (f^{-1} \, g \, p) = \sigma (q_0 \, p) \in A \  $ 
for all $p \in A[x]$. 
Then we set 
\[
q:= \sum_m \sigma ({{f}}^{-1} \, {g} \,{a}_m)
\, c_m  \in A[x].
\]
Tate trace formula implies that ${q} \equiv {q}_0 \pmod{(F)}$, 
and so
$ \ {q} \, {f} = {g}$ in $B$.

\smallskip

Let $J \in A[x]$ 
denote the Jacobian determinant 
of the
complete intersection $F$ with respect to the group of 
variables $x$.
This is a non-zero divisor because of the 
Jacobian criterion, and so it is also 
invertible in $M$. 

Let $(J\,f)^*$ be the adjoint polynomial of 
$J\, f$ and set 
$$
\Lambda_m:= \Tr_V( {(J\,f)^*} {g}\,{a_m}) \in A.
$$ 
We have $J\, f \, (J\, f)^* = \No (J\, f) \in A\setminus \{0\} $, and so 
$$
\Lambda_m / \No(J\,f)
= \Tr (({{J}} \, {{f}})^{-1} \, {g} \,{a}_m)
= \sigma ({{f}}^{-1} \, {g} \,{a}_m)\in A .
$$
In particular 
$ \ \No(J\,f) \,|\, \Lambda_m \ $  in  $A$, and we have the expression 
$$ 
q= {1\over \No(J\, f)} \, \sum_m \Lambda_m \, c_m . 
$$

\medskip

Clearly
$$
\deg_{x} q \leq \max_m \deg c_m \leq n(d-1)\le n\,d. 
$$

Next we analyze the total degree of $q$. 
Let $g:=\sum_\alpha p_\alpha \,  x^\alpha$  be the monomial expansion 
of $g$ with respect to $x$. 
Then 
\begin{eqnarray} \label{Lambdam}
\Lambda_m= \sum_\alpha p_\alpha \,  \Tr ((Jf)^*x^\alpha\,a_m),
\end{eqnarray}
as $\Tr$ is a $A$-linear map. 
We have the estimates 
$\deg (J\,f)  \le  n\, (d-1)+ d  \le (n+1) \, d $ and $\deg (x^\alpha\,a_m)
\le \deg_x g + \deg a_m$, from where we get 
$$
\deg \Tr((Jf)^*x^\alpha\,a_m)  \le  \max\{(n+1) \, d,
\deg_x g + \deg a_m\} \deg V
$$
by Lemma \ref{traza}. 
Thus
\begin{eqnarray*}
\deg q &\le& \deg_t  g + \max_m\{\max\{ (n+1)\, d, \deg_x g + \deg a_m
\}\,  \deg V + \deg c_m\}\\
 &\le& \deg_t  g + \max_m\{\max\{ (n+1) \, d +\deg c_m, \deg_x g + \deg a_m
+ \deg c_m\} \}\, \deg V \\
&\le & \deg_t  g + \max\{ (n+1)\, d +n\, d, \deg_x g + n\, d\} \deg V \\[1mm]
&\le & \deg_t  g + (n\,d + \max\{ (n+1)\, d , \deg_x g \}) \deg V .
\end{eqnarray*}

\medskip

For the rest of the proof, we will use several times the 
following basic estimates:
\begin{eqnarray*}
\max\{ \deg (J\,f),\deg (x^\alpha\,a_m)\}  & \le & n\,d +  \max\{d,
\deg_x g \}  , \\[2mm]
\deg \Tr((Jf)^*x^\alpha\,a_m)  & \le &  (n\,d +  \max\{d,
\deg_x g \}) \deg V . 
\end{eqnarray*}

Finally we estimate the local height of $q$. 
Let $v\in M_K^\infty$.
We have $h_v(\partial F_i/\partial x_j ) \leq h_v + \log d$ and so 
$$
 h_v(J ) \leq n\,(h_v+ \log d ) + (n-1) \, \log (n+r+1) \, (d-1)
 + n\, \log n\le n\,(h_v + \log(n+r+1) \, d + \log d). 
$$
Therefore
\begin{eqnarray} \label{Jf}
h_v(J\,f)&\le& n(h_v +\log (n+r+1) \, d + \log d) + h_v + 
\log(n+r+1) \, d 
\nonumber \\[1mm]
& \le&  (n+1) \,h_v + (n+1)\,\log (n+r+1) \,d + n\,\log d 
\end{eqnarray}
by Lemma \ref{hprod}(b). 
We recall that 
$ \ h_v(x^\alpha\,a_m)\le    n\,(h_v + \log (n+r+1)\,d + \log d) \ $
by Inequality (\ref{hvam}) and so 
$$
\max\{ h_v(J\,f) , h_v(x^\alpha\,a_m)\} 
 \le (n+1) \,h_v + (n+1)\,\log (n+r+1) \,d + n\,\log d . 
$$

Then 
\begin{eqnarray*}
h_v(\Tr ((J\,f)^* 
x^\alpha a_m)) &\le & (n\,d + \max\{ d,\deg_x g\})
\, h_v (V) \\[1mm]
& & +\  ((n+1) h_v + (n+1) \,\log(n+r+1) \,d +n\, \log d + \log 2) 
\, \deg V\\[1mm]
& & + \ (r+1) \, \log(n+r+1) \, (n\,d + \max\{ d,\deg_x g\}) \, \deg V\\[2mm]
&\le & (n\,d + \max\{d,\deg_x g\}) \,  h_v(V) 
+ ((n+1)\, \,h_v + ((2 \, n+1)\,\log (n+r+1)\,d\, ) \, \deg V \\[1mm]
& &  + \  (r+1) \, \log(n+r+1) 
 \, (n\,d + \max\{ d,\deg_x g\}) \,\deg V \\[2mm]
&\le & (n\,d + \max\{d,\deg_x g\}) \, h_v(V) \\[1mm]
&& + \ ((n+1)\,h_v + (r+2) \, \log(n+r+1) \, (n\,d + \max\{(n+1) \, d,\deg_x g\})) 
\,\deg V 
\end{eqnarray*}
by Lemma \ref{traza}. Hence 
\begin{eqnarray*}
h_v(\Lambda_m) &\le & \max_\alpha \{
h_v(p_\alpha \,  \Tr ((J\,f)^* x^\alpha a_m) \} 
+ \log (n+1) \, \deg_x g \\ [2mm]
& \le & h_v (g) + \max_\alpha \{h_v(\Tr ((J\,f)^* x^\alpha a_m)) \}\\[1mm]
& & +\ \log(r+1)\,(n\,d + \max\{ d,\deg_x g\}) \, \deg V
+ \log (n+1) \, \deg_x g \\[2mm]
& \le & h_v (g) + 
(n\,d + \max\{ d,\deg_x g\})h_v(V) \\[1mm]
&& + \ ((n+1)\,
h_v + \ (r+2) \, \log(n+r+1)\, (n\,d + \max\{ (n+1) \, d,\deg_x g\}) ) 
\,\deg V \\ [1mm]
& & +\ \log(r+1)\,(n\,d + \max\{ d,\deg_x g\}) \, \deg V
+ \log (n+1) \, \deg_x g \\[2mm]
&\le& h_v (g) +  (n\,d + \max\{ d,\deg_x g\}) \, h_v(V) \\[1mm]
&& + \ ((n+1)\,h_v + (r+4) \, \log(n+r+1)\, (n\,d + \max\{ (n+1)d,\deg_x g\}) ) 
\,\deg V
\end{eqnarray*}

by application of Identity (\ref{Lambdam}) and Lemma \ref{hprod}(b).
We have 
$$
h_v(q) \le  \max_m \, \{h_v (\Lambda_m / \No (J\,f))\}
$$
as each $c_m$ is a different monomial in $x$. 
Thus it only remains to estimate the local height of each 
$\Lambda_m / N(J\, f)$.
Let 
 $\xi\in K^*$ be any non-zero coefficient of $\No(J\,f)$. 
Then 
\begin{eqnarray} \label{Lambdam/N(Jf)}
\log |\Lambda_m / \No (J\,f)|_v  
&\le & h_v(\Lambda_m ) +  2\,\log(r+1) \, ( \deg _t g+ 
(n\,d + \max\{ d,\deg_x g\})\, \deg V ) - \log |\No(J\, f) |_v \nonumber \\[2mm]
&\le &  h_v(g)+ (n\,d + \max\{d,\deg_x g\}) \, h_v(V) \nonumber \\[1mm]
&& + \ ( (n+1) \,h_v  +  (r+6) \, \log(n+r+1) \, 
(n\,d + \max\{(n+1) d,\deg_x g\}) ) \,\deg V \nonumber \\[1mm]
&& + \  2\,\log(r+1) \,  \deg _t g - \log |\xi |_v
\end{eqnarray}
by Lemma \ref{hprod}(d) and the fact that $\log |\xi|_v  \le \log |\No(J\, f)|_v$. 
{}From Lemma \ref{norma} and Inequality (\ref{Jf}) we obtain 
\begin{eqnarray}\label{normaJf}
\log|\xi|_v & \le & h_v(\No(J\,f))\nonumber \\[2mm] 
& \le &  (n+1)\,d\,h_v(V) +  
((n+1)\,h_v +   (n+1) \, \log (n+r+1) \,  d + n\,\log d)\deg V \nonumber\\[1mm]
& &  + \ (r+1)\, (n+1) \,\log (n+r+1)\,d\,\deg V \nonumber\\[1mm]
&\le & (n+1)\,d\,h_v(V) +  
((n+1)\,h_v + (r+3)\,(n+1)\,\log(n+r+1)\,d  )\,\deg V  
\end{eqnarray}
This implies that the right hand side of Inequality (\ref{Lambdam/N(Jf)})
is non-negative.
So the inequality also holds for $h_v(\Lambda_m/\No(J\, f))$, and thus for 
$h_v(q)$.

\smallskip

The case $v \notin M_K^\infty$ is treated analogously. 
We remark that the election of $\xi$ is independent of 
$v$, and so it can be done uniformly. 

\end{proof}

\begin{rem} \label{n=0}

Let notations be as in the previous lemma. 
In case $n=0$ we have the sharper estimates
\begin{itemize}

\item $ \deg q \le \deg g $, 

\item $ h_v(q) \le h_v(g) + h_v+ 2\, \log (r+1) \, \deg g -\log |\xi|_v \  $
for $v \in M_K^\infty$, 

\item $ h_v(q) \le h_v(g)  + h_v-\log |\xi|_v \  $
for $v \notin M_K^\infty$.

\end{itemize}

Here $\xi \in K^*$ denotes any non-zero coefficient of $f$. 
The local 
height estimates follow from Lemma \ref{hprod}(d) and the fact that 
$h_v -\log|\xi|_v  \ge 0 $. 

\end{rem}


\subsection{An effective arithmetic Nullstellensatz}


\subsubsection{Estimates for the complete intersection case}

The following result gives estimates for 
the degree and local height  of  polynomials 
arising in the Nullstellensatz over a number field $K$
in case the input is a reduced weak regular sequence,
i.e. when the input is a reduced regular sequence
which eventually  may  have no common zeros in $\A^n$. 
It is a direct consequence of the division lemma above. 

These estimates depend mainly on the degree and height of the varieties 
successively cut out by the input polynomials. 
They are quite flexible apply to other situations as 
 we will see in Chapter 4. 

\begin{lem} \label{nullstlocal}
Let $n\ge 2$ and let $f_1,\ldots,f_s \in K[x_1,\ldots,x_n]$ 
be polynomials without 
common zeros in  $\A^n$ which form 
a reduced  weak regular sequence.
Furthermore assume that for $j= 1, \dots, s-1$,
 $V_j:=V(f_1, \ldots, f_j) $ satisfies 
Assumption \ref{assumption}.

Set $d:= \max_i \{\,\deg f_i, 2\,\} $ and 
$h_v:= \max_i h_v(f_i)$ for $v\in M_K$. 
 
Then there exist  
$p_1,\ldots,p_s \in K[x_1,\ldots,x_n]$ and $\xi \in K^*$
such that 
\begin{itemize}
\item $1 = {p}_1 \, {f}_1 
+ \cdots + {p}_s \,
{f}_s $,
\item $ \deg p_i \le 2\,n\,d \,(1+{\sum_{j=1}^{\min\{n,s\} -1}
 \deg V_j)}$,
\item 
$
h_v (p_i) \leq 
2\,n\,d\,\sum_{j=1}^{s-1}  h_v(V_j) 
\ + \ ((n+1)\,h_v + 2 \, n \, (2\, n+5) \,\log(n+1)\, d \, ) 
\, (1 + \sum_{j=1}^{s-1} \deg V_j)
- \log|\xi|_v $ \\[2mm] 
for $v \in M_K^\infty$, 
\item 
$
h_v (p_i) 
 \leq   
2\,n\,d\,\sum_{j=1}^{s-1}  h_v(V_j) 
+ (n+1)\,h_v\,(1+\sum_{j=1}^{s-1}  \deg V_j) 
- \log|\xi|_v $ \ \ 
for $v \notin M_K^\infty$.
\end{itemize}

\end{lem}
 
\smallskip
\begin{proof}{Proof.--}
Set  
$I_i:= I(V_i) = (f_1, \dots, f_i)$ for $i=0, \dots, s-1$.  
Also set $ \ A_i:=K[x_1,\ldots,x_{n-i}] \ $ and 
$ \ B_i:= K[V_i]= K[x_1,\dots,x_n]/I_i $. 
The fact that $V_i$ satisfies Assumption \ref{assumption} implies that the 
inclusion
$ \ A_i \hookrightarrow B_i \ $ 
is integral. 

\smallskip

We note that the sets of free and dependent variables of $B_i$ 
have cardinality $n-i$ and $i$  respectively. 
Also the set of dependent variables of $B_j$
is contained in that of $B_i$  for $i\le j$. 

\smallskip

For $f\in K[x_{1},\ldots,x_n]$ we denote by 
$\deg_{x(i)}f$ the degree of $f$ in the dependent variables 
$x_{n-i+1},\ldots,x_n$ of $B_i$ with respect to the integral inclusion 
$ A_i \hookrightarrow B_i$.
For $i\le j$, the previous observation implies that 
$  \ \deg_{x(j)}f \le   \deg_{x(i)}f $.

\medskip

Applying the Division Lemma \ref{Division Lemma}, we  
will construct inductively polynomials
$p_1, \dots, p_s$: 
first we take $p_s$ such that 
$$ 
p_s \, f_s \equiv 1 \ \pmod{I_{s-1}} . 
$$
For $0\le i \le s-2 $ we assume
that $p_{i+2}, \dots, p_s$ are already constructed and  we set
$$
b_{i+1}:= 1- (p_{i+2} \, f_{i+2} 
+ \cdots + p_s \, f_s) .
$$
Then $f_{i+1} $ is a non-zero divisor and $f_{i+1}\, |\, b_{i+1} $ in 
$B_i$. 
We apply again Division Lemma to obtain $p_{i+1} $ such that 
$$
p_{i+1} \, f_{i+1} \equiv b_{i+1}  \ \pmod{I_i}, 
$$
Continuing this procedure until  $i=0$,
we get $1=p_1\,f_1+\cdots +p_s\,f_s$ in $K[x_1,\dots,x_n]$.

\medskip

Let us  analyze  degrees.

First we consider the case
$s\le n$. 
Again we proceed by induction.

First we have 
$ \ \deg_{x(s-1)} p_s \le (s-1) \, d\le (n-1)\,d \ $ and 
$
\deg p _s \le (2\,s-1) \,d\,\deg V_{s-1} .
$

Now let $1\le i\le s-2$. Then 
$ \ \deg_{x(i)} p_{i+1}  \leq i \, d \ $ and 
$$
\deg p_{i+1}  \leq  \deg b_{i+1} +  (i \, d + \max \{(i+1)\, d , 
\deg_{x(i)}  b_{i+1} \}) \, \deg V_i .
$$

where
$$
\deg_{x(i)}b_{i+1}\le \max_{j\ge i+2} \{ \deg_{x(i)} p_j
+\deg f_j \} \le
\max_{j\ge i+2} \deg_{x(j-1)}p_j + d \le s\,d.
$$
Hence
\begin{eqnarray*}
\deg p_{i+1} & \le&  \max_{j\ge i+2}  \deg p_j +d+  (s+i)
\, d\,\deg V_i \nonumber \\ [-1mm]
& \leq & (2\,s-1) \,d\,\deg V_{s-1} + 
 \sum_{j=i}^{s-2}(d+ (s+j)\, d\,\deg V_j ) \\[-2mm]
& = &  (s-i-1)\, d + \sum_{j=i}^{s-1}(s+j)\, d\,\deg V_j .
\end{eqnarray*}
For $i=0$  we have $p_1\,|\,b_1$ and 
therefore  $\deg p_1\le \deg b_1\le \max_{j\ge 2} \deg p_j + d$. 
Then for all $i$:
\begin{eqnarray*}
\deg p_i \le (s-1)\,d +
  \sum_{j=1}^{s-1} (s+j) \, d\, \deg V_j \le 2\,n\,d\,
(1+\sum_{j=1}^{s-1}\deg V_j) .
\end{eqnarray*}


Next we consider the case $s=n+1$. 
In this case $V_s$ is a 0-dimensional 
variety and so 
$$
\deg p_{n+1}= \deg_{x(n)} p_{n+1} \le n \, d . 
$$

Let $1\le i \le n-1$. Then 
$ \ \deg_{x(i)} p_{i+1} \le i\, d \ $ and 
\begin{eqnarray*}
\deg p_{i+1} & \leq & \max_{j\ge i+2} \deg p_j + d + (n+1+i) \, d \, \deg V_i
\nonumber \\ [-1mm]
& \le & n \, d 
+ \sum_{j=i}^{n-1} (d + (n+1+j) \, d \, \deg V_j ) \nonumber \\[-2mm]
& =  &   (2\, n-i) \, d + \sum_{j=i}^{n-1} (n+1+j) \, d \, \deg V_j . 
\end{eqnarray*}

We have also $ \ \deg p_1\le \deg b_1\le \max\deg_{j\ge 2}
\deg p_j + d$. We conclude for all $i$:
\begin{eqnarray*}
\deg p_i \le 2\,n\,d + 
\sum_{j=1}^{n-1} (n+1+j) \, d \, \deg V_j \le 2\,n\,d\,(1
+\sum_{j=1}^{n-1}\deg V_j) . 
\end{eqnarray*}


Finally we estimate the local height of these polynomials. 
In the rest of the proof we will make repeated use of the following 
degree bounds: 
\begin{eqnarray*}
\deg_{x(i-1)} p_{i} & \le&  n \, d ,  \nonumber \\ [-2mm]
\deg p_i & \le  & 2\, n \, d \, 
(1+ \sum_{j=i-1}^{\min\{ n, s\}-1}   \deg V_j ) .
\end{eqnarray*}

As usual, we consider only the case $v\in M_K^\infty$, the case 
$v\notin M_K^\infty$ can be treated analogously. 
From Division Lemma we obtain 
\begin{eqnarray*}
h_v(p_s)  & \leq  &   s\, d \,h_v(V_{s-1}) + ( s\,h_v+
 (n - (s-1) +6)\,(s+(s-1))\,  \log(n+1) \, d ) \, \deg V_{s-1} - 
\log |\xi_{s-1}|_v 
\end{eqnarray*}  

\smallskip

for some $\xi_{s-1} \in K^*$.

Let $1\le i \le s-2 $
and set $n_0:= \min\{ n, s\}$. Then there exists $\xi_i \in K^*$ such that 
\begin{eqnarray*}
h_v(p_{i+1})  & \leq  &   h_v(b_{i+1}) + 
(i\,d + \max\{d,\deg_{x(i)} b_{i+1}\}) \, h_v(V_i) \\[1mm]
&& + \ ((i+1)\,h_v + \ (n-i+6)\,\log(n+1)\,(i\,d\, + \max
\{ (i+1) \, d,\deg_{x(i)}b_{i+1}\})) \, \deg V_i \\[1mm]
&& + \ 2\, \log (n-i+1) \, \deg b_{i+1} 
- \log |\xi_i|_v \\[2mm]
&\le & \max_{j\ge i+2} h_v (p_j) + h_v + \log(n+1)\,d + \log (s-i) 
 + (s+i)\,d\,h_v(V_i) + (i+1)\,h_v \, \deg V_i\\[-1mm]
&& + \ (n-i+6)\,(s+i)\,\log(n+1)\, d  \, \deg V_i 
+ 2\, \log (n+1) \, ( 2\, n \, d \, ( 1+ \sum_{j=i+1}^{n_0 -1}\deg V_j) + d) 
 \\[-1mm] && 
- \log |\xi_i|_v  .
\end{eqnarray*}
Applying  the inductive hypothesis  we  obtain
\begin{eqnarray*}
h_v(p_{i+1})  & \leq  &   
 s\, d\, h_v(V_{s-1})+  d\, \sum_{j=i}^{s-2} (s+j) \, h_v(V_j) 
 + (s-i-1) \, h_v + h_v \, 
\sum_{j=i}^{s-1} (j+1) \, \deg V_j \\ [-1mm]
& & + \ 4\, (s-i-1) \, (n+1)  \, \log (n+1) \, d 
+  \log(n+1)\, d \,
\sum_{j=i}^{s-1} (n-j+6)\, (s+j)\,\deg V_j \\[-1mm]
& & 
+ \ 4\, n \, \log (n+1) \, d \sum_{j=i+1}^{n_0 -1} (j-i ) \, \deg V_j 
- \sum_{j=i}^{s-1} \log |\xi_j|_v . 
\end{eqnarray*}

For  $i=0$ we apply Remark \ref{n=0}: 
 there exists $\xi_0 \in K^*$ such that  
\begin{eqnarray*}
h_v(p_1) & \le & h_v(b_1) +h_v+ 2\,\log(n+1)\deg b_1 -
\log |\xi_0|_v \\
& \le & 
\max_{j\ge 2} h_v (p_j) + 2\, h_v + \log(n+1)\,d + \log s
+ 2\, \log (n+1) \, ( 2\, n \, d \, ( 1+ \sum_{j=1}^{n_0-1}\deg V_j) + d) 
\\
&  & -\ \log |\xi_0|_v.
\end{eqnarray*}

We  set $\xi:= \prod_{j=0}^{s-1} \xi_j$. Then 
\begin{eqnarray*}
h_v(p_{1})  
&\le & 2\, n\, d\, \sum_{j=1}^{s-1} h_v(V_j) + (n+1) \, h_v \, 
(1+ \sum_{j=1}^{s-1} \deg V_j ) + 4\, n \, (n+1) \, \log (n+1) \, d 
\\
& & 
+ \ \log(n+1)\, d \,
\sum_{j=1}^{s-1} (n-j+6)\, (n+1+j) \, \deg V_j 
+ 
4\, n \, \log (n+1) \, d \, \sum_{j=1}^{n_0 -1}  j \, \deg V_j 
- \log |\xi|_v \\
&\le & 2\, n\, d\, \sum_{j=1}^{s-1} h(V_j) + ( (n+1) \, h_v 
+  2\, n \, (2\, n +5) \, \log(n+1)\, d ) \,
(1+ \sum_{j=1}^{s-1} \deg V_j ) - \log |\xi|_v .
\end{eqnarray*}  
This last inequality follows from the facts 
that $ \ 4\, n \, j + (n-j+6)\, (j+s) \le 2\, n \, (2\, n+5) \ $ 
for $j\le n-1,$ and 
$6\, (2\, n+1) \le 2\, n \, (2\, n+5) \ $ as $n\ge 2$.

To conclude the proof, observe that for
$i=1,dots,s-1$, Inequality \ref{normaJf} guarantees
that the obtained estimate  for $p_i$
differs from the one for $p_{i+1}$ by a positive term.
Thus, the same estimate holds for $h_v(p_i)$, $1\le
i\le s$.
\end{proof}

By means of B\'ezout inequality, we can now 
estimate the degree and height of the 
varieties $V_j$.
In this way we obtain an estimate which only depends on 
the degree and height 
of the input polynomials.

\begin{cor} \label{extrinsecolocal}

Let notations and assumptions be as in Lemma \ref{nullstlocal}. Then 
there exist  
 $p_1,\ldots,p_s \in K[x_1,\ldots,x_n]$ and $\gamma \in K^*$
 such that 
 
\begin{itemize}
\item $1 = {p}_1 \, {f}_1 
+ \cdots + {p}_s \,
{f}_s $, 
\item $ \deg p_i \le 4\,n\,d^n$, 
\item 
$
h_v (p_i) \leq 
4\,n\,(n+1)\,d^n\,h_v + 4\, n \, (4\, n+5) \,\log
(n+1)\,d^{n+1} - \log |\gamma|_v $ \ \ for $v \in M_K^\infty$,
\item 
$h_v (p_i) \leq   
4\,n\,(n+1)\,d^n\,h_v - \log |\gamma|_v $ \ \ for $v \notin M_K^\infty$.

\end{itemize}
\end{cor}

\begin{proof}{Proof.--}
Let us first consider  degrees. 
We assume without loss of generality $d\ge 2$. 
{}From
the preceding result we obtain 
$$
\deg (p_i) \le
  2\,n\,d \,(1+\displaystyle{\sum_{j=1}^{\min\{n,s\} -1}
 \deg V_j)}
\le  2\,n\,d\, (1+\cdots+d^{n-1})
\le  4\,n\,d^n . 
$$
Next we consider the local height estimates. 
Let $v\in M_K^\infty$. 
We have 
$$
h_v (p_i) \leq 
2\,n\,d\,\sum_{j=1}^{s-1}  h_v(V_j) 
+ ((n+1)\,h_v +
2\, n\, (2\, n +5)  \, \log(n+1) d\,) \,(1 + \sum_{j=1}^{s-1} \deg V_j)
- \log|\xi|_v 
$$
for some $\xi \in K^*$. 
Applying
Corollary \ref{bezloc1}, 
$\ h_v(V_j)\le
j\,d^{j-1}\,h_v+(n+j)\log(n+1)\,d^j - \log |\lambda_j|_v \ $ for some 
$\lambda_j \in K^*$.
Therefore
\begin{eqnarray*}
h_v (p_i) &\leq &  
2\,n\,d\,\sum_{j=1}^{s-1} (j\,d^{j-1}\,h_v+(n+j)\log(n+1)\,d^j - 
\log |\lambda_j|_v)\\ [-2mm]
& &  + \  ((n+1)\,h_v +
2\, n\, (2\,  n+5) \,\log(n+1) \, d\, ) \,\sum_{j=0}^n d^i 
 - \log|\xi|_v\\[0mm]
&\le & 4\,n^2\,d^n\,h_v + 8\,n^2\,\log(n+1)\,d^{n+1} \\ [-1mm]
& & +  \ 2\,(n+1)\,d^n\,h_v  + 4\, n\, (2\, n+5)\, \log(n+1)\,d^{n+1}
 -\ 2\,n\,d\,\sum_{j=1}^{s-1} \log|\lambda_j|_v \ - \
\log |\xi|_v\\
& \le & 4\,n\,(n+1)\,d^n\,h_v 
+ 4\, n\, (4\, n+5)\, \log(n+1)\,d^{n+1} - \log|\gamma|_v, 
\end{eqnarray*}
where $\gamma\in K^*$ is defined as 
$ \ \gamma:=\xi\,\prod_{j=1}^{s-1}\lambda_j^{2\, n\, d}$. 

The case $v\notin M_K^\infty$ follows anagolously.
\end{proof}


\subsubsection{Proof of Theorem 1}

In order to prove Theorem 1, 
it only remains to put the case of a general input into the 
hypothesis of Corollary \ref{extrinsecolocal}. 
This is accomplished by replacing the input polynomials and variables
by generic linear combinations. 
The coefficients of the linear 
combinations will be chosen to be roots of 1. 
Amazingly enough,  we will see using next lemma that we 
don't need to control  
the degree of the involved number field extension. 

\medskip

Let $L$ be a finite extension of $K$, and let 
$\cB:= \{e_1, \dots, e_N\}$ be a basis 
of $L$ as a $K$-linear space. 
We recall that $\cB^*:= \{e_1^*, \dots, e_N^*\}$ is 
the {\em dual basis} of $\cB$ 
if 
$ \ \Tr_K^L(e_i \, e^*_j) =1 \ $ for $ i=j $ and $0$ otherwise.

\begin{lem} \label{lang}
Let $\omega \in \Qbarra$ be a primitive $p$-root of 1 for some prime $p$. 
Then the basis $ \ \cB^* := 
\{ \,  (\omega^{-j} - \omega) \, /\, p \, : \, j =0, \dots, p-2 \, \} \ $
of $\Q(\omega)$ is dual to 
$ \ \cB := \{  \, \omega^{i} \, : \, i=0, \dots, p-2\, \}$. 

\end{lem}

\begin{proof}{Proof.--}
A direct computation shows that for $i, j =0, \dots, p-2$
$$
\Tr (\omega^i\, (\omega^{-j} -\omega ) )  =
 \sum_{l=1}^{p-1} \omega^{l\, i }  \, (\omega^{ - l \, j} -\omega^l ) =
\left \{ 
\begin{array}{lcl}
p , &\mbox{ for }& i=j   , \\ [1mm]
0  , & \mbox{ for }&  i \not= j .
\end{array}
\right. 
$$

\end{proof}

We will use this result in the following way:
let $\omega $ be a primitive $p$-root of 1 and set 
$L:= K(\omega)$. 
Let us assume that $\Q(\omega) $ and $K$ are linearly independent and that 
$p$ does not divide the discriminant of $K$. 
Both conditions are satisfied by all but a finite number of $p$. 
Then $[L:K]=p-1$ and ${\cal O}_L={\cal O}_K[\omega]$ \ 
\cite[Ch. III, Prop. 17]{Lang70}.

Now, let $\nu \in L \setminus \{0\}$. Then
$$
\nu= {1\over p} \, \Tr (\nu \, ( 1-\omega ) ) + \cdots + 
{1\over p} \, \Tr (\nu \, (\omega^{2-p} -\omega ) ) \, \omega^{2-p} 
\ \in K[\omega] \setminus \{ 0\}
$$ 

and so  
there exists $0\le j \le p-2$ such that 
$ \ \Tr ( \nu \, (\omega^{-j} -\omega ) ) \, 
/ \,  p \ \in K \setminus \{0\}$. 

Moreover, if $\nu
\in {\cal O}_L \setminus \{0\}$, as every coefficient belongs 
to $\cO_K$, there exists $0\le j \le p-2$ such that 
$ \ \Tr ( \nu \, (\omega^{-j} -\omega ) ) \, 
/ \,  p \ \in {\cal O}_K \setminus \{0\}$.

\bigskip

\begin{teo}{(Effective arithmetic Nullstellensatz)} \label{arithnullst}

Let  $K$ be a number field and let
$f_1,\ldots,f_s \in {\cO}_K[x_1,\ldots,x_n]$
be polynomials without 
common zeros in  $\A^n$. 
Set $d:= \max_i \deg f_i $ and $h:= h(f_1, \dots , f_s)$. 
 
Then there exist 
$ a\in {\cal O}_K\setminus\{0\}$ and 
$g_1,\ldots,g_s \in {\cal O}_K  [x_1,\ldots,x_n]$ 
such that 

\begin{itemize}
\item $a = {g}_1 \, {f}_1 
+ \cdots + {g}_s \,
{f}_s $, 
\item
$ \deg g_i \le 4\, n \, d^n $, 
\item
$h(a, g_1, \dots, g_s) \leq   4\, n\,(n+1) \, d^n \,(h + \log s
+ 
(n+7) \, \log (n+1) \, d).
$
\end{itemize}
\end{teo}

The extremal cases   $n=1$ and $d=1$ are treated directly
in the following results.

\begin{lem} \label{d=1}
Let $\ell_1, \dots, \ell_s \in \cO_K[x_1, \dots, x_n] $ 
be polynomials of degree 
bounded by 1 without common zeros in $\A^n$.
Set $h:=  h(\ell_1,\dots,\ell_s)$.

Then there exist $a\in \cO_K\setminus \{0\}$ and $a_1, \dots, a_s \in \cO_K$
 such that 
\begin{itemize}
\item $ \ a=a_1 \, \ell_1 + \cdots + a_s \, \ell_s $, 

\item $h(a,a_1,\dots,a_s) \le (n+1)\, ( h + \log (n+1) )$.
\end{itemize}
\end{lem}

\begin{proof}{Proof.--}
Equation
$ \ a=a_1 \, \ell_1 + \cdots + a_s \, \ell_s \ $ is equivalent 
to a $\cO_K$-linear system of $n+1$ equations in $s$
unknowns, which can be solved applying
 Cramer rule. The integer $a$ is the determinant of a non-singular
submatrix of the system.
\end{proof}

\begin{lem} \label{n=1}
Let $f_1, \dots, f_s \in \cO_K[x]$ be polynomials 
without common zeros in $\A^1$. 
Set $d:= \max_i \deg f_i$ and $h:=  h (f_1,\dots,f_s))$. 

Then there exist $a \in \cO_K\setminus \{0\}$ and 
$g_1, \dots, g_s  \in \cO_K[x]$ such that 
\begin{itemize}

\item $a = {g}_1 \, {f}_1 + \cdots + {g}_s \, {f}_s $,

\item $ \deg g_i \le d-1 $,

\item $ h(a,g_1,\dots,g_s) \leq 2\, d \, (h + d)$.
\end{itemize}
\end{lem}

\begin{proof}{Proof.--}

Let 
$ \ f:= \sum_i a_i \, f_i  , \ g:= \sum_i b_i \, f_i \in K[x] \ $
be generic linear combinations of $f_1, \dots, f_s$. 
Then $f$ and $g$ are coprime polynomials, and so there exist $p, q  \in K[x]$ 
with $\deg p < \deg g $ and $\deg q < \deg f$ such that 
$
1 = p\, f + q\, g. 
$

Expanding this identity there exists 
$p_1, \dots, p_s  \in K[x]$ with $\deg p_i \le d-1$ such that 
$$
1 = {p}_1 \, {f}_1 + \cdots + {p}_s \, {f}_s .
$$
Thus the above B\'ezout identity translates to 
a {\em consistent} system of $K$-linear equations. 
The number of equations and variables equal $2\, d$ and $s\, d$ respectively. 
This system can be solved by Cramer rule. The integer
$a$ is  the determinant of a non-singular $2\, d
\times 2\, d-$submatrix 
of the matrix of the linear system.
\end{proof}

\begin{proof}{Proof of Theorem 1.--} 

We assume $n>1$ and $d>1$.

\smallskip
Let $G_p \subset \Qbarra $ denote the group of
$p$-roots of 1, for a prime $p$. 
For $ a_{i j } \in G_p$ and $i=1, \dots, \min\{n+1, s\}$ we set 
$$
q_i := a_{i\, 1} \, f_1+\dots+ a_{i\, s} \, f_s . 
$$
Also, for $b_{k l} \in G_p$ and $k=1, \dots, n$ we set 
$$
y_k:= b_{k\, 0}+ b_{k \, 1} \, x_1 + \cdots + b_{k \, n} \, x_n .
$$ 

We will assume that for a specific choice of $p$, $a_{ij}$
and $b_{kl}$  there exists 
$t \le \min\{n+1, s\}$ such that  
$ (q_1, \dots, q_i)\subset K[x_1, \dots, x_n]$
is a radical ideal of dimension $n-i$ for 
$i=1, \dots, t-1$ and $1\in (q_1, \dots, q_t)$. 
We also assume that $y_1, \dots, y_n$ is a linear change of variables, 
and that 
$V_i:= V( q_1, \dots, q_i) \subset \A^n$ 
satisfies Assumption \ref{assumption} for $i=1, \dots, t-1$
with respect to $y_1,\dots,y_{n-i}$. 

\smallskip 

This is guaranteed by the fact that these conditions are generically 
satisfied:
there exists a hypersurface $H$ of the coefficient space
such that $(a_{i\, j} , b_{k\, l}) \notin H$ 
implies that $q_1, \dots, q_s$ satisfy the stated  conditions
with respect to the variables $y_1, \dots, y_n$ \cite{GiHeSa93}, 
\cite{SaSo95}. 
As $\cup_p G_p$ is Zariski dense in $\A^1$, it follows that 
these coefficients can be chosen to lie in $G_p$ for 
some $p$. Moreover, $p$ can be chosen such that 
for  $\omega$  a primitive $p$-root of 1 and  $L:= K(\omega)$,
$\Q(\omega) $ and $K$ are linearly independent and 
$p$ does not divide the discriminant of $K$. 

\smallskip

We refer the reader 
to Section \ref{Equations in general position}, where  
we give a self-contained treatment 
of this topic.

\medskip

Set $ \ b:= (b_{k \,0})_k \in G_p^n \ $ and 
$\ B:=(b_{k\, l} )_{ k,l \ge 1 } \in \GL_n(\Qbarra)$. 
For $j=1, \dots, t$ set
$$
F_j(y) := q_j(x)= q_j (B^{-1}(y -b)) \ \in L[y_1,\dots,y_n] . 
$$
Then $F_1, \dots, F_t$ satisfy the hypothesis
of Corollary \ref{extrinsecolocal}. Let 
$\gamma \in L^*$ and 
$P_1, \dots, P_t \in L[y_1, \dots, y_n]$ be the non-zero element 
and the polynomials satisfying B\'ezout identity
 we obtain there.

Now, for $i=1,\dots, s$, set
$$
p_i:= \sum_{j=1}^t a_{i\, j} \, P_j(B\, x+b) 
\in L[x_1, \dots, x_n] 
$$
so that $  \ 1=p_1 \, f_1+ \cdots + p_s \, f_s$ holds.

Finally set
$  \ \mu:= (\det B)^{4\, n\,(n+1)\,d^{n+1}}\, \gamma \ \in L^* \ $. 
By Lemma \ref{lang} there exists $0\le \ell \le p-2$ such that 
$\ \Tr(\mu \, (\omega^{-\ell} -\omega) ) \ne 0$. 

We define
$$
a:= \Tr (\mu\,(\omega^{-\ell}-\omega)) / p \, \in K^* \quad \quad ,
\quad \quad g_i :=
\Tr (\mu\,p_i\,(\omega^{-\ell}-\omega))  / p \, \in K[x_1, \dots, x_n]
$$
for $i=1, \dots , s $.

Then
$$
a=g_1\,f_1+\cdots+g_s\,f_s
$$
as $f_1,\dots,f_s \in K[x_1,\dots,x_n]$ and $\Tr$ is a $K$-linear map.

Aside from the degree and height bounds, we will show
that since $f_1,\dots, f_s \in \cO_K[x_1,\dots,x_n]$, 
 $a\in {\cal O}_K$ and $g_i\in {\cal O}_K[x_1,\dots,x_n]$. 

\medskip

Let us first analyze degrees and local heights.

As  $\deg F_j \le d$,
$
\deg g_i\le  \deg p_i \le \max_j \deg P_j \le 4\, n\, d^n
$.

Now let $v\in M_K^\infty$ and let $w\in M_L$ such that $w\, |\, v$. 
We have 
$\ h_w(B^{-1} \, (y- b))  \le n\, \log n - \log |\det B|_w $
and so 
\begin{eqnarray*}
h_w(F_j) & \le & h_w (q_j) + (n \, \log n -
\log |\det B|_w + 2\, \log (n+1) ) \, d \\[2mm]
&\le  &  h_v  + \log s + (n+2) \, \log (n+1) \, d - 
\log|\det B|_w \, d
\end{eqnarray*}
by Lemma \ref{hprod}(c). 
From Corollary \ref{extrinsecolocal} 
\begin{eqnarray*}
h_w (P_j) &\leq &  
4\,n\,(n+1)\,d^n  \max_k h_w(F_k) +
4\, n \, (4\, n+5) \,\log (n+1)\,d^{n+1} - \log|\gamma|_w \\[2mm]
&\le&4\,n\,(n+1)\,d^n(h_v  + \log s + (n+2) \, \log (n+1) \, d - 
\log|\det B|_w \, d) \\[1mm]
& & +
\ 4\, n \, (4\, n+5) \,\log (n+1)\,d^{n+1} - \log|\gamma|_w \\[2mm]
& = & 4\,n\,(n+1)\,d^n\,(h_v  + \log s)  + 
4\, n\, (n^2 +7\, n +7) \,\log (n+1)\,d^{n+1} - \log |\mu|_w . 
\end{eqnarray*}

Therefore
\begin{eqnarray} \label{uniforme} 
h_w(\mu \, p_i) &\le& \max_j h_w(P_j) + 2\, \log (n+1) \, \max_j \deg P_j
 + \log t + \log |\mu|_w  \nonumber \\
&\le & 4\, n\,(n+1) \, d^n \,  (h_v +\log s)  + 4\, n\, (n^2 +7\, n +7) 
\, \log (n+1) \, d^{n+1}  \nonumber \\[1mm]
& & +\, 8\,n\,\log (n+1)\,d^n + \log (n+1) \nonumber \\[1mm]
&\le & 4\, n\,(n+1) \, d^n \,(  h_v + \log s + 
(n+7)\,
\log(n+1)\,d ) - \log 2  
\end{eqnarray}
again by Lemma \ref{hprod}(c) and the fact $d,n\ge 2$.
We have 
$$  
g_i= {1\over p}\,  \Tr (\mu\,p_i\,(\omega^{-\ell}-\omega))  = 
 {1\over p}\, \sum_{\sigma \in \Gal_{L/K}}
\sigma(\mu\,p_i\,(\omega^{-\ell}-\omega)) \\[-4mm]
$$
and so 
\begin{eqnarray*}
h_v(g_i) 
& \le & 
\max_{w\, | \, v} \, h_w ( \mu\, p_i) + \log 2 \\
& \le  & 4\, n\,(n+1) \, d^n \,(  h_v + \log s + 
(n+7) \, 
\log(n+1)\,d ) . 
\end{eqnarray*}

We have $ \ h_w(\mu ) \le 
4\, n\,(n+1) \, d^n \,(  h_v + \log s) + 
4\,n\,(n^2+7\, n +7 )\,
\log(n+1)\,d^{n+1}  \ $ and so the previous estimate also 
holds for $h_v(a)$.

\smallskip

Now let $v\notin M_K^\infty$ and $w\, | \, v$. 
Analogously we have 
$$
h_w( \mu ), \ h_v( \mu \, p_i)  \le 4\, n\,(n+1) \, d^n \, h_v =0
$$
as $f_1, \dots, f_s \in \cO_K[x_1, \dots, x_n]$. 
Then 
$\ \mu \in {\cal O}_L\setminus \{0\} \ $ and  $ \ \mu \, p_i \in {\cal
O}_L[x_1,\dots,x_n] \ $, which in term implies that 
$\ a\in {\cal O}_K\setminus \{0\} \ $ and  $ \ g_i \in {\cal
O}_K[x_1,\dots,x_n] \ $ as desired. 
 
\smallskip 

The global height estimate follows then from the expression
$$
h(a,g_1,\dots,g_s)={1\over [K:\Q]}\sum_{v\in M_K^\infty}N_v\,
\max\{ h_v(a), h_v(g_1), \dots, h_v(g_s)\}. 
$$

\end{proof}

\begin{rem} \label{root of 1}
The fact that the bound \ref{uniforme} is uniform on $w$ for 
$w\, |\, v$ is the key that 
allows us to get rid of the roots of 1. 
This will  no longer be the case in our treatment of the more refined 
arithmetic Nullstellens\"atze in Chapter 4.
\end{rem}

\smallskip

The following example, although  not verifying the worst
case bound for the degrees,   improves the 
lower bound stated in the introduction for a general height estimate 
 and shows that the term $d^n h$ is unavoidable.

\begin{exmpl} \label{d^n h}
Set 
$$
f_1:= x_1- H, \ \ f_2:= x_2 - x_1^d, \dots,
\ \ f_n:= x_n- x_{n-1}^d ,  \ \ f_{n+1}:= x_n^d
$$
for any $d, \, H \in \N$.
These are polynomials 
without common zeros in $\A^n$
of degree and height bounded by $d$ and $h:= \log H$
respectively.

\smallskip

Let $a\in \Z\setminus \{0\}$ and 
$g_1, \dots , g_{n+1} \in \Z[x_1, \dots, x_n]$
such that 
$ \ a = g_1 \, f_1 + \cdots + g_{n+1} \, f_{n+1}$. 
We evaluate this identity 
in $(H, H^d, \cdots, H^{d^{n-1}})$ and we obtain
$$
a =g_{n+1}(H, H^d, \cdots, H^{d^{n-1}}) \, H^{d^{n}}
$$
from where we deduce $h(a) \ge d^n \, h$.
\end{exmpl}

%% file: seccion4.tex
\typeout{Seccion 4}

\section{Intrinsic type estimates}
\label{Intrinsic type estimates}

\setcounter{equation}{0}

Theorem 1 is essentially optimal in the 
general case. 
There are however many 
particular instances in which these estimates can be 
improved. 
Consider the following example: 
\begin{eqnarray*}
f_1:= x_1- 1, \ \ f_2:= x_2 - x_1^d, \dots,
\ \ f_n:= x_n- x_{n-1}^d ,  \ \ f_{n+1}:= H- x_n^d
\end{eqnarray*}
for any $d, H \in \N$.
These are polynomials 
without common zeros in $\A^n$
of degree and height bounded by $d$ and $h:= \log H$
respectively.
Theorem 1 says 
there exist $a\in \Z\setminus \{0\}$ and 
$g_1, \dots , g_{n+1} \in \Z[x_1, \dots, x_n]$
such that 
$$ 
a = g_1 \, f_1 + \cdots + g_{n+1} \, f_{n+1} 
$$
with $ \ \deg g_i  \le
4\,n\,d^n \ $ and 
$ \ h(a) , h(g_i)  \le
4 \, n \, (n+1) \, d^n\,( h+ (n+7)\,\log(n+1)\,d )$.
However the following 
B\'ezout identity holds: 
$$
H-1 = {x_1^d-1\over x_1 - 1} \cdots {x_{n}^d- 1 \over x_{n} -1} \, f_1
+ \cdots +
{x_{n}^d- 1 \over x_{n} -1}  \, f_n + f_{n+1}.
$$
Note that the polynomials arising in 
this identity have degree  and height bounded by 
$n \, (d-1)$ and $h$ respectively. 

There is in this case an exponential gap between 
the a priori general estimates
and the actual ones. 
The explanation is somewhat simple: 
for $i=1, \dots, n$, the varieties
$$
V_i:= V(f_1, \ldots, f_i) = V(x_1-1, x_2-1, \dots, x_i-1) \subset \A^n,
$$
verify $\deg(V_i)=1$ and $h(V_i)\le  2\,n\,\log(n+1)$.
Namely, both the degree and the height of the varieties
successively cut out by the input polynomials are much smaller
than the corresponding B\'ezout estimate. 

As the varieties $V_i$ verify the assumptions of
Lemma \ref{nullstlocal}, a direct application  together with
Lemma \ref{hacer entero} produces  the more realistic 
estimates:
$$
\deg g_i  \le  2\,n^2\,d \quad , \quad 
h(a) , h(g_i)  \le  (n+1)^2 (h+ 8\, n\, \log (n+1) \, d ) .
$$

\smallskip

Based on this idea, we devote this chapter to the study of 
more refined arithmetic Nullstellens\"atze
which can deal with such situations.


\subsection{Equations in general position}
\label{Equations in general position}

This section deals with the preparation of the input data.
To apply Lemma \ref{nullstlocal}, we need to prepare 
the polynomials and the 
variables of the ambient space. 

\smallskip

Let $f_1, \dots, f_s \in K[x_1, \dots, x_n]$
be polynomials without common zeros in $\A^n$.
For $i=1, \dots,  s$ and $a_{i\, j} \in \Z$ we set 
$$
q_i := a_{i\, 1} \, f_1 + \cdots + a_{i\, s } \, f_s
.
$$
We will estimate the height of rational integers $a_{ij}$
in order that there exists 
$t \le \min \{ n+1, s\}$ such that
$(q_1, \dots , q_i) \subset  K[x_1, \dots, x_n]$ is a radical ideal of 
dimension $n-i$ 
for $i=1, \dots, t-1$ and 
$  1\in (q_1, \dots, q_t)$.

Also we set 
$$
y_k:= b_{k \, 0} + b_{k \, 1} \, x_1 + \cdots + b_{k \, n} \, x_n
$$
for $k=1, \dots, n$ and $b_{k\, l} \in \Z$. 
Again we want to estimate the height of  rational
integers $b_{kl}$ such that 
$V_i:= V(q_1, \dots, q_i) \subset \A^n$ 
satisfies 
Assumption \ref{assumption} with respect to this set of variables
for $i=1, \dots, t-1$. 
Namely, the projection 
$$
\pi_i:V_i \to \A^{n-i} \quad \quad , \quad \quad 
x\mapsto (y_1, \dots, y_{n-i})
$$
must verify $\# \pi_i^{-1}(0) = \deg V_i$, that is
$ \ \# \, V_i \cap V(y_1,  \cdots , y_{n-i}) =\deg V_i \ $ 
for $i=1, \dots, t-1$.
Lemma \ref{evaluacion} implies that the variables $y_1, \dots, y_{n-i}$
are in Noether normal position with respect to $V_i$.

It is well-known that these conditions are satisfied by a
generic election of $a_{i j}$ and $b_{k l}$, see for instance
\cite{GiHeSa93}, \cite{SaSo95}.

\smallskip

We have already applied such a preparation to obtain
the  classic style
version of the effective arithmetic Nullstellensatz
presented in 
Theorem \ref{arithnullst}.  There, we chose
roots of 1 as coefficients of the linear combinations
since their existence was sufficient in our proof.
However, technical reasons (see Remark \ref{root of 1}) prevent
us to apply the same principle in this chapter, and we
need to 
carry out
a more careful analysis. 

\smallskip

We note that all  aspects of this preparation 
were previously covered in 
the research papers 
\cite{BeYg91}, \cite{GiHeSa93}, \cite{KrPa96},
\cite{HaMoPaSo98}. 
However the bounds presented therein are 
 either non-explicit or
 not precise enough for our purposes. 
Here we chose   to give a self-contained presentation,
which yields another  proof of the existence of such linear 
combinations, with good  control of the integers height.


\subsubsection{An effective Bertini theorem}

This subsection concerns  the preparation of the polynomials. 
We will first show  some auxiliary results.

\smallskip

The following is a version of the so-called shape lemma representation
of a 0-dimensional radical ideal.
The main difference here is that  we choose a generic linear form ---
instead of a particular one --- as a primitive element. 

\smallskip

For a polynomial 
$f=c_D \, t^D +\cdots + c_0 \in k[t]$ we denote its  discriminant
by $\discr (f) \in k$. 
We recall that $\discr \, (f) \not=0$ if and only if 
$c_D\not=0$ and $f$ is squarefree, that is when $f$ has exactly $D$ 
distinct roots. 

\begin{lem} \label{shape lemma}
{(Shape Lemma)}

Let $V\subset \A^n$ be a 0-dimensional
variety defined over $k$. 
Let $U:=(U_0,\dots,U_n)$ be a group of $n+1$ variables, and set
$
L:= U_0+ U_1 \, x_1 + \cdots + U_n \, x_n
$
for the associated generic linear form.

Let $P:= P_V \in k[U][T]$ be a characteristic polynomial of $V$. 
Set $P^\prime:= \partial P/\partial T \in k[U][T]$ and 
$\rho:=\discr_T P \in k[U]\setminus \{0\}$. 
Also set $I$ for the extension of $I(V)$ to $k[U][x]$. 

Then there exist $v_1, \dots, v_n \in k[U][T]$ with 
$\deg v_i \le \deg V-1$ such that
$$
I_\rho = (P(L) , P^\prime (L) \, x_1 - v_1(L),
\dots, P^\prime (L) \, x_n - v_n(L))_\rho  \ \subset k[U]_\rho [x].
$$
\end{lem}

\begin{proof}{Proof.--}
We note first that $I(V)$ is a radical ideal, and so 
$I= k[U] \otimes_k I(V)$ is also radical. 
We readily obtain from the definition of $P:= P_V$ 
that $I\cap k[U][L] = (P(L))$, and so 
$P(L) \in I$. 

\smallskip

We can write $P(L)=\sum_\alpha a_\alpha(x) \, U^\alpha$
with $a_\alpha(x) \in I(V)$.
Therefore  $ \partial P(L)/ \partial U_i $ also lies
in $I$ for all $i$.
A direct computation shows that for $i=1, \dots, n$ 
$$
\partial P(L) / \partial U_i = P^\prime(L) \, x_i - v_i(L)
$$
for some $v_i \in k[U][T]$ with
$\deg v_i \le \deg P-1 = \deg V-1$.

Set 
$$
 J:=(P(L) , P^\prime (L) \, x_1 - v_1(L),
\dots, P^\prime (L) \, x_n - v_n(L))  \subset k[U][x].
$$
The previous argument shows the inclusion
$
I\supset J$. 

\medskip

On the other hand,  $\rho = A \, P + B \, P^\prime$ 
for some $A, B \in k[U][T]$. 
Set $w_i := B \, v_i$. 
Then $x_i\equiv w_i(L)/\rho \pmod{J_\rho}$ and so for every $f\in k[U][x]$ 
we have that 
$ \ f\equiv f(U,w_1(L)/\rho,\dots,w_n(L)/\rho) \ $
modulo $J_\rho$, and hence modulo $I_\rho$.

For $f\in I$, 
$$
\rho^{\deg f} \, f(U,w_1(L)/\rho,\dots,w_n(L)/\rho) \in
I \cap k[U][L]= (P(L))
$$ 
which implies 
$I_\rho \subset J_\rho$ as desired.
\end{proof}

Let $\nu \in k^{n+1}$ such that
$\rho(\nu)\ne 0$. 
It follows  that $I(V)$ can be represented as
$$
I(V) = (P(L),  P^\prime (L) \,x_1 - v_1(L),
\dots, P^\prime (L) \, x_n - v_n(L))(\nu) \ \subset k[x].
$$

\smallskip

Now let $f_1, \dots, f_s \in k[x_1, \dots, x_n]$
be polynomials without common zeros in $\A^n$.
For $i = 1, \dots,  s$ we let 
$Z_i:=(Z_{i1},\dots,Z_{is})$ denote a group of $s$
variables, and we set
$$
Q_i:= Z_{i 1} \, f_1 + \cdots + Z_{i s} \, f_s \in k[Z][x]
$$
for the associated generic linear combination of $f_1, \dots, f_s$.

\begin{lem} \label{ideal primo}
For $\ell=1,\dots,s$, the ideal $(Q_1, \dots, Q_\ell)$
 is a complete intersection prime ideal of $k[Z][x]$.
\end{lem}

\begin{proof}{Proof.--}
Set $I:= (Q_1, \dots, Q_\ell) $ and $V:= V(I) \subset \A^{s \ell} \times \A^n$.
First we observe that $V$ is a linear bundle over $\A^n$:
the projection
$$
\pi : V\to \A^n \quad \quad , \quad \quad (z, x) \mapsto x
$$
is surjective, and the fibers are affine spaces of dimension $(s-1) \, \ell$.
This follows from the assumption that the $f_j$ have no common zeros.
This implies 
that
$$
\dim V =  (s-1) \, \ell + n
$$
because of the theorem of dimension of fibers.
Namely $Q_1, \dots, Q_\ell$ is a complete intersection, and 
in particular the ideal $I$ is unmixed.

\smallskip

Set
$
I= I_1 \cap \dots \cap I_m
$
for the primary decomposition of this ideal.
We will show that $I_j$ is prime for all $j$, and then that $m=1$.

\smallskip

First we have that 
$I_{f_j}= (Q_1/f_j , \dots, Q_\ell/f_j ) = (Z_{1j}+H_{1j},\dots, 
Z_{\ell j}+H_{ \ell j})$
 where $H_{ij} \in k[Z_i][x]_{f_j}$ does not
depend on $Z_{i j} $. Therefore
$$
(k[\A^{s \ell} \times \A^n] /I)_{f_j} \cong
k[\A^{(s-1) \ell} \times \A^n]_{f_j}
$$
is a domain, that is $I_{f_j}$ is prime.
We have $I_{f_j}= (I_1)_{f_j} \cap \dots \cap
(I_m)_{f_j}$, and so there exists
$1\le n(j)\le m$ such that
$$ 
 I_{f_j} = (I_{n(j)})_{f_j} \quad \quad , \quad \quad 
V(I_i)\subset \{ f_j=0\} \quad \mbox{ for } i  \neq n(j).
$$

In particular $ \ I_{n(j)}= I_{f_j} \cap k[\A^{s \ell} \times \A^n] \ $
is prime.
The fact that $\cap_j\{ f_j = 0 \}=\emptyset$
ensures that $n(j)$ runs over all $1\le i \le m$, and so $I$ is radical.

\smallskip

The expression   $I_{f_j}=
(Z_{1j}+H_{1j},\dots, Z_{\ell j}+H_{ \ell j})$
implies that 
$ \pi (V(I_{f_j})) \subset \A^n
$
contains the dense open set $\{f_j\not= 0\}$.
In particular $V(I_{f_j})$ is not contained in any of the hypersurfaces
$\{ f_i=0\}$ and so $n(j)=n(1)$ for all $j$. 
This  implies that $m=1$, and so $I=I_1$ is prime. 
\end{proof}

\smallskip
The following proposition shows that $(Q_1(a_1), 
\dots, Q_\ell(a_\ell))$ is a radical ideal for a generic election of 
$a_i:=(a_{i1},\dots,a_{is})$. 
Unlike  Lemmas \ref{shape lemma} and \ref{ideal primo}, 
this result does note hold for
 arbitrary characteristic. 
For instance, let  $ x^p,  1-x^p \in \F_p[x]$ for some 
prime $p$. 
Then $ \ Q_1(a_1) = b+ c \, x^p \ $ for some $b, c \in \overline{\F_p}$
and so 
$$
Q_1(a_1)= (b^{1/p}+c^{1/p}\, x)^p
$$ 
is not squarefree.

\begin{prop} \label{radical}
Let $\char (k)=0 $ and set $I:= (Q_1, \dots, Q_\ell) \subset k[Z][x]$. 
\begin{itemize}
\item In case $I\cap k[Z] \ne \{0\} $ there exists $F \in k[Z] \setminus \{0\}$
with $ \deg F \le (d+1)^{\ell}$ such that
$  F(a_1, \dots, a_\ell) \neq 0  $ for 
$a_1, \dots, a_\ell \in k^{s} $ implies that 
$\ 1\in (Q_1(a_1), \dots, Q_\ell (a_\ell)) $.  
\item In case $I\cap k[Z] = \{0\} $ there exists $F \in k[Z] \setminus \{0\}$
with $ \deg F \le 2\, (d+1)^{2\, \ell}$ such that
$  F(a_1, \dots, a_\ell) \neq 0  $ for 
$a_1, \dots, a_\ell \in k^{s} $ implies that 
$ \ (Q_1(a_1), \dots, Q_\ell(a_\ell)) \subset k[x]\ $ is a radical ideal
of dimension
$n-\ell$.
\end{itemize}
\end{prop}

\begin{proof}{Proof.--}
Set $V:=V(I)\subset \A^{s \ell}\times \A^n$.
We have $\dim V = (s-1)\, \ell +n$ and $\deg V\le (d+1)^\ell$.

\medskip

First we consider  the case $I \cap k[Z] \neq \{ 0\} $.
This occurs, for instance, when $\ell\ge n+1$, since then 
$ \ \dim
I=s\,\ell + n -\ell < \dim k[Z]=s\, \ell$.

Let $ \pi: \A^{s \ell} \times \A^n \to \A^{s \ell} $ be the canonical projection.
Then $\overline{\pi (V)}
$
is a proper subvariety 
of $\A^{ s \ell}$, and thus it is contained in a hypersurface of 
degree bounded by $\deg V$. 
This can be seen by 
taking a generic projection of
this variety 
into an affine space of dimension
$s\,\ell + n -\ell +1$ \cite[Remark 4]{Heintz83}.
Let $F \in k[Z]$ be a defining equation of this hypersurface. 
Then $F\in I$ as $I$ is prime, and we have $\deg F \le (d+1)^\ell$. 
Thus 
$$
1\in I_F  \subset k[Z]_F[x],
$$
and therefore $1\in I(a):=
(Q_1(a_1),\dots,Q_\ell(a_\ell)) $
for 
$a \in k^{s\ell} $ such that 
 $  F(a) \neq 0  $. 

\medskip
Next we consider the case   $\ I \,  \cap \, k[Z] = \{ 0\} $.

We adopt the following convention:
for an ideal  $J\subset k[x]$ and  for $\zeta$  any new
group of variables, we 
 denote by $J^{[\zeta]}$ and $J^{(\zeta)}$ 
the extension of $J$ to the 
polynomial rings $k[\zeta][x]$ and $k(\zeta)[x]$ respectively. 

\smallskip

We assume for the moment $\ell=n$. 
Then $\dim I= s\, \ell $ and so the extended ideal 
$I^{(Z)} \subset k(Z)[x]$ is a 0-dimensional 
prime ideal. 
We have then that $\overline{k(Z)} \otimes_k I^{(Z)}
\subset  \overline{k(Z)}[x]$ 
is a radical ideal, as 
$\char (k)= 0$ \cite[Thm. 26.3]{Matsumura86}. 

\smallskip

Our approach to this case is based on Shape Lemma
\ref{shape lemma}. 
We will determine a polynomial $F \in k[Z]$ 
such that $F(a) \neq 0$ implies that the shape lemma 
representation of $I^{(Z)}$ can be transferred to a 
shape lemma representation of $I(a)$. 

\smallskip

Let $U$ be a group of $n+1$ variables and set

$
\ L:= U_0 + U_1 \, x_1 + \cdots + U_n \, x_n
\ $
for the associated generic linear form.
Consider the morphism 
$$
\Psi : \A^{n+1}  \times \A^{s \ell} \times \A^n \to
 \A^{n+1} \times \A^{s \ell}  \times\A^1
\quad \quad , \quad \quad 
(u,z,x) \mapsto (u,z,L(x)).
$$
and let $W$ be the variety defined by $I$ 
in $ \A^{n+1}  \times \A^{s \ell} \times \A^n $, that is 
$W= \A^{n+1} \times V$. 
The Zariski closure $\overline{\Psi(W)}$ is then an irreducible
hypersurface. We
set $P \in k[U,Z][T]$ for one of its defining equations. 

If $I^{[U](Z)}$ is  the extension of $I^{(Z)}$ to $k[U](Z)[x]$,
the polynomial $P$ can be equivalently defined through the condition
that $P(L)$ is a generator of the principal ideal 
$I^{[U](Z)} \cap k[U,Z][L] $. 
Namely, $P$ is a
characteristic polynomial of 
the $0-$dimensional variety $W_0$ defined by $I^{(Z)}$ in 
$\A^n(\overline{k(Z)})$. 

\smallskip

Let 
$ v_1, \dots, v_n \in k[U](Z)[T]$ denote the polynomials
arising  in 
Shape Lemma applied to $W_0$. 
{}From the proof of this lemma we have that 
$$
\partial P(L) / \partial U_i = P^\prime(L) \, x_i - v_i(L)
\in k[U,Z][L]
$$
and so $v_i \in k[U,Z][T]$. 
Set 
$J:=(P(L), P^\prime(L) \, x_1 - v_1(L), \dots,
P^\prime(L) \, x_n - v_n(L)) \subset k[U,Z][x]$
and $\rho:= \discr_T P\in k[U,Z]\setminus \{0\}$. 
Then 
$$
(I^{[U](Z)})_\rho = (J^{[U](Z)})_\rho \ \subset k[U](Z)_\rho[x].
$$

We have that both $I^{[U,Z]}_\rho$ and $J^{[U,Z]}_\rho$ are prime
ideals of $k[U,Z]_\rho[x]$ with trivial intersection 
with the ring $k[U,Z]$. 
Thus they coincide with the contraction of 
$I^{[U](Z)}_\rho$ and $ J^{[U](Z)}_\rho $ to 
 $k[U,Z]_\rho[x]$ respectively, and so 
$$
I^{[U,Z]}_\rho = J^{[U,Z]}_\rho  \ \subset k[U,Z]_\rho[x].
$$

Define  $F\in k[Z]\setminus \{0\}$ as any of the non-zero coefficients of 
the monomial expansion of $\rho$ with respect to $U$.
Let $a\in k^{s\ell}$ such that $F(a)\not= 0 $. 
Then $\rho(U,a) \neq 0$ and so $P(U,a)[T]$ is squarefree. 
Then 
$$
(I(a)^{[U]})_{\rho(U,a)} =(P(L), P^\prime(L) \, x_1 - v_1(L), \dots,
P^\prime(L) \, x_n - v_n(L))(a) \  \subset k[U][x]
$$
is radical, which implies 
in turn that $I(a) = (I(a)^{[U]})_{\rho(U,a)} \cap k[x] $ 
is a radical ideal of $k[x]$ as desired. 

\smallskip

It remains
to estimate the degree of $F$.
To this end, it suffices to bound the degree of 
$\rho$ with respect to the group of variables $Z$. 
We recall that $P$ was defined as a defining equation of 
the hypersurface $\overline{\Psi(W)}$. 
The map $\Psi$ is linear in the variables $Z$ and $x$, 
and so 
$$
\deg_Z P \le \deg W = \deg V \le (d+1)^n. 
$$
This implies that $ \deg F \le 
\deg_Z \rho \le \deg_Z P\, (2\, \deg_Z P-1)  \le 2\, (d+1)^{2n}$.

\medskip

Finally we consider 
the case $\ell < n$ for  $I\cap k[Z]=\{0\}$.

Let $U_1, \dots, U_{n-\ell}$ be groups of $n+1$
variables each, and set
$$
L_i:=U_{i\, 0}+U_{i\, 1} \, x_1+\cdots
+ U_{i\, n}\, x_n.
$$ 
for $i=1, \dots, n-\ell$. 
Set $U:=(U_1, \dots, U_{n-\ell})$, 
$L:=(L_1, \dots, L_{n-\ell})$
 and $k_0:= k(U,L)$.
The extended ideal $I_0 \subset k_0[Z][x_1, \dots, x_\ell]$
verifies $I_0 \cap k_0[Z] =\{0\}$  and thus falls
into the previously considered case. 

\smallskip

Thus there exists $F_0 \in k_0[Z] \setminus \{0\}$
with $\deg F_0\le 2\, (d+1)^{2\ell}$ such that 
$F_0(a) \ne 0$ for $a\in k^{s\ell}$ implies 
that $I_0(a)$ is a radical ideal of 
$k_0[Z][x_1, \dots, x_\ell]$. 
This implies in turn that 
$I(a)$ is a radical ideal of $k[x]$, as 
$$
I(a)= I_0(a)\cap k[x].
$$
We can assume without loss of generality that $F_0$ lies
in $k[U,L][Z]$. 
We conclude by taking $F$ as any non-zero coefficient 
of the monomial expansion of $F_0$ 
with respect to the variables $U$ and $L$. 
\end{proof}

\begin{cor}
\label{bertini}
Let $\char (k)=0 $, and let
$f_1,\dots, f_s \in k[x_1,\dots,x_n]$ be polynomials
without common zeros in
$\A^n$.
Set $d:= \max_i \deg f_i$. 

Then there exist $ t\le \min\{ n+1, s\} $ and 
$a_1, \dots, a_t\in \Z^s$ such that 
\begin{itemize}
\item $(Q_1(a_1),\dots,Q_i(a_i))$ is a radical ideal
of dimension $n-i$ for $1\le i \le t-1$, 
\item $1 \in (Q_1(a_1),\dots,Q_t(a_t))$,
\item $h(a_i) \le 2\, (n+1)  \, \log(d+1) $.
\end{itemize}
\end{cor}

\begin{proof}{Proof.--}
Set $t$ for the minimal $i$ such that
$I_i:=(Q_1,\dots,Q_i) \cap k[Z]
\not= \{0\}$. Then $t\le n+1$,  and by the previous result 
there exists $F_t\in
k[Z]$ with $\deg F_t\le (d+1)^t$
 such that $F_t(a)\not=
0$ 
 implies that $1\in (Q_1(a_1),\dots,Q_t(a_t))$.

On the other hand, for $i < t$ we take a polynomial
$F_i\in
k[Z]$ of degree bounded by
$2\,(d+1)^{2\,i}$ such that $F_i(a)\not=
0$ implies that $(Q_1(a_1),\dots,Q_i(a_i))$ is a
radical ideal of dimension $n-i$.
Then we take $F:= F_1 \cdots F_t$ and so 
\begin{eqnarray*}
\deg F &\le &  2\,(d+1)^2 + \cdots + 2\,(d+1)^{2(t-1)} +
(d+1)^t\\
& \le & (d+1)^{2n} + 2\,(d+1)^{2n} + (d+1)^{n+1}\\
&\le & 4\,(d+1)^{2n}
.
\end{eqnarray*}
Finally,  $F\ne 0$ implies there 
exist $a_1, \dots, a_t \in \Z^s$ such that $h(a_i)\le 
\log (\deg F) $
and $F(a) \ne 0$. 

\end{proof}


\subsubsection{Effective Noether normal position}

Now we devote to the preparation of the variables.
For $k= 0, \dots,  n$  we let $U_k:=(U_{k0},\dots,U_{kn})$ 
be a group of $n+1$ variables
and we set
$$
Y_k := U_{k\, 0} + U_{k\, 1} \, x_1 + \cdots +U_{k \, n} \, x_n .
$$

\begin{prop} \label{noether}

Let $V\subset \A^n$ be an equidimensional
 variety of dimension $r$ defined over $k$.

Then there exists $G \in k[U_1,\dots,U_r] \setminus \{0\}$
with  $\deg_{U_k}G\le 2\,(\deg V)^2 $ such that
  $G(b_1, \dots, b_r) \neq 0$ for $b_1, \dots, b_r\in k^{n+1}$ 
  implies that
$$
\# \ V \cap V(Y_1(b_1), \dots , Y_r(b_r)) = \deg V.
$$
\end{prop}

\begin{proof}{Proof.--}

Let $f_V$ be a Chow form of $V$ and $P_V\in k[U,T]$ be
the
characteristic polynomial of $V$ associated to $f_V$
given by Lemma \ref{caracteristico}.

Set $D:= \deg V$ and let
$ \ P_V = c_D \, T_0^D + \cdots + c_0
\ $
be its expansion with respect to $T_0$. 
Also set 
$$
\rho := \discr_{T_0} P_V\in k[U_0, \dots, U_r][T_1, \dots, T_r]
\setminus \{0\}
$$
for the discriminant of $P_V$  with respect to
$ T_0$.

\smallskip

Observe that as $P_V$ is  multihomogeneous 
of degree $D$ in each group of variables 
$U_i \cup \{T_i\}$,
the degree of $\rho $ in 
each of these group of variables is bounded by 
$D\,(2D-1)$.

Now let $\nu_1, \dots, \nu_r \in \overline{k}^{n+1}$ such that
$V(\nu) := V\cap V(Y_1(\nu_1) ,  \dots , Y_r(\nu_r)) $
is a $0$-dimensional variety of cardinality $D$, and 
 $f_{V(\nu)}$ be a Chow forms of $V(\nu)$.

Set $\zeta_0:= ( T_0 - U_{0 0}, U_{0 1}. \dots, U_{0 n})$.
 Then applying
 \cite[Prop. 2.4]{Philippon86}, there exists 
$\lambda \in k^*$ such that:

$$ 
P_V(U_0, \nu_1, \dots, \nu_r)(T_0, 0, \dots, 0)= 
f_V(\zeta_0(U_0, T_0) , \nu_1, \dots, \nu_r)= 
 \lambda \, f_{V(\nu)}(\zeta_0(U_0, T_0)) = 
\lambda \, P_{V(\nu)}(U_0)(T_0)
$$
where $P_{V(\nu)}$ is a characteristic polynomial of $V(\nu)$.

This implies 
$P_V(U)(T_0, 0, \dots, 0) \in k[U][T_0]$ is a squarefree 
polynomial and so $\rho (U)(0) \neq 0$. 

\smallskip

We take $G \in k[U_1, \dots, U_r]$ as any non-zero coefficient
of the expansion of $\rho(U)(0)$ with respect to $U_0$. 
Therefore 
$$
\deg G \le \deg_{U_i}\rho(U)(0) \le D\, (2\, D-1).
$$ 

The condition $G(b)\not= 0$ implies that $\rho(U_0, b_1, \dots, b_r)(0)\ne 0$, 
and so $\# V(b) =D$. 
\end{proof}

As we noted before, this implies that the variables 
$Y_1(b_1), \dots, Y_r(b_r)$ are in Noether normal position with respect to 
the variety $V$. 

\begin{cor} \label{variables}
Let $\char(k)=0$ and let 
$q_1,\dots,q_t \in k[x_1,\dots,x_n]$ be polynomials
without common zeros in
$\A^n$ which form a reduced
weak 
regular sequence.
Set $d:= \max_i \deg f_i$. 

Then there exist
$b_1, \dots, b_n\in \Z^{n+1}$ such that for $i=1,\dots,t$,
$V(q_1,\dots,q_i)$ satisfies Assumption \ref{assumption}
with respect to the variables $Y_1(b_1), \dots, Y_{n-i}(b_{n-i})$
and
$$
h(b_k) \le 2\,( n+1)  \, \log(d+1). 
$$
\end{cor}

\begin{proof}{Proof.--} 
This follows readily from the previous result. 
We take $G_i $ as the polynomial corresponding to 
the variety $V(q_1, \dots, q_i)$ and we set 
$ \ G:= G_1 \cdots G_{t-1} \in k[U_1, \dots, U_n]$. 
We have $\deg_{U_j} G_i \le 2\, d^{2\, i}$ and so 
$$
\deg_{U_j} G \le 2 \, d^2 + \cdots + 2\, d^{2\, (t-1)} 
 \le 4\, d^{2\, (t-1)} \le 4\, d^{2\, n}.
$$
We conclude by taking $b_1, \dots, b_n \in \Z^{n+1}$ such that 
$h(b_i) \le \log (\deg G) $ and $G(b) \ne 0$. 
\end{proof}


\subsection{An intrinsic arithmetic Nullstellensatz}
\label{An intrinsic arithmetic Nullstellensatz}

In this section  we introduce the  notions of {\em degree 
and height of a polynomial system} defined over a
number field $K$.
Modulo setting the input equations in general position, 
these parameters measure 
the degree and height of the varieties successively cut out.

The resulting estimates for the arithmetic Nullstellensatz 
are {\em linear } in these parameters.

As an important particular case, we derive a sparse arithmetic 
Nullstellensatz. 


\subsubsection{Intrinsic parameters}
\label{Intrinsic parameters}

Let $f_1,\ldots,f_s \in K[x_1,\ldots,x_n]$
be polynomials of degree bounded by $d$ without
common zeros in  $\A^n$.
For $i = 1, \dots,  s$ we let 
$Z_i$ denote a group of $s$
variables and we set
$$
Q_i(Z):= Z_{i 1} \, f_1 + \cdots + Z_{i s} \, f_s \in K[Z][x]
$$
for the associated generic linear combination of $f_1, \dots, f_s$.

\smallskip

Let $\Gamma $ be the set of integer $s\times s-$matrices
$a=(a_{ij})_{ij} \in \Z^{s\times s}$ 
of height bounded by 
$\, 2\,(n+1) \, \log(d+1) $
such that 
$$
I_i(a):= (Q_1(a_1) , \dots, Q_i(a_i)) \subset K[x_1, \dots, x_n]
$$
is a radical ideal of dimension $n-i$ for $i=1, \dots, t-1$ and 
$1 \in I_t(a)$ for some $t \le \min\{ n+1, s\}$.

Corollary \ref{bertini} implies that $\Gamma\neq \emptyset$.

\smallskip

For  $a \in \Gamma$ we set
\begin{eqnarray*}
\delta(a)&:=& \max \ \{ \, \deg V(I_i(a))
\, ; \,  1\le i \le  \min\{t, n\}-1\, 
\},  \\
[1mm]
\eta(a)&:= & \max \ \{ \, h(V(I_i(a))
) \, ; \,  1\le i \le  t-1\, 
\,\}
.
\end{eqnarray*}

We set $\Gamma_{\min} \subset \Z^{s\times s} $ for the subset 
of matrices $ a\in \Gamma$ such that 
$ \  \eta(a) + d\, \delta(a) \ $ 
is minimum.
Finally let $a_{\, \min} \in \Gamma_{\min}$ be a matrix which attains
the minimum of $\delta (a) $ for $a\in \Gamma_{\min}$. 

\begin{defn}\label{defintr}
Let notations be as in the previous paragraph. 
Then we define the {\em degree} and the {\em height} of 
the polynomial system $f_1, \ldots, f_s $
respectively
as
$$
\delta(f_1, \ldots, f_s)
 :=  \delta(a_{\, \min})  \quad \quad , \quad \quad 
\eta(f_1, \ldots, f_s) := \eta(a_{\, \min}).
$$

\end{defn}

We restrict ourselves to integer matrices of bounded
height
in order to keep control of the height of
$Q_1(a_1), \dots, Q_t(a_t)$. 
The election of  $\eta(a)+ d\, \delta(a) $ 
as the defining invariant comes from the need
of estimating the degree and  height simultaneously. 

\smallskip

We note that in case $f_1, \ldots, f_s$
is already a reduced weak regular sequence
we have
$$
\eta(f_1,\dots,f_s)+d\, \delta(f_1,\dots,f_s)
\le \eta(\Id) + d\, \delta(\Id) .
$$

\smallskip

We can estimate this parameters through the arithmetic B\'ezout inequality:

\begin{lem} \label{d^n}
Let 
$f_1,\ldots,f_s \in K[x_1,\ldots,x_n]$
be polynomials without
common zeros in  $\A^n$.
Set $d_i:= \deg f_i $ and assume that 
$d_1 \ge \cdots \ge d_s$ holds. 
Set $d:= d_1 = \max_i \deg f_i$ and  $h:= h(f_1, \dots , f_s)$. 
Also set $n_0:= \min\{ n, s\}$
and 
$n_1:= \min\{ n+1, s\}$.
Then 
\begin{itemize}

\item $\delta(f_1, \dots, f_s)  \le \prod_{j=1}^{n_0 -1} d_j$, 

\item $\eta(f_1, \dots, f_s)  \le 
n\, \prod_{j=1}^{n_1-2} d_j \, (h + \log s + 3\, n \,  (n+1) \, d) $.

\end{itemize}

\end{lem}

\begin{proof}{Proof.--}
Let $ a:=a_{\, \min} = (a_{i j})_{i j} \in \Z^{s\times s}$ be 
the coefficient 
matrix which realizes the degree and height of 
$f_1, \dots, f_s$ and set 
$$
q_i:= a_{i 1} \, f_1 + \cdots + a_{i s} \, f_s \ ,
\qquad 1\le i \le s. 
$$

Let $t\le n_1 = \min\{ n+1, s\}$ 
be minimum such that $1 \in (q_1, \dots , q_t)$.
Let 
$\widetilde{a} \in \Z^{(t-1) \times s } $ be the matrix 
formed by the first 
$t-1$ rows of $a$ and let
$c \in \Z^{(t-1) \times s } $ be a staircase
matrix equivalent to $\widetilde{a} $.

The polynomial system 
$$
\widetilde{q}_i:= c_{i 1} \, f_1 + \cdots + c_{i s} \, f_s
$$
is then  equivalent to 
$q_1, \dots, q_{t-1}$, that is  
$(\widetilde{q}_1, \dots, \widetilde{q}_i)= 
({q}_1, \dots, {q}_i)$ for $i=1, \dots , t-1$. 
Also we have $\deg \widetilde{q}_i \le d_i$, and so 
$$
\delta := \max\, \{ \deg V_i ; \ 1\le i \le \min \{ n, t\} -1\} \le
\prod_{j=1}^{n_0-1} d_j .
$$

We have also that each coefficient of $c$ is a subdeterminant of 
$\widetilde{a}$. Thus 
\begin{eqnarray*}
\widetilde{h}:=h(\widetilde{q}_1,\dots,\widetilde{q}_{t-1})
 & \le & h + \log s + h (c) \\[1mm]
& \le & h + \log s + (t-1) \, (2\, (n+1) \, \log (d+1) + \log (t-1)) \\[1mm]
& \le & h + \log s + n \,  (3\, n+1) \, d
\end{eqnarray*}
and so, applying Corollary \ref{inters-global},  
\begin{eqnarray*}
\eta & \le & \max\, \{ h(V_i ) : 1\le i \le min \{ n+1, t\} -1\} \\
& \le & (\prod_{j=1}^{n_1-1} d_j) \, 
(\sum_{j=1}^{n_1-1} \widetilde{h})/d_j + (n+{n_1-1}) \, \log (n+1) ) \\[-2mm]
& \le & 
(\prod_{j=1}^{n_1-2} d_j)\,
(n\, (h + \log s + n \,  (3\, n+1) \, d) 
+ 2\, n \, \log (n+1) ) \\[-2mm]
& \le & 
n\, (\prod_{j=1}^{{n_1-2}} d_j )\, 
(h + \log s + 3\,n\,  (n+1)\,d) .
\end{eqnarray*}
\end{proof}

We can also estimate these parameters through the arithmetic 
Bernstein-Kushnirenko inequality: 

\begin{lem} \label{cota esparsa}
Let  
$f_1,\ldots,f_s \in K[x_1,\ldots,x_n]$
be polynomials without
common zeros in  $\A^n$.
Set $d:= \max_i \deg f_i $ and $h:= h(f_1, \dots , f_s)$.
Also let $\cV $ denote the volume of $1, x_1, \dots, x_n, 
f_1, \dots , f_s$.
Then 
\begin{itemize}

\item $\delta (f_1, \dots, f_s)\le \cV $, 

\item $\eta (f_1, \dots, f_s) 
\le n\, \cV \, (h+\log s + 2^{2\, n +4} \, d) $.

\end{itemize}
\end{lem}

\begin{proof}{Proof.--}

Let $a:=a_{\, \min } = (a_{i\, j})_{i\, j} \in \Z^{s\times s} $ and set
$ \ q_i := a_{i\, 1} \, f_1 + \cdots + a_{i\, s} \, f_s \ $ 
for $i=1, \dots, s$. 

Then $\Supp(q_i ) \subset \Supp(f_1, \dots , f_s) $ and so 
$ \ \cV(1, x_1, \dots x_n, q_1, \dots, q_s) \le \cV$. 

Applying  Proposition \ref{bernstein}
  we 
obtain $ \delta \le \cV $ and 
\begin{eqnarray*}
\eta & \le & (n\, \max_i h(q_i) + 2^{2\, n +3} \, \log (n+1) \, d \, ) \, 
\cV\\[1mm]
&\le & (n\, (h+ \log s + 2\, (n+1)\, \log (d+1) ) 
+ 2^{2\, n +3} \, \log (n+1) \, d \, ) \, \cV\\[2mm]
& \le & n  \, \cV \, (h+\log s + 2^{2\, n +4} \, d \, ).
\end{eqnarray*}
\end{proof}


\subsubsection{Proof of Theorem 2}
\label{Proof of Theorem 2}

Modulo the preparation of the input data, the proof of Theorem 2 follows 
the 
lines of the example introduced at the beginning of
Chapter 4. 

\begin{teo}{(Intrinsic arithmetic Nullstellensatz)} \label{intrinsic}

Let $K$ be a number field and let 
$f_1,\ldots,f_s \in \cO_K[x_1,\ldots,x_n]$
be polynomials without
common zeros in  $\A^n$.
Set $d:= \max_i \deg f_i $ and $h:= h(f_1, \dots , f_s)$.
Also let $\delta$ and $\eta$ denote the degree and the
height of
the polynomial system $f_1, \dots , f_s$.

Then there exist $a\in \Z\setminus \{0\}$ and
$g_1,\ldots,g_s \in {\cal O}_K  [x_1,\ldots,x_n]$
such that

\begin{itemize}
\item $a = {g}_1 \, {f}_1
+ \cdots + {g}_s \,
{f}_s $
,
\item
$ \deg g_i \le 2\, n^2 \, d\,\delta $
,
\item
$h(a, g_1, \dots, g_s) \leq   (n+1)^2 \, [K:\Q] \, d\, (
2\, \eta +  (h +\log s) \, \delta + 
21\,(n+1)^2\, d\, \log (d+1) \, \delta ).
$
\end{itemize}
\end{teo}

\begin{proof}{Proof.--}
Let $a_{\, \min}=  (a_{ij})_{ij} \in \Z^{ s \times s} $ 
be a coefficient matrix which realizes 
the degree and height of $f_1, \dots, f_s$. 
We have then $\delta=\delta(a_{\, \min } )$, 
$\eta= \eta(a_{\, \min } )$ and $\ h(a_{\, min}) \le 2\, (n+1)\, \log (d+1)$.
We set 
$$
q_i:= a_{i\, 1} \, f_1 + \cdots + a_{i\, s}\, f_s
$$
for $i=1, \dots, s$. 
Then $(q_1, \dots , q_i) $ is a radical ideal of dimension $n-i$ 
for $i=1, \dots, t-1$ and 
$1\in (q_1, \dots , q_t)$ for some 
$t\le \min\{ n+1, s\}$.

\smallskip

Also let $b_{kl} \in \Z$ be integers with $h(b_{k\, l} ) 
\le 2\, (n+1)\, \log (d+1)$ such that 
$V_i:= V(q_1, \dots, q_i) $ 
satisfies Assumption \ref{assumption} with respect to 
the variables
$$
y_k:= b_{k\, 0} + b_{k\, 1} \, x_1 + \cdots + b_{k\, n} \, x_n
$$
for $i=1, \dots , t-1$. 
Set $b:= (b_{k\, 0} )_k \in \Z^n$ and  $ B:= (b_{k\, l} )_{k, l \ge 1} 
\in \GL_n(\Q)$, and set 
$\varphi:\A^n \to \A^n$ for the affine map 
$ \ \varphi(x) := B\, x +b$. 
For $j=1, \dots, t$ we then set 
$$
F_j(y) := q_j(x) = q_j (\varphi^{-1}(y)) \ \in K[y_1,\dots,y_n]
$$
Thus $F_1, \dots, F_t$ are in the hypothesis of Lemma \ref{nullstlocal}
with respect to $y_1, \dots, y_n$ and we let 
$P_1, \dots, P_t\in K[x_1,\dots,x_n]$ be the polynomials 
satisfying B\'eszout identity we obtain there.

Finally, for $i=1, \dots, s$, we set 
$$
p_i := \sum_{j=1}^t a_{i\, j } P_j (\varphi(x)) \ \in K[x_1, \dots, x_n].
$$
 We have 
$ \ 1= p_1 \, f_1 + \cdots + p_s\, f_s$. 

\medskip

Now we analyze the degree and the height of these polynomials.
We will assume $n, d\ge 2$ as the remaining cases 
have already been considered in Lemmas \ref{d=1} and \ref{n=1}. 

\smallskip

Set $W_l:= V(F_1, \dots, F_l) \subset \A^n$ for $l=1, \dots, t-1$. 
We have $W_l= \varphi(V_l)$ and so 
$ \deg W_l= \deg V_l$.
We have also $\deg F_j =\deg q_j \le d$ and so 
$$
\deg p_i \le \max_j \deg P_j \le 2 \, n\, d\, 
(1 + \sum_{l=1}^{\min\{n,s\}-1} 
\deg W_l) \le 2\, n^2 \, d \, \delta
$$
as $\deg W_l \le \delta $ for $l \le n-1$. 

Now let $v\in M_K^\infty$ and set $h_v := \max_i h_v(f_i)$. 
We have 
\begin{eqnarray*}
h_\infty (\varphi^{-1}) & \le & n \, (h_\infty (\varphi) + \log n) 
 - \log |\det B|_\infty \\[1mm]
& \le & n \, ( 2\, (n+1)\, \log (d+1)  + \log n) - \log |\det B|_\infty \\[1mm]
& \le &  3\, n \, (n+1)\, \log (d+1)  - \log |\det B|_\infty 
\end{eqnarray*}

Then 
\begin{eqnarray*}
h_v(F_i) &\le &  h_v(q_i)  + (h_\infty(\varphi^{-1})
+2\, \log (n+1) ) \,\deg q_i\\
[1mm]
&\le &  h_v+ 2\, (n+1) \, \log(d+1) + \log s
+ (3\,n \, ( n+1) \, \log (d+1) - \log |\det B|_\infty 
+2\, \log (n+1) ) \, d \\
[1mm]
&\le &  h_v+ \log s + (n+1 + 3\,n \, ( n+1) + 2\, n) \, d\, \log (d+1) 
- \log |\det B|_\infty \, d \\
[1mm]
&\le& h_v + \log s + 3\,( n+1)^2 \, d\, \log (d+1)
-
\log|\det B|_\infty
\, d.
\end{eqnarray*}
by Lemma \ref{hprod}(c) and the facts that $\log
(n+1)\le n$ and $\log (d+1)\ge 1$ for $d\ge 2$. 
Next, applying Lemma \ref{inversible}, 
we obtain 
\begin{eqnarray*}
h(W_l) &\le &
h(V_l) + (n-l+1)(h(\varphi) + 5\,\log(n+1))\,\deg V_l\\
[2mm]
&\le & 
h(V_l) + n\, (2\, (n+1) \, \log (d+1) +\, 5\,\log(n+1)) \, \deg V_l \\[2mm]
&\le & 
\eta + n\, (7\,n+2) \, d\, \log (d+1) \, \delta 
\end{eqnarray*}
as $\deg W_l = \deg V_l \le d\, \delta$ and $h(V_l) \le \eta$
for $l=1, \dots, t-1$. 
By Lemma \ref{nullstlocal} there exists $\xi \in K^*$ such that 
\begin{eqnarray*}
h_v(P_j )&\le&  2\, n\, d\,
\sum_{l=1}^{t-1} h_v(W_l) + ((n+1)\,\max_l h_v(F_l)
 + \ 2\, n \, (2\, n+5) \, \log(n+1) \, d ) \,(1+ \sum_{l=1}^{t-1}
\deg W_l )  \\
& & - \ \log|\xi|_v\\
&\le& 
2\, n\, d\,
\sum_{l=1}^{t-1} h_v(W_l) + (n+1)^2\,
(h_v + \log s) \, d\, \delta  \\[1mm]
& & 
+ (3\,( n+1)^4 + \ 2\, n^2 \, (2\, n+5) \, (n+1) )\, d^2\, \log (d+1) \, 
\delta  
- \log |\mu|_v
\end{eqnarray*}
with $\mu:= (\det B)^{(n+1)^2\,d^2 \,\delta} \,  \xi\in K^*$.
Then 
\begin{eqnarray*}
h_v(p_i) & \le& \max_j h_v (P_j) + ( h_\infty(\varphi)+ 2\,\log(n+1)) 
\max_j \deg
P_j + 2\, (n+1) \, \log (d+1) + \log t \\
[2mm]
&\le & 2\,n\,d \,\sum_l h_v(W_l) + (n+1)^2\,(h_v+\log s) \,d\, \delta \\[-1mm]
&& + (3\,( n+1)^4 + \ 2\, n^2 \, (2\, n+5) \, (n+1) ) 
\,d^2\,\log (d+1) \, \delta
- \log |\mu|_v \\ [2mm]
& & + \  (2\, (n+1) \, \log (d+1) + 2\, \log (n+1) ) \, 2 \, n^2 \, d\, \delta 
+ 2\, (n+1) \, \log (d+1) + \log (n+1) \\[3mm]
&\le & 2\,n\,d \, \sum_l h_v(W_l) + (n+1)^2\,(h_v+\log s) \,d\, \delta
+ 7 \, (n+1)^3 \, (n+2)  \,d^2\,\log (d+1) \, \delta
- \log |\mu|_v . 
\end{eqnarray*}

Analogously 
$ \ h_v(p_i)  \le
2\,n\,d \, \sum_l h_v(W_l) + (n+1)^2\,h_v \,  d\, \delta
- \log |\mu|_v  \ $ for $v \notin M_K^\infty$. 

Hence
\begin{eqnarray*}
h(p_1,\dots,p_s) & \le & 
2\,n\,d \, \sum_l h(W_l) + (n+1)^2\,(h +\log s) \,d\, \delta
+ 7 \, (n+1)^3 \, (n+2)  \,d^2\,\log (d+1) \, \delta
\\
& \le & 
2\,n^2 \,d \, \eta + 2\, n^3 \, (7\, n+2) \, \,d^2\,\log (d+1) \, 
\delta+  (n+1)^2\,(h_v+\log s) \,d\, \delta \\[1mm]
&& + \ 7 \, (n+1)^3 \, (n+2)  \,d^2\,\log (d+1) \, \delta
\\[2mm]
& \le & 
2\,n^2\,d \, \eta + (n+1)^2\,(h_v+\log s) \, d\, \delta
+ 21 \, (n+1)^4 \,d^2\,\log (d+1) \, \delta
 .
\end{eqnarray*}

\smallskip

Finally 
we apply Lemma \ref{hacer entero} to obtain $a\in \Z\setminus \{0\}$ 
such that $g_i:= a\, p_i \in \cO_K[x_1,\dots,x_n]$. Thus 
$$
a=g_1\,f_1+\cdots+g_s\,f_s
$$
and the corresponding height estimates are 
 multiplied by  $[K:\Q]$.
\end{proof}

We derive from this result and the B\'ezout inequality
\ref{d^n}
the following estimate in terms of the degree and the height of the input 
polynomials:

\begin{cor}
Let 
$f_1,\ldots,f_s \in \cO_K[x_1,\ldots,x_n]$
be polynomials without
common zeros in  $\A^n$.
Set $d_i:= \deg f_i $ and assume that 
$d_1 \ge \cdots \ge d_s$ holds. 
Also set $d:= d_1= \max_i \deg f_i $, $h:= h(f_1, \dots , f_s)$, and 
$n_0:= \min\{ n,s\}$.

Then there exist $a\in \Z\setminus \{0\}$ and
$g_1,\ldots,g_s \in {\cal O}_K  [x_1,\ldots,x_n]$
such that

\begin{itemize}
\item $a = {g}_1 \, {f}_1
+ \cdots + {g}_s \,
{f}_s $
,
\item
$ \deg g_i \le 2\, n^2 \, d \, \prod_{j=1}^{n_0-1} d_j $
,
\item
$h(a, g_1, \dots, g_s) \leq   2\, (n+1)^3 \,  [K:\Q] 
\, d\, 
\prod_{j=1}^{n_0-1} d_j 
\, \, ( h +\log s  + 3\, n  (n+7) \, d\, \log (d+1))$.
\end{itemize}
\end{cor}


\subsubsection{Estimates for the sparse case}
\label{Estimates for the sparse case}

Bernstein-Kushnirenko inequality \ref{cota esparsa}
shows that both the degree and the height of a system are 
controlled by its volume. 
We derive then from Theorem \ref{intrinsic} an 
arithmetic Nullstellensatz for sparse polynomial systems.

\begin{cor}{(Sparse arithmetic Nullstellensatz)} 
\label{Sparse arithmetic Nullstellensatz}

Let  
$f_1,\ldots,f_s \in \cO_K[x_1,\ldots,x_n]$
be polynomials without
common zeros in  $\A^n$.
Set $d:= \max_i \deg f_i $ and $h:= h(f_1, \dots , f_s)$.
Also let $\cV $ denote the volume of the
polynomial system $1, x_1, \dots, x_n, 
f_1, \dots , f_s$.

Then there exist $a\in \Z\setminus \{0\}$ and
$g_1,\ldots,g_s \in {\cal O}_K  [x_1,\ldots,x_n]$
such that

\begin{itemize}
\item $a = {g}_1 \, {f}_1
+ \cdots + {g}_s \,
{f}_s $
\item
$ \deg g_i \le 2\, n^2 \, d\,\,\cV$
,
\item
$h(a, g_1, \dots, g_s) \leq  2\, (n+1)^3\,  [K:\Q]  \,d\, \, \cV\,(\,h
+\log s    +  2^{2n+4} d\,  \log(d+1) )
.
$
\end{itemize}
\end{cor}

\begin{exmpl} \label{ejemplosparse}
For $1\le i \le s$ we  let
$$f_i:=a_{i\, 0}+a_{i\, 1}x_1+\cdots + a_{i\, n}x_n + b_{i\, 1}x_1\cdots x_n
+\cdots + b_{i\, d} (x_1 \cdots x_n)^d \in \Z[x_1, \ldots, x_n]
$$
be polynomials of degree bounded by $n\, d$ 
without common zeros in $\A^n$. 
Set $h:= \max_i h(f_i)$. 
Also set 
$\ \cP_d := \Conv (0, e_1, \ldots, e_n, d\, (e_1+\cdots+e_n)) \subset \R^n$,
so that $\cP_d$  contains the
Newton  polytope  of  the polynomials
$1, x_1, \ldots, x_n, f_1, \ldots, f_s$.
Then 
$$
\cV\le \Vol (\cP_d) = n! \,d/(n-1)!\,=\, n\,d . 
$$
We conclude that there exist $a\in \Z\setminus \{ 0\}$
and $g_1, \ldots, g_s \in \Z[x_1, \ldots, x_n]$ such that 

\begin{itemize}

\item $a=g_1
\, f_1 + \cdots + g_s \, f_s $, 

\item $  \deg g_i  \le 2\, n^4\,
d^2 $ , 

\item $ h(a) , h(g_i) \leq 2 \, n^2 \, (n+1)^3  \, d^2
\,  (\,h+ \log s + n\, 2^{2\, n+4} \, d \, 
\log(n\, d+1)) $ . 

\end{itemize} 

This estimate is sharper than the one 
given by
Theorem 1.
\end{exmpl}

%% file: ref.tex
%
%

\typeout{Referencias}

%% file: direc.tex
\noindent {\sc Teresa Krick: } 
 Departamento de Matem\'atica, 
Universidad de Buenos Aires, 
 Ciudad Universitaria, 
1428 Buenos Aires, Argentina  \\ 
{\tt E-mail : krick@dm.uba.ar}  

\vspace{4mm}

\noindent {\sc Luis Miguel Pardo: }
Departamento de Matem\'aticas, Estad\'\i stica y Computaci\'on, 
Universidad de Cantabria, 
E-39071 Santander, Espa\~na \\ 
{\tt E-mail: pardo@matesco.unican.es} 

\vspace{4mm}

\noindent {\sc Mart\'\i n Sombra: } 
Departamento de Matem\'atica, 
Universidad Nacional de La Plata,
Calle 50 y 115, 
1900 La Plata, Argentina, and 
School of Mathematics,
Institute for Advanced Study,
Princeton NJ 08540, USA, \\ 
{\tt E-mail: sombra@mate.unlp.edu.ar, sombra@ias.edu}